\documentstyle[12pt]{article}
\def\sq{\hbox {\rlap{$\sqcap$}$\sqcup$}}
\def\sq{\hbox {\rlap{$\sqcap$}$\sqcup$}}
\def\R{ {\rm R \kern -.31cm I \kern .15cm}}
\def\C{ {\rm C \kern -.15cm \vrule width.5pt \kern .12cm}}
\def\Z{ {\rm Z \kern -.27cm \angle \kern .02cm}}
\def\N{ {\rm N \kern -.26cm \vrule width.4pt \kern .10cm}}
\def\1{{\rm 1\mskip-4.5mu l} }
\def\lsim{\raise0.3ex\hbox{$<$\kern-0.75em\raise-1.1ex\hbox{$\sim$}}}
\def\gsim{\raise0.3ex\hbox{$>$\kern-0.75em\raise-1.1ex\hbox{$\sim$}}}
\def\noi{\noindent}

\def\beq{\begin{equation}}   \def\eeq{\end{equation}}
\def\bea{\begin{eqnarray}}  \def\eea{\end{eqnarray}}
\def\nn{\nonumber}
\def\noi{\noindent}
\def\beeq{\begin{eqnarray}} \def\eeeq{\end{eqnarray}}
\newcommand\mysection{\setcounter{equation}{0}\section}

\newcounter{hran}

\begin{document}
\centerline{\Large\bf Long Range Scattering and Modified }
  \vskip 3 truemm \centerline{\Large\bf  Wave Operators for the 
Maxwell-Schr\"odinger
System} \vskip 3 truemm
\centerline{\Large\bf I. The case of vanishing asymptotic magnetic field}
\vskip 0.5
truecm

\centerline{\bf J. Ginibre}
\centerline{Laboratoire de Physique Th\'eorique\footnote{Unit\'e Mixte de
Recherche (CNRS) UMR 8627}}  \centerline{Universit\'e de Paris XI, B\^atiment
210, F-91405 ORSAY Cedex, France}
\vskip 3 truemm
\centerline{\bf G. Velo}
\centerline{Dipartimento di Fisica, Universit\`a di Bologna} 
\centerline{and INFN, Sezione di
Bologna, Italy}

\vskip 1 truecm

\begin{abstract}
We study the theory of scattering for the Maxwell-Schr\"odinger 
system in space dimension 3,
in the Coulomb gauge. In the special case of vanishing asymptotic 
magnetic field, we prove
the existence of modified wave operators for that system with no size 
restriction on the
Schr\"odinger data and we determine the asymptotic behaviour in time 
of solutions in the
range of the wave operators. The method consists in partially solving 
the Maxwell equations
for the potentials, substituting the result into the Schr\"odinger 
equation, which then
becomes both nonlinear and nonlocal in time, and treating the latter 
by the method
previously used for the Hartree equation and for the 
Wave-Schr\"odinger system.   \end{abstract}

\vskip 3 truecm
\noi AMS Classification : Primary 35P25. Secondary 35B40, 35Q40, 81U99.  \par
\noi Key words : Long range scattering, modified wave operators, 
Maxwell-Schr\"odinger system.
\vskip 1 truecm

\noindent LPT Orsay 02-53 \par
\noindent June 2002 \par

\newpage
\pagestyle{plain}
\baselineskip 18pt

\mysection{Introduction}
\hspace*{\parindent}
This paper is devoted to the theory of scattering and more precisely 
to the existence of
modified wave operators for the Maxwell-Schr\"odinger (MS) system in 
$3 + 1$ dimensional
space time. This system describes the evolution of a charged 
nonrelativistic quantum
mechanical particle interacting with the (classical) electromagnetic 
field it generates.\par

We write that system in Lorentz covariant notation~: greek indices 
run from 0 to 3, latin
indices run from 1 to 3, indices are raised and lowered with the 
metric tensor $g_{\mu\nu} =
(1,-1,-1,-1)$ so that, for any vector field $v = (v^{\mu})$, $v_0 = 
v^0$ and $v_j = - v^j$, and
we use the standard summation convention on repeated indices. The MS 
Lagrangian density can be
written as
\beq \label{1.1e}
{\cal L} = - (1/2) \ F_{\mu\nu} \ F^{\mu\nu} - \ {\rm Im} \ \bar{u} \ 
D_0u + (1/2)
\overline{D_ju}\ D^ju \eeq

\noi where
$$F_{\mu\nu} = \partial_{\mu} A_{\nu} - \partial_{\nu} A_{\mu} \qquad 
, \qquad D_{\mu} =
\partial_{\mu} + i \ A_{\mu}$$

\noi and $u$ and $A_{\mu}$ are respectively a complex scalar valued 
and a real vector valued
function defined in ${I\hskip-1truemm R}^{3+1}$. Here $\{F_{0j}\}$ 
are the components of the electric field and $\{F_{ij}\}$ are the
components of the magnetic field. The associated variational system is  \beq
\label{1.2e}
\left \{ \begin{array}{l} i\partial_0 u = (1/2) \ D_j D^ju +
A_0u \\ \\ \partial^{\nu} \ F_{\nu\mu} = J_{\mu}\end{array} \right . \eeq

\noi where the current density $J_{\mu}$ is given by
\beq
\label{1.3e}
J_0 = |u|^2 \qquad , \qquad J_j = - \ {\rm Im}\ \bar{u} \ D_ju
\eeq

\noi and satisfies the local conservation condition $\partial_{\mu} 
J^{\mu} = 0$. Formally the
$L^2$-norm is conserved as well as the energy
  \beq
\label{1.4e}
E(u,A) = \int dx \left \{ (1 /2) \left ( \sum_{\mu < \nu} 
F_{\mu\nu}^2 + \sum_j |D_ju|^2 \right ) +
A_0 |u|^2 \right \} \ . \eeq

\noi The system (\ref{1.2e}) is gauge invariant and we shall study it 
in the Coulomb gauge
which experience shows to be the most convenient one for the purpose 
of analysis. \par

The MS system (\ref{1.2e}) in ${I\hskip-1truemm R}^{3+1}$ is known to 
be locally well posed in
sufficiently regular spaces in the Lorentz gauge $\partial_{\mu} 
A^{\mu} = 0$ \cite{14r} and to
have global weak solutions in the energy space in various gauges 
including the Lorentz and
Coulomb gauges \cite{7r}. However that system is so far not known to 
be globally well posed in
any space whatever the gauge is. \par

A large amount of work has been devoted to the theory of scattering 
for nonlinear equations
and systems centering on the Schr\"odinger equation, in particular 
for nonlinear Schr\"odinger
(NLS) equations, Hartree equations, Klein-Gordon-Schr\"odinger (KGS) 
systems, Wave-Schr\"odinger
(WS) systems and Maxwell-Schr\"odinger (MS) systems. As in the case 
of the linear Schr\"odinger
equation, one must distinguish the short range case from the long 
range case. In the former
case, ordinary wave operators are expected and in a number of cases 
proved to exist, describing
solutions where the Schr\"odinger function behaves asymptotically 
like a solution of the free
Schr\"odinger equation. In the latter case, ordinary wave operators 
do not exist and have to be
replaced by modified wave operators including a suitable phase in 
their definition. In that
respect, the MS system (\ref{1.2e}) in ${I\hskip-1truemm R}^{3+1}$ 
belongs to the borderline
(Coulomb) long range case, because of the $t^{-1}$ decay in 
$L^{\infty}$ norm of solutions of
the wave equation. Such is the case also for the Hartree equation 
with $|x|^{-1}$ potential,
for the WS system in ${I\hskip-1truemm R}^{3+1}$ and for the KGS 
system in ${I\hskip-1truemm
R}^{2+1}$. \par

Whereas a well developed theory of long range scattering exists for 
the linear Schr\"odinger
equation (see \cite{1r} for a recent treatment and for an extensive 
bibliography), there exist
only few results on nonlinear long range scattering. The existence of 
modified wave operators
in the borderline Coulomb case has been proved first for the NLS 
equation in space dimension $n
= 1$ \cite{17r}, then for the NLS equation in dimensions $n = 2,3$ 
and for the Hartree equation
in dimension $n \geq 2$ \cite{2r}, for the derivative NLS equation in 
dimension $n = 1$
\cite{10r}, for the KGS system in dimension 2 \cite{18r} and for the 
MS system in dimension 3
\cite{19r}. All those results are restricted to the case of small 
data. In the case of
arbitrarily large data, the existence of modified wave operators has 
been proved for a family
of Hartree type equations with general (not only Coulomb) long range 
interactions \cite{3r}
\cite{4r} \cite{5r} by a method inspired by a previous series of 
papers by Hayashi et al
\cite{8r} \cite{9r} on the Hartree equation. Part of the results have 
been improved as
regards regularity \cite{15r} \cite{16r}. Finally the existence of 
modified wave operators
without any size restriction on the data has been proved for the WS 
system in dimension $n =
3$ \cite{6r}, by an extension of the method of \cite{3r} \cite{4r}. \par

The present paper is devoted to the extension of the results of 
\cite{6r} to the MS system in
the Coulomb gauge, and in particular to the proof of the existence of 
modified wave operators
for that system without any size restriction on the Schr\"odinger 
data, in the special case
of vanishing asymptotic magnetic field, namely in the special case 
where the asymptotic
state for the vector potential is zero. The method consists in 
replacing the Maxwell
equation for the vector potential by the associated integral equation 
and substituting the
latter into the Schr\"odinger equation, thereby obtaining a new 
Schr\"odinger equation which
is both nonlinear and nonlocal in time. The latter is then treated as 
in \cite{6r}, namely
$u$ is expressed in terms of an amplitude $w$ and a phase $\varphi$
satisfying an auxiliary system similar to that introduced in 
\cite{6r}. Wave operators are
constructed first for that auxiliary system, and then used to 
construct modified wave
operators for the original system (\ref{1.2e}). The detailed 
construction is too complicated to
allow for a more precise description at this stage, and will be 
described in heuristic terms in
Section 2 below. \par

The results of this paper improve over those of \cite{19r} by the 
fact that we do not
require smallness conditions on the Schr\"odinger asymptotic data. On 
the other hand, in
contrast with \cite{19r}, we restrict our attention to the case of 
vanishing asymptotic
data for the vector potential in order to avoid the difficulties 
coming from the
difference of propagation properties of the Wave and Schr\"odinger 
equations that occur
both in \cite{6r} and in \cite{19r}. In a subsequent paper we hope to 
extend the results
of the present one to the general case. \par

We now give a brief outline of the contents of this paper. A more 
detailed description of the
technical parts will be given at the end of Section 2. After 
collecting some notation and
preliminary estimates in Section 3, we study the asymptotic dynamics 
for the auxiliary system
in Section 4 and we construct the wave operators for that system by 
solving the local Cauchy
problem at infinity associated with it in Section 5, which contains 
the main technical results
of this paper.  We finally come back from the auxiliary system to the 
original one and
construct the modified wave operators for the latter in Section 6, 
where the final result is
stated in Proposition 6.1. \par
 
We conclude this section with some general notation which will be 
used freely throughout this
paper. We denote by $\parallel \cdot \parallel_r$ the norm in $L^r 
\equiv L^r({I\hskip-1truemm
R}^3)$ and we define $\delta (r) = 3/2 - 3/r$. For any interval $I$ 
and any Banach space $X$,
we denote by ${\cal C}(I,X)$ the space of strongly continuous 
functions from $I$ to $X$ and by
$L^{\infty}(I,X)$ (resp. $L^{\infty}_{loc}(I,X))$ the space of 
measurable essentially bounded
(resp. locally essentially bounded) functions from $I$ to $X$. For 
real numbers $a$ and $b$, we
use the notation $a \vee b = {\rm Max}(a,b)$, and $a \wedge b = {\rm 
Min} (a, b)$. In the
estimates of solutions of the relevant equations, we shall use the 
letter $C$ to denote
constants, possibly different from an estimate to the next, depending 
on various parameters,
but not on the solutions themselves or on their initial data. We 
shall use the notation
$C(a_1,a_2, \cdots )$ for estimating functions, also possibly 
different from an estimate to the
next, depending in addition on suitable norms $a_1, a_2, \cdots $ of 
the solutions or of their
initial data. Additional notation will be given in Section 3.

\mysection{Heuristics and formal computations} \hspace*{\parindent} 
In this section, we discuss
in heuristic terms the construction of the modifed wave operators for 
the MS system as it will
be performed in this paper and we derive the equations needed for 
that purpose. \par

We first recast the MS system in the Coulomb gauge $\partial_j A^j = 
0$ in the form that will
be used later on. We introduce the noncovariant notation
$$\partial_t = \partial_0 \quad , \quad \nabla = \{ \partial_j\} 
\quad , \quad A = \{ A^j\}
\quad , \quad J = \{ J^j\}$$

\noi so that
$$\{D_j\} = \nabla - i A \quad , \quad J \equiv J(u, A) = \ {\rm Im} 
\ \bar{u} (\nabla - iA)u$$

\noi and the Coulomb gauge condition becomes $\nabla \cdot A = 0$. 
With that notation, the MS
system in the Coulomb gauge becomes
\bea \label{2.1e}
&&i\partial_t u = - (1/2) (\nabla - i A)^2 u + A_0 u\\
\label{2.2e}
&&\Delta A_0 = - J_0 \\
&&\sq A + \nabla \left (\partial_t A_0 \right ) = J
\label{2.3e}
\eea

\noi where $\sq = \partial_t^2 - \Delta$. We replace that system by a 
formally equivalent one
in the following standard way. We solve (\ref{2.2e}) for $A_0$ as
\beq
\label{2.4e}
A_0 = - \Delta^{-1} \ J_0 = (4 \pi |x|)^{-1} \ * \  |u|^2 \equiv g(u)
\eeq

\noi where $*$ denotes the convolution in ${I\hskip-1truemm R}^3$, so 
that by the current
conservation $\partial_t J_0 + \nabla \cdot J = 0$,
\beq
\label{2.5e}
\partial_t \ A_0 = \Delta^{-1} \nabla \cdot J \ .
\eeq

\noi Substituting (\ref{2.4e}) into (\ref{2.1e}) and (\ref{2.5e}) 
into (\ref{2.3e}), we obtain
the new system
  \bea
\label{2.6e}
&&i \partial_t u = - (1/2) (\nabla - iA)^2u + g(u) u \\
\label{2.7e}
&&\sq A = PJ \equiv P\ {\rm Im}\ \bar{u}(\nabla - iA)u
\eea

\noi where $P = \1 - \nabla \Delta^{-1} \nabla$ is the projector on 
divergence free vector
fields. The system (\ref{2.6e}) (\ref{2.7e}) is the starting point of 
our investigation. We
want to address the problem of classifying the asymptotic behaviours 
in time of the solutions
of the system (\ref{2.6e}) (\ref{2.7e}) by relating them to a set of 
model functions ${\cal
V} = \{ v = v(v_+)\}$ parametrized by some data $v_+$ and with 
suitably chosen and preferably
simple asymptotic behaviour in time. For each $v \in {\cal V}$, one 
tries to construct a
solution $(u, A)$ of the system (\ref{2.6e}) (\ref{2.7e}) defined at 
least in a neighborhood
of infinity in time and such that $(u, A)(t)$ behaves as $v(t)$ when 
$t \to \infty$ in a suitable sense. We then
define the wave operator as the map $\Omega : v_+ \to (u, A)$ thereby 
obtained. A similar
question can be asked for $t \to - \infty$. We restrict our attention 
to positive time. The
more standard definition of the wave operator is to define it as the 
map $v_+ \to (u,A)(0)$,
but what really matters is the solution $(u, A)$ in the neighborhood 
of infinity in time,
namely in some interval $[T, \infty )$. Continuing such a solution 
down to $t = 0$ is
a somewhat different question which we shall not touch here, 
especially since the MS system is not known to be globally
well posed in any reasonable space. \par

In the case of the MS system, which is long range, it is known that 
one cannot take for ${\cal V}$ the set of solutions
of the linear problem underlying (\ref{2.6e}) (\ref{2.7e}), namely
\beq
\label{2.8e}
\left \{ \begin{array}{l} i\partial_t u = - (1/2) \ \Delta u\\ \\ \sq 
A = 0\end{array} \right . \eeq

\noi and one of the tasks that will be performed in this paper will 
be to construct a better set ${\cal V}$ of model
asymptotic functions. The same situation prevails for long range 
Hartree equations and for the WS system and we refer
to Sections 2 of \cite{3r} \cite{4r} \cite{6r} for more details on 
that point.\par

Constructing the wave operators essentially amounts to solving the 
Cauchy problem with infinite initial time. The
system (\ref{2.6e}) (\ref{2.7e}) in this form is not well suited for 
that purpose, and we now perform a number of
transformations leading to an auxiliary system for which that problem 
can be handled. We first replace the
equation (\ref{2.7e}) by the associated integral equation namely
\beq
\label{2.9e}
A = A_0^{\infty} + A_1(u, A)
\eeq

\noi where
\beq
\label{2.10e}
A_0^{\infty} = \dot{K}(t) A_+ + K(t) \dot{A}_+
\eeq
 
\beq
\label{2.11e}
A_1(u, A) = - \int_t^{\infty} dt'\ K(t-t') PJ (u, A)(t')
\eeq

\noi with
$$\omega = (-\Delta)^{1/2} \quad , \quad K(t) = \omega^{-1} \sin 
\omega t \quad , \quad \dot{K}(t) = \cos \omega t\ .$$

\noi In particular, $A_0^{\infty}$ is a solution of the free (vector 
valued) wave equation with initial data $(A_+ ,
\dot{A}_+)$ at $t = 0$, and $(A_+ ,\dot{A}_+)$ is naturally 
interpreted as the asymptotic state for $A$. In order to
ensure the condition $\nabla\cdot A = 0$, we assume that $\nabla 
\cdot A_+ = \nabla \cdot \dot{A}_+ = 0$. \par

We next perform the same change of variables as for the Hartree 
equation and for the WS system. That change of
variables is well adapted to the study of the asymptotic behaviour in 
time. The unitary
group
\beq
\label{2.12e}
U(t) = \exp (i(t/2)\Delta )
\eeq
\noi which solves the free Schr\"odinger equation can be written as
\beq
\label{2.13e}
U(t) = M(t) \ D(t) \ F \ M(t)
\eeq
\noi where $M(t)$ is the operator of multiplication by the function
\beq
\label{2.14e}
M(t) = \exp  \left ( i x^2/2t \right ) \ ,
\eeq
\noi $F$ is the Fourier transform and $D(t)$ is the dilation operator
\beq
\label{2.15e}
(D(t)\ f)(x) = (it)^{-n/2} \ f(x/t)
\eeq
\noi normalized to be unitary in $L^2$. We shall also need the 
operator $D_0(t)$ defined by
\beq
\label{2.16e}
\left ( D_0(t) f\right )(x) = f(x/t) \ .
\eeq
\noi We parametrize the Schr\"odinger function $u$ in terms of an 
amplitude $w$ and of a real phase $\varphi$ as
\beq
\label{2.17e}
u(t) = M(t) \ D(t) \exp [- i \varphi (t)] w(t) \ .
\eeq

\noi Correspondingly we change the variable for the vector potential 
from $A$ to $B$ according to
\beq
\label{2.18e}
A(t) = t^{-1} \ D_0(t) \ B(t)
\eeq

\noi and similarly for $A_0^{\infty}$ and $A_1$. \par

We now perform the change of variables (\ref{2.17e}) (\ref{2.18e}) on 
the system (\ref{2.6e}) (\ref{2.7e}).
Substituting (\ref{2.17e}) (\ref{2.18e}) into (\ref{2.6e}) and 
commuting the Schr\"odinger operator with $MD$, we obtain
\beq
\label{2.19e}
\left \{ i \partial_t + (2t^2)^{-1} (\nabla - i B)^2 + t^{-1} x \cdot 
B - t^{-1} \ g(w) \right \} \left ( \exp (- i
\varphi ) w \right ) = 0  \eeq

\noi where $g$ is defined by (\ref{2.4e}). Expanding the derivatives, 
using the condition $\nabla\cdot B = 0$ and introducing the notation
$s = \nabla \varphi$, we rewrite (\ref{2.19e}) as
\bea
\label{2.20e}
&&\Big \{ i \partial_t + \partial_t \varphi + (2t^2)^{-1} \Delta - i 
t^{-2} \left ( (s + B) \cdot \nabla + (1/2)
(\nabla \cdot s)\right ) - (2t^2)^{-1} (s +  B)^2  \nn \\
&&\hskip 4 truecm  + t^{-1} \ x \cdot B - t^{-1}\ g(w) \Big \} w = 0  \ .\eea

We next turn to the Maxwell equation (\ref{2.7e}). Substituting 
(\ref{2.17e}) (\ref{2.18e}) into the definition of
$J$, we obtain
\bea
\label{2.21e}
J(t) &=& t^{-3}\ D_0(t) \left \{ x|w|^2 + t^{-1} \ {\rm Im} \ \bar{w} 
\nabla w - (s + B) |w|^2 \right \} \nn \\
&=& t^{-3} \ D_0(t) \left ( M_a + t^{-1}  M_b \right )
\eea

\noi where
\bea
\label{2.22e}
&&M_a = x |w|^2 \ , \\
&&M_b = {\rm Im} \ \bar{w} \nabla w - (s + B) |w|^2 \ .
\label{2.23e}
\eea

\noi Using the identity
\beq
\label{2.24e}
- \int_t^{\infty} dt'\ K(t-t') t'^{-3-j} \ D_0(t') \ PM(t') = 
t^{-1-j} \ D_0 (t) \ F_j(M)
\eeq

\noi where $F_j$ is defined by
\beq
\label{2.25e}
F_j(M) = \int_1^{\infty} d\nu \ \omega^{-1} \sin (\omega (\nu - 1)) \ 
\nu^{-3-j} \ D_0 (\nu ) \ M(\nu t) \ ,
\eeq

\noi we rewrite (\ref{2.7e}) as
\beq
\label{2.26e}
B = B_0^{\infty} + B_1
\eeq

\noi with
\bea
\label{2.27e}
&&B_1 = B_a + B_b \\
\label{2.28e}
&&B_a \equiv B_a (w) = F_0(M_a) \\
&&B_b = t^{-1} \ F_1(M_b) \ .
\label{2.29e}
\eea

This completes the change of variables for the system (\ref{2.6e}) 
(\ref{2.7e}). Now however we have parametrized $u$
in terms of an amplitude $w$ and a phase $\varphi$ and we have only 
one equation for the two functions $(w, \varphi )$.
We arbitrarily impose a second equation, namely a Hamilton-Jacobi (or 
eikonal) equation for the phase $\varphi$, thereby
splitting the equation (\ref{2.20e})  into a system of two equations, 
the other one of which being a transport equation
for the amplitude $w$. There is a large amount of freedom in the 
choice of the equation for the phase. The role of the
phase is to cancel the long range terms in (\ref{2.20e}) coming from 
the interaction. All the interaction terms in
(\ref{2.20e}) having an explicit $t^{-2}$ factor are expected to be 
short range and are therefore left in the equation
for $w$. Such is also the case for the contribution of $B_b$ to 
$x\cdot B$ because of the $t^{-1}$ factor in
(\ref{2.29e}). The term $t^{-1} g(w)$ is clearly long range (of 
Hartree type), and is therefore included in the
$\varphi$ equation. The term $t^{-1}(x \cdot B_a)$ is also long 
range, but since it is less regular than the previous
one, it is convenient to split $x \cdot B_a$ into a short range and a 
long range part, namely
   \beq \label{2.30e}
x \cdot B_a = (x \cdot B_a)_S + (x \cdot B_a)_L
\eeq

\noi in the following way. We take $0 < \beta < 1$ and we define
\beq
\label{2.31e}
\left \{ \begin{array}{l} (x \cdot B_a)_S = F^{-1} \ \chi (|\xi | > 
t^{\beta}) \ F(x \cdot B_a)\\  \\(x \cdot B_a)_L =
F^{-1} \ \chi (|\xi | \leq t^{\beta}) \ F(x \cdot B_a)\\ \end{array} 
\right .  \eeq

\noi where $\chi ({\cal P})$ is the characteristic function of the 
set where ${\cal P}$ holds. (In practice we shall
need $0 < \beta < 1/2$ in most of the applications). Finally the 
contribution of $B_0^{\infty}$ to $B$ should be
considered short range. The Hamilton-Jacobi equation for the phase is 
then taken to be
  \beq \label{2.32e}
\partial_t \varphi = (2t^2)^{-1} \ s^2 + t^{-1} g(w) - t^{-1} (x 
\cdot B_a)_L (w) \ .
\eeq

We are now in a position to introduce the auxiliary system which will 
be used to study the MS system. From now on, we
restrict our attention to the case $B_0^{\infty} = 0$, namely to the 
case where the asymptotic state $(A_+,
\dot{A}_+)$ for the vector potential is zero. This is the technical 
meaning of the expression ``vanishing asymptotic
magnetic field'' used in the title and in the introduction of this paper. \par

The equation for $w$ is obtained by substituting (\ref{2.32e}) into 
(\ref{2.20e}). The resulting equation contains
$\varphi$ only through $s = \nabla \varphi$. The same property holds 
for the RHS of (\ref{2.32e}) and for
(\ref{2.26e})-(\ref{2.29e}). It is then convenient to replace 
(\ref{2.32e}) by its gradient so as to obtain a final
system containing only $s$. The Maxwell equation now has $B = B_1$. 
Furthermore we shall regard (\ref{2.28e}) as the
definition of $B_a$ as a function of $w$, we shall take $B_b$ as the 
dynamical variable and we shall regard
(\ref{2.27e}) (with $B = B_1$) as a change of variable from $B$ to 
$B_b$. The actual equation is then (\ref{2.29e})
regarded as an equation for $B_b$. We finally introduce the notation
\beq \label{2.33e}
Q(s, w) = s \cdot \nabla w + (1/2) (\nabla \cdot s) w
\eeq

\noi for the transport term. We then obtain the auxiliary system in 
the following form~:

\beq
\label{2.34e}
\left \{ \begin{array}{l} \partial_t w = i(2t^2)^{-1} \Delta w + 
t^{-2} Q(s + B, w) - i(2t^2)^{-1} (2B\cdot s + B^2 )
w \\ \\ \qquad + i t^{-1} \left ((x \cdot B_a)_S + x \cdot B_b \right 
) w \equiv L(w,s,B_b) w\\  \\ \partial_t s =
t^{-2} s \cdot \nabla s + t^{-1} \nabla g(w) - t^{-1} \nabla (x \cdot 
B_a)_L \\ \end{array} \right .
\eeq

\beq
\label{2.35e}
B_b = t^{-1} \ F_1 \left ( {\rm Im}\ \bar{w} \nabla w - (s + B) |w|^2 \right )
\eeq

\noi with $B = B_a + B_b$ and $B_a$ {\it defined} by (\ref{2.28e}) 
(\ref{2.22e}) (\ref{2.25e}). The phase $\varphi$ is
regarded as a derived quantity to be recovered from (\ref{2.32e}). 
The linear operator $L(w,s,B_b)$ is defined in an
obvious way by (\ref{2.34e}). Its dependence on $s$, $B_b$ is 
explicit, while its dependence on $w$ occurs through
$B_a$. Since $B_a$ and a fortiori $x\cdot B_a$ contains $xw$, it will 
be useful to have an explicit evolution equation
for $xw$. From (\ref{2.34e}) we obtain immediately
\beq
\label{2.36e}
\partial_t xw = L(w,s,B_b) xw - i t^{-2} \nabla w - t^{-2} (s+ B) w \ .
\eeq

The Cauchy problem for the system (\ref{2.34e}) (\ref{2.35e}) with 
initial data $(w, s)(t_0) = (w_0, s_0)$ for some
time $t_0$ is no longer a usual PDE Cauchy problem since $B_a$ 
depends on $w$ nonlocally in time and since
(\ref{2.35e}) is an integral equation in time. A convenient way to 
handle that difficulty is to replace that problem
by a partly linearized version thereof, namely
\bea
\label{2.37e}
&&\left \{ \begin{array}{l} \partial_t w' = L(w,s,B_b) w'\\ \\ 
\partial_t s' = t^{-2} s \cdot \nabla s' + t^{-1} \nabla
g(w) - t^{-1} \nabla (x \cdot B_a)_L \\ \end{array} \right . \\
&& \nn \\
&&B'_b = t^{-1} \ F_1 \left ( {\rm Im} \ \bar{w} \nabla w - (s+B) |w|^2 \right )\label{2.38e} \eea\noi still with $B = B_a + B_b$ and $B 
_a = B_a(w)$. Correspondingly, the evolution equation for $xw'$ 
becomes
\beq
\label{2.39e}
\partial_t xw' = L(w, s, B_b) xw' - it^{-2} \nabla w' - t^{-2} (s+ B) w' \ .
\eeq

\noi Solving (\ref{2.37e}) with suitable initial data for given $(w, 
s, B_b)$, together with (\ref{2.38e}), defines a
map $\Gamma : (w, s, B_b) \to (w', s', B'_b)$ and solving 
(\ref{2.34e}) (\ref{2.35e}) reduces to finding a fixed
point of $\Gamma$, which in favourable cases can be done for instance 
by contraction.\par

The first problem that we shall consider is whether the auxiliary 
system (\ref{2.34e}) (\ref{2.35e}) defines a
dynamics for large time. This will be the subject of Section 4 below. 
In particular we shall prove that the Cauchy
problem for that system is locally well posed in a neighborhood of 
infinity in time, namely that (\ref{2.34e})
(\ref{2.35e}) with initial data $(w, s) (t_0) = (w_0, s_0)$ has a 
unique solution defined in $[T, \infty )$ for some
$T \leq t_0$, with $t_0$ and $T$ suitably large depending on $(w_0, 
s_0)$, and with continuous dependence on the
data. Those results will be obtained through the use of the 
linearized system (\ref{2.37e})(\ref{2.38e}) by
following the method sketched above. In addition, we shall derive 
some asymptotic properties of the solutions
thereby obtained. In particular, for those solutions, $w(t)$ tends to 
a limit $w_+$ when $t \to \infty$.\par

The previous results are insufficient to construct the wave 
operators, namely to solve the Cauchy problem for
(\ref{2.34e}) (\ref{2.35e}) with infinite initial time because (i) $T 
\to \infty$ when $t_0 \to \infty$ and (ii)
the solutions are not estimated uniformly in $t_0$, so that it is not 
possible to perform the limit $t_0 \to \infty$ in
those results. In order to construct the wave operators, we follow 
instead the procedure explained at the beginning of
this section. We choose an asymptotic function $v$ and we look for 
solutions that behave asymptotically as $v$ when $t
\to \infty$. The asymptotic $v$ will be taken as a pair $(W,S)$ with 
$S = \nabla \phi$ and with $W(t)$ tending to a
limit $w_+$ as $t \to \infty$. This will provide the asymptotic form 
of $(w, s)$. No asymptotic form is needed for
$B_b$, because $B_b$ tends to zero at infinity. In order for $(W,S)$ 
to be an adequate asymptotic function, it has to
be an approximate solution of the system (\ref{2.34e}). This is 
achieved by solving that system approximately by
iteration and taking for $(W,S)$ an iterate of suitable order. In the 
present case, the first iteration turns out to be
sufficient, and we shall take accordingly
\beq \label{2.40e}
\left \{ \begin{array}{l} W(t) = U^*(1/t) w_+\\ \\ S(t) = 
\displaystyle{\int_1^t} dt' \ t'^{-1} \left ( \nabla g(W) -
\nabla (x\cdot B_a)_L(W) \right ) \\ \end{array} \right . \eeq

\noi for some given $w_+$ and for $t \geq 1$. Actually, at the cost 
of some loss of regularity, we could also make the
simpler choice
\beq \label{2.41e}
\left \{ \begin{array}{l} W(t) = w_+\\ \\ S(t) = 
\displaystyle{\int_1^t} dt' \ t'^{-1} \left ( \nabla g(w_+) - \nabla
(x\cdot B_a)_L(w_+) \right ) \ .\\ \end{array} \right . \eeq

\noi The latter $S(t)$ can be explicitly computed~:
\beq
\label{2.42e}
S(t) = \ell n \ t\ \nabla g(w_+) - \ell n \left ( t(\omega \vee 
1)^{-1/\beta} \vee 1 \right ) \nabla (x\cdot B_a) (w_+)
\eeq

\noi where the second $\ell n$ is defined in an obvious way in 
Fourier transformed variables. \par

In order to solve the system (\ref{2.34e}) (\ref{2.35e}) with $(w, 
s)$ behaving asymptotically as a given $(W, S)$
(possibly but not necessarily (\ref{2.40e})), we make a change of 
variables from $(w,s)$ to $(q, \sigma )$ defined by
\beq \label{2.43e}
(q, \sigma ) = (w, s) - (W, S) \ .
\eeq

\noi Substituting (\ref{2.43e}) into (\ref{2.34e}) (\ref{2.35e}) will 
yield a new auxiliary system for the variables
$(q, \sigma, B_b)$. For that purpose we introduce the following 
additional notation. We define
  \bea
\label{2.44e}
&&g(w_1, w_2) = (4 \pi |x|^{-1}) \ * \ {\rm Re}\ \bar{w}_1 \ w_2 \\
&&B_a(w_1, w_2) = F_0 (x\ {\rm Re}\ \bar{w}_1 \ w_2)
\label{2.45e}
\eea

\noi so that $g(w) = g(w, w)$ and $B_a (w) = B_a (w, w)$. We next 
define $B_* = B_a (W)$ and
\beq
\label{2.46e}
G \equiv G(q, W) = B_a (w) - B_* = B_a (q, q + 2W) \ .
\eeq

\noi We finally define the remainders
\bea
\label{2.47e}
&&\left \{ \begin{array}{l} R_1(w,s,B_b) = - \partial_t w + L(w, s, 
B_b) w\\ \\ R_2(w, s) = - \partial_t s + t^{-2} s
\cdot \nabla s + t^{-1} \nabla g(w) - t^{-1} \nabla (x \cdot 
B_a)_L(w)\\ \end{array} \right . \\
&& \nn \\
&&R_3(w, s, B_b) = - B_b + t^{-1} \ F_1({\rm Im}\ \bar{w} \nabla w - 
(s+ B) |w|^2)
\label{2.48e}
\eea

\noi so that the system (\ref{2.34e}) (\ref{2.35e}) can be rewritten 
as $R_i = 0$, $i = 1, 2, 3$, and for general $(w,
s, B_b)$, the remainders measure its failure to satisfy that system. 
Using the previous notation, we can rewrite (\ref{2.34e})
(\ref{2.35e}) in the new variables as follows
\beq
\label{2.49e}
\left \{ \begin{array}{l} \partial_t q = L(w,s,B_b)q + t^{-2} 
Q(\sigma + G + B_b, W) \\ \\ \quad - i t^{-2} \left ( B
\cdot \sigma + (G + B_b) (S + (B + B_*)/2)\right ) W\\ \\ \quad + i 
t^{-1} \left ((x \cdot G)_S + x \cdot B_b \right ) W
+ R_1(W,S,0)  \equiv L(w,s,B_b)q + \widetilde{R}_1\\  \\ \partial_t 
\sigma = t^{-2} (s \cdot \nabla \sigma + \sigma \cdot
\nabla S)  + t^{-1} \nabla g (q, q+2W) - t^{-1} \nabla (x \cdot G)_L 
\\ \\ \qquad + R_2(W,S) \equiv t^{-2} s
\cdot \nabla \sigma + \widetilde{R}_2\end{array} \right . \eeq

\bea
\label{2.50e}
B_b &=& t^{-1} \ F_1 \Big ( {\rm Im}(\bar{w} \nabla q + \bar{q} 
\nabla W) - (s + B) \left ( |q|^2 + 2 {\rm Re}\ \bar{q} W \right
)\nn \\    &&- ( \sigma + G + B_b ) |W|^2 \Big ) + R_3 (W,S,0) \equiv 
\widetilde{R}_3 \ .
\eea

\noi For the same reason as above, we also need the evolution 
equation for $xq$. We define the linear operator $L_1$
by rewriting the evolution equation for $q$ in the form
\beq
\label{2.51e}
\partial_t q = L(w,s,B_b) q + L_1 W + R_1 (W,S, 0) \ .
\eeq

\noi The evolution equation for $xq$ then becomes
\bea
\label{2.52e}
\partial_t xq &=& L(w,s,B_b) xq + L_1 xW - it^{-2} \nabla q - t^{-2} 
(s+ B) q \nn \\
&& - t^{-2} (\sigma + G + B_b) W + x R_1(W,S,0) \ .
\eea

We shall also need to rewrite the partly linearized system 
(\ref{2.37e}) (\ref{2.38e}) in terms of the new variables
$(q, \sigma )$ defined by (\ref{2.43e})  and $(q', \sigma ')$ defined 
similarly by
\beq
\label{2.53e}
(q', \sigma ') = (w' , s') - (W, S) \ .
\eeq

The system (\ref{2.37e}) (\ref{2.38e}) then becomes
\bea
\label{2.54e}
&&\left \{ \begin{array}{l} \partial_t q' = L(w, s, B_b) q' + 
\widetilde{R}_1\\ \\\partial_t \sigma ' = t^{-2}
s\cdot \nabla \sigma ' +\widetilde{R}_2 \\ \end{array} \right . \\ && \nn \\
&&B'_b = \widetilde{R}_3
\label{2.55e}
\eea

\noi where $\widetilde{R}_i$, $i = 1,2,3$, are defined by 
(\ref{2.49e}) (\ref{2.50e}), and the equation (\ref{2.39e})
becomes
\bea
\label{2.56e}
\partial_t xq' &=& L(w,s,B_b) xq' + L_1 xW - it^{-2} \nabla q' - 
t^{-2} (s+ B) q' \nn \\
&& - t^{-2} (\sigma + G + B_b) W + x R_1(W,S,0) \ .
\eea

The main technical result of this paper is the construction of 
solutions $(q, \sigma, B_b)$ of the auxiliary system (\ref{2.49e})
(\ref{2.50e}) defined for large time and tending to zero at infinity 
in time. That construction is performed by
solving the Cauchy problem for the linearized system (\ref{2.54e}) 
(\ref{2.55e}), first for finite intial time $t_0$,
and then for infinite initial time by taking the limit $t_0 \to 
\infty$. One then proves the existence of a fixed
point for the map $\Gamma : (q, \sigma, B_b) \to (q', \sigma ', 
B'_b)$ thereby defined by a contraction method, as
mentioned above. With that result available, it is an easy matter to 
construct the modified wave operator $\Omega$ for
the MS system in the form (\ref{2.6e}) (\ref{2.7e}). We start from 
the asymptotic state $u_+$ for $u$ and we define
$w_+ = Fu_+$. The asymptotic state ($A_+, \dot{A}_+)$ for $A$ is 
taken to be zero. We define $(W,S)$ by (\ref{2.40e}).
We solve the system (\ref{2.49e}) (\ref{2.50e}) for $(q, \sigma , 
B_b)$ as indicated above. Through (\ref{2.43e}), this
yields a solution $(w, s, B_b)$ of the auxiliary system (\ref{2.34e}) 
(\ref{2.35e}) defined for large time. From $s$ we
reconstruct the phase $\varphi$ by using (\ref{2.32e}). We finally 
substitute $(w, \varphi , B_b)$ into (\ref{2.17e})
(\ref{2.18e}) with $B = B_a + B_b$ and $B_a$ defined by 
(\ref{2.28e}). This yields a solution $(u, A)$ of the system
(\ref{2.6e}) (\ref{2.7e}) defined for large time. The modified wave 
operator is the map $u_+ \to (u, A)$ thereby
obtained. \par

The main result of this paper is the construction of $(u, A)$ from 
$u_+$, as described above, together with the
asymptotic properties of $(u, A)$ that follow from that construction. 
It will be stated below in full mathematical
detail in Proposition 6.1. We give here only a heuristic preview of 
that result, stripped from most technicalities. \\

\noi {\bf Proposition 2.1.} {\it Let $k > 3/2$, $0 < \beta < 1/2$ and 
let $\alpha > 1$ be such that $\beta(\alpha + 1) \geq
1$. Let $u_+$ be such that $w_+ = Fu_+ \in H^{k+\alpha + 1}$ and $x 
w_+ \in H^{k+\alpha}$. Define $(W, S)$ by
(\ref{2.40e}). Then \par

(1) There exists $T = T(w_+)$, $1 \leq T < \infty$, such that the 
auxiliary system (\ref{2.34e}) (\ref{2.35e}) has a
unique solution $(w, s, B_b)$ in a suitable space, defined for $t 
\geq T$, and such that $(w - W, s- S, B_b)$ tends to
zero in suitable norms when $t \to \infty$. \par

(2) There exist $\varphi$ and $\phi$ such that $s = \nabla \varphi$, 
$S = \nabla \phi$, $\phi (1) = 0$ and such that
$\varphi - \phi$ tends to zero in suitable norms when $t \to \infty$. 
Define $(u, A)$ by (\ref{2.17e}) (\ref{2.18e})
with $B = B_a + B_b$ and $B_a$ defined by (\ref{2.28e}). Then $(u, 
A)$ solves the system (\ref{2.6e})
(\ref{2.7e}) for $t \geq T$ and $(u, A)$ behaves asymptotically as 
$(MD \exp (- i \phi ) W, t^{-1} D_0 B_a
(W))$ in the sense that the difference tends to zero in suitable 
norms (for which each term separately is
$O(1)$) when $t \to \infty$.} \\

We now describe the contents of the technical parts of this paper, 
namely Sections 3-6. In
Section 3, we introduce some notation, define the relevant function 
spaces and collect a number
of preliminary estimates. In Section 4, we study the Cauchy problem 
for large time for the
auxiliary system (\ref{2.34e}) (\ref{2.35e}). We solve the Cauchy 
problem with finite initial time for the
linearized system (\ref{2.37e}) (Proposition 4.1), we prove a number 
of uniqueness results
for the system (\ref{2.34e}) (\ref{2.35e}) (Proposition 4.2), we 
solve the Cauchy problem for the system (\ref{2.34e}) (\ref{2.35e})
for large but finite $t_0$ in the special case $A_0^{\infty} = 0$ 
(Proposition 4.3) and we finally derive some
asymptotic properties of the solutions thereby obtained (Proposition 
4.4). In Section 5, we study the Cauchy problem
at infinity for the auxiliary system (\ref{2.34e}) (\ref{2.35e}) in 
the difference form (\ref{2.49e}) (\ref{2.50e}).
We prove the existence of solutions first for the linearized system 
(\ref{2.54e}) (\ref{2.55e}), for $t_0$ finite
(Proposition 5.1) and infinite (Proposition 5.2), and then for the 
nonlinear system (\ref{2.49e}) (\ref{2.50e}) for
$t_0$ infinite (Proposition 5.3). Finally in Section 6, we construct 
the modified wave operators for the system
(\ref{2.6e}) (\ref{2.7e}) from the results previously obtained for 
the system (\ref{2.49e}) (\ref{2.50e}) and we
derive the asymptotic estimates for the solutions $(u,A)$ in their 
range that follow from the previous estimates
(Proposition 6.1).

\mysection{Notation and preliminary estimates}
\hspace*{\parindent} In this section, we define the function spaces 
where we shall study the auxiliary system (\ref{2.34e})
(\ref{2.35e}) and we collect a number of estimates which will be used 
throughout this paper. In addition to the
standard Sobolev spaces $H^k$, we shall use the associated 
homogeneous spaces $\dot{H}^k$ with norm $\parallel u ;
\dot{H}^k\parallel \ = \ \parallel \omega^k u \parallel_2$, where 
$\omega = (- \Delta)^{1/2}$ and the spaces
$$K^k = \dot{H}^1 \cap \dot{H}^k \ ,$$

\noi where it is understood that $\dot{H}^1 \subset L^6$. We shall 
use the notation
\beq
\label{3.1e}
|w|_k = \ \parallel w; H^k\parallel\quad , \quad  |s|_k^{\dot{}} \ = 
\ \parallel s ; K^k
\parallel \ .  \eeq

We shall look for solutions of the auxiliary system in spaces of the 
type ${\cal C}(I, X^{k})$ where $I$ is an
interval and
\beq
\label{3.2e}
X^{k} = \left \{ (w, s, B) : w \in H^{k+1} \ , \ xw \in H^{k} \ , \ s 
\in K^{k+2} \ , \ B \in K^{k+1} \ , \ x\cdot B \in
K^{k+1} \right \} \ .  \eeq

\noi For the needs of this paper, we could have replaced $\dot{H}^1$ 
in the definition of $K^k$ by $\dot{H}^{k_0}$ for
some $k_0$ with $1/2 < k_0 < 3/2$. We have chosen $k_0 = 1$ for 
simplicity. \par

We shall use extensively the following Sobolev inequalities, stated 
here in ${I\hskip-1truemm
R}^n$, but to be used only for $n = 3$. \\

\noi {\bf Lemma 3.1.} {\it Let $1 < q$, $r < \infty$, $1 < p \leq 
\infty$ and $0 \leq j < k$.
If $p = \infty$, assume that $k - j > n/r$. Let $\sigma$ satisfy $j/k 
\leq \sigma \leq 1$ and
$$n/p - j = (1 - \sigma )n/q + \sigma (n/r - k) \ .$$
\noi Then the following inequality holds
\beq
\label{3.3e}
\parallel \omega^j u \parallel_p \ \leq C \parallel u \parallel_q^{1 
- \sigma} \ \parallel
\omega^ku \parallel_r^{\sigma} \ .  \eeq} \\
\indent The proof follows from the Hardy-Littlewood-Sobolev (HLS) 
inequality (\cite{11r},
p.~117) (from the Young inequality if $p = \infty$), from 
Paley-Littlewood theory and
interpolation.\\

We shall also use extensively the following Leibnitz and commutator 
estimates.\\

\noi {\bf Lemma 3.2.} {\it Let $1 < r, r_1, r_3 < \infty$ and
$$1/r = 1/r_1 + 1/r_2 = 1/r_3 + 1/r_4 \ .$$
\noi Then the following estimates hold
\beq
\label{3.4e}
\parallel \omega^m (uv) \parallel_r \ \leq C \left ( \parallel 
\omega^m u \parallel_{r_1}
  \ \parallel v \parallel_{r_2} + \parallel  \omega^m v \parallel_{r_3} \
  \parallel u \parallel_{r_4} \right ) \eeq
\noi for $m \geq 0$,
\beq
\label{3.5e}
\parallel [\omega^m , u]v\parallel_r \ \leq C \left ( \parallel 
\omega^m u \parallel_{r_1} \
  \parallel v \parallel_{r_2} + \parallel  \omega^{m -1} v \parallel_{r_3} \
  \parallel \nabla u \parallel_{r_4} \right ) \eeq
\noi for $m \geq 1$, where $[\ ,\ ]$ denotes the commutator, and
\beq
\label{3.6e}
\parallel [\omega^m , u]v\parallel_2 \ \leq C \parallel  F(\omega^m 
u) \parallel_{1} \
  \parallel v \parallel_{2} \eeq
\noi for $0 \leq m \leq 1$.}\\

\noi {\bf Proof.} The proof of (\ref{3.4e}) (\ref{3.5e}) is given in 
\cite{12r} \cite{13r} with $\omega$
replaced by $<\omega >$ and follows therefrom by a scaling argument. 
The proof of (\ref{3.6e}) follows from
\begin{eqnarray*}
|F([\omega^m, u]v)(\xi ) | &=& \left |\int d\eta | |\xi|^m - |\eta 
|^m |\widehat{u}(\xi - \eta) \widehat{v}(\eta
) \right |\\
&\leq & \int d\eta |\xi - \eta|^m |\widehat{u}(\xi - \eta)|\ 
|\widehat{v}(\eta )|
\end{eqnarray*}

\noi and from the Young inequality and the Plancherel theorem.\par\nobreak
\hfill $\sq$\par

We shall also need the following consequence of Lemma 3.2. \\

\noi {\bf Lemma 3.3.} {\it Let $m \geq 0$ and $1 < r < \infty$. Then 
the following estimate
holds}
\beq
\label{3.7e}
\parallel \omega^m (e^{\varphi} - 1) \parallel_r \ \leq \ \parallel 
\omega^m \varphi
\parallel_{r} \exp \left ( C \parallel \varphi  \parallel_{\infty} 
\right ) \ . \eeq

\vskip 3 truemm
\noi {\bf Proof.} For any integer $n \geq 2$, we estimate
\begin{eqnarray*}
a_n &\equiv& \parallel \omega^m \ \varphi^n\parallel_r \ \leq C \left 
( \parallel  \omega^m
\varphi \parallel_{r} \
  \parallel \varphi \parallel_{\infty}^{n-1} + \parallel  \omega^m 
\varphi^{n-1} \parallel_{r} \
  \parallel \varphi \parallel_{\infty} \right )  \nn \\
&=& C \left ( a_1 \ b^{n-1} + a_{n-1} \ b \right )
\end{eqnarray*}

\noi by (3.4), where $b = \parallel \varphi \parallel_{\infty}$ and 
we can assume $C \geq 1$
without loss of generality. It follows easily from that inequality that
$$a_n \leq n(Cb)^{n-1} \ a_1$$
\noi for all $n \geq 1$, from which (\ref{3.6e}) follows by expanding 
the exponential.\par\nobreak
\hfill $\sq$\par

We shall apply Lemma 3.2 in the form of Lemma 3.4 below which we 
state for clarity in general dimension $n$ and with
general $k_0 < n/2$ and $K^m = \dot{H}^m \cap \dot{H}^{k_0}$. \\

\noi {\bf Lemma 3.4.} {\it Let $\bar{m} > n/2$. The following 
inequalities hold~: \par

(1) Let $0 \leq m \leq \bar{m}$. Then
\bea
\label{3.8e}
&\parallel fg; \dot{H}^m \parallel \ \leq \ C |f|_{\bar{m}}^{\dot{}} 
\ \parallel g; \dot{H}^m\parallel&\quad \hbox{for} \
m < n/2 \ , \nn \\
&\qquad \leq \ C |f|_{\bar{m}}^{\dot{}} \ | g|_m^{\dot{}} &\quad 
\hbox{for} \ n/2 \leq m \leq \bar{m} \ , \\
\label{3.9e}
&|fg|_m^{\dot{}} \leq \ C |f|_{\bar{m}}^{\dot{}} \ 
|g|_m^{\dot{}}\hskip 2 truecm  &\quad \hbox{for} \ k_0 \leq m \leq
\bar{m} \ , \\
&|fg|_m \leq \ C |f|_{\bar{m}}^{\dot{}} \ |g|_m\hskip 2 truecm 
&\quad \hbox{for} \ 0 \leq m \leq \bar{m} \ .
\label{3.10e}
\eea

(2) Let $0 \leq m \leq \bar{m} +1$. Then
\bea
\label{3.11e}
\parallel [\omega^m, f] \nabla g \parallel_2 \ &\leq \ C |\nabla 
f|_{\bar{m}}^{\dot{}} \ \parallel \omega^m g\parallel_2 &
\hbox{for} \ 0 \leq m < n/2 + 1 \ , \\
\label{3.12e}
&\leq \ C |\nabla f|_{\bar{m}}^{\dot{}} \ | \nabla g|_{m-1}^{\dot{}}\ 
&\hbox{for} \ n/2+1 \leq m \leq
\bar{m} + 1 \ , \\
\label{3.13e}
&\leq \ C |\nabla f|_{\bar{m}}^{\dot{}} \ | g|_{m}^{\dot{}}\qquad 
&\hbox{for} \ k_0 \leq m \leq \bar{m} + 1 \ ,  \\
&\leq \ C |\nabla f|_{\bar{m}}^{\dot{}} \ | g|_{m}\qquad  &\hbox{for} 
\ 0 \leq m \leq \bar{m} + 1 \ .
\label{3.14e}
\eea

(3) Let $n \geq 3$ and $m \geq k_0$, $1 < m \leq \bar{m} + 1$. Then}

\beq
\label{3.15e}
\parallel [\omega^m, f] g \parallel_2 \ \leq \ C |f|_{m}^{\dot{}}\ | 
g|_{\bar{m}}^{\dot{}} \ .
\eeq

\vskip 3 truemm

\noi {\bf Proof.} \underbar{Part (1)}. By Lemma 3.2, we estimate
\beq
\label{3.16e}
\parallel \omega^m fg \parallel_2 \ \leq \ C \left ( \parallel f 
\parallel_{\infty} \  \parallel \omega^m g \parallel_2 \ + \ \parallel
\omega^m f \parallel_{n/ \delta} \  \parallel g \parallel_{r}  \right )
\eeq

\noi with $0 < \delta = \delta (r) \leq n/2$. \par

For $ m < n/2$, we choose $\delta = m$ and continue (\ref{3.16e}) by
$$\cdots \leq \ C \left ( \parallel f \parallel_{\infty} \ + \ 
\parallel \omega^{n/2} f \parallel_2 \right )
\parallel \omega^m g \parallel_{2} \ \leq \ C |f|_{\bar{m}}^{\dot{}} 
\  \parallel \omega^m g \parallel_{2}  $$

\noi by Sobolev inequalities, which yields (\ref{3.8e}) in this case. \par

For $m > n/2$, we choose $\delta = n/2$ and continue (\ref{3.16e}) by
$$\cdots \leq \ C \left ( \parallel f \parallel_{\infty} \ \parallel 
\omega^{m} g \parallel_2 \ + \
\parallel \omega^m f \parallel_{2} \  \parallel g \parallel_{\infty} 
\right ) \leq \ C |f|_{m}^{\dot{}} \ |g|_{m}^{\dot{}}$$

\noi which yields (\ref{3.8e}) since $m \leq \bar{m}$. \par

For $m = n/2$, we take $k_0 \vee (n - \bar{m}) \leq \delta < n/2$ and 
estimate the last term in (\ref{3.16e}) by
$$C \parallel \omega^{n-\delta} f \parallel_{2} \  \parallel 
\omega^{\delta} g \parallel_{2}$$

\noi by Sobolev inequalities, which yields again (\ref{3.8e}) in that 
case. \par

Finally (\ref{3.9e}) and (\ref{3.10e})  are immediate consequences of 
(\ref{3.8e}) .\\

\noi \underbar{Part (2)}. For $m > 1$, we estimate by Lemma 3.2
\beq
\label{3.17e}
\parallel [\omega^m ,f] \nabla g \parallel_2 \ \leq \ C \left ( 
\parallel \nabla f \parallel_{\infty} \  \parallel \omega^m g 
\parallel_2
\ + \ \parallel \omega^m f \parallel_{n/ \delta} \  \parallel \nabla 
g \parallel_{r}  \right )
\eeq

\noi with $0 < \delta = \delta (r) \leq n/2$. \par

For $m < n/2 + 1$, we choose $\delta = m - 1$ and continue (\ref{3.17e}) by
$$\cdots \leq \ C \left ( \parallel f \parallel_{\infty} \ + \ 
\parallel \omega^{n/2} \nabla f \parallel_2 \right )
\parallel \omega^m g \parallel_{2} $$

\noi by Sobolev inequalities, which yields (\ref{3.11e}) in this case. \par

For $m > n/2 + 1$, we choose $\delta = n/2$ and continue (\ref{3.17e}) by
$$\cdots \leq \ C \left ( \parallel \nabla f \parallel_{\infty} \ 
\parallel \omega^{m} g \parallel_2 \ + \
\parallel \omega^m f \parallel_{2} \  \parallel \nabla g 
\parallel_{\infty}  \right ) \leq \ C |\nabla f|_{m-1}^{\dot{}} \
|\nabla g|_{m-1}^{\dot{}}$$

\noi which yields (\ref{3.12e}) since $m \leq \bar{m} + 1$. \par

For $m = n/2 + 1$, we take $k_0 \vee (n - \bar{m}) \leq \delta < n/2$ 
and estimate the last term in (\ref{3.17e}) by
$$C \parallel \omega^{n-\delta} \nabla f \parallel_{2} \  \parallel 
\omega^{\delta} \nabla g \parallel_{2}$$

\noi by Sobolev inequalities, which yields again (\ref{3.12e}) in 
this case. \par

For $0 \leq m \leq 1$, (\ref{3.11e}) follows immediately from Lemma 
3.2 and from the inclusion $K^{\bar{m}} \subset F(L^1)$. \par

Finally (\ref{3.13e}) and (\ref{3.14e}) are immediate consequences of 
(\ref{3.11e}) (\ref{3.12e}). \\

\noi \underbar{Part (3)}. By Lemma 3.2, we estimate
\beq
\label{3.18e}
\parallel [\omega^m ,f] g \parallel_2 \ \leq \ C \left ( \parallel 
\omega^m f \parallel_2 \ \parallel g \parallel_{\infty} \  +
\ \parallel \nabla f \parallel_{n/ \delta} \  \parallel \omega^{m-1} 
g \parallel_{r}  \right )
\eeq

\noi with $0 \leq \delta = \delta (r) < n/2$. We estimate the last 
term in (\ref{3.18e}) by Sobolev inequalities as
$$C \parallel \omega^{1+ n/2 -\delta} f \parallel_{2} \  \parallel 
\omega^{m-1 + \delta} \parallel_{2} \ \leq \ C |f|_{m}^{\dot{}} \
|g|_{\bar{m}}^{\dot{}}$$

\noi provided
\begin{eqnarray*}
&& 1 + n/2 - m \leq \delta \leq 1 + n/2 - k_0 \\
&&k_0 - m + 1 \leq \delta \leq \bar{m} + 1 - m \ ,
\end{eqnarray*}
 
\noi and the various conditions on $\delta$ are compatible for $m 
\geq k_0$ and $1 < m \leq \bar{m} + 1$. This proves
(\ref{3.15e}).\par\nobreak \hfill $\sq$\par

We next give some estimates of the various components of $B_1$, 
defined by (\ref{2.27e})-(\ref{2.29e}) and (\ref{2.31e}). It follows
immediately from (\ref{2.31e}) that
\beq \label{3.19e} \parallel \omega^{m} (x \cdot B_a)_S \parallel_{2} 
\  \leq \ t^{\beta (m-p)} \parallel \omega^{p} (x \cdot
B_a)_S\parallel_{2} \ \leq \ \ t^{\beta (m-p)} \parallel \omega^{p} 
(x \cdot B_a)\parallel_{2} \eeq

\noi for $m \leq p$, and similarly
\beq \label{3.20e}
\parallel \omega^{m} (x \cdot B_a)_L \parallel_{2} \  \leq \ t^{\beta 
(m-p)} \parallel \omega^{p} (x \cdot
B_a)_L\parallel_{2} \ \leq \ \ t^{\beta (m-p)} \parallel \omega^{p} 
(x \cdot B_a)\parallel_{2} \eeq

\noi for $m \geq p$. \par

We next estimate $F_j (M)$ defined by (\ref{2.25e}). From 
(\ref{2.25e}) and from the dilation identity
\beq
\label{3.21e}
\parallel \omega^{m} D_0(\nu) f \parallel_{2} \ = \nu^{-m+3/2} 
\parallel \omega^{m} f\parallel_{2} \ ,
\eeq

\noi it follows immediately that
\beq
\label{3.22e}
\parallel \omega^{m+1} F_j(M) \parallel_2 \ \leq \ C \ I_{m+j} \left 
( \parallel \omega^m M \parallel_2 \right )
\eeq

\noi for $j = 0, 1$, where
\beq
\label{3.23e}
(I_m(f))(t) = \int_1^{\infty} d\nu \ \nu^{-m-3/2} \ f(\nu t)
\eeq

\noi or equivalently
\beq
\label{3.24e}
(I_m(f))(t) = t^{m+1/2} \int_t^{\infty} dt'\ t'^{-m-3/2} \ f(t')
\eeq

\noi for $t > 0$. In particular for $m \geq 1$ and $j = 0,1$,
\beq
\label{3.25e}
|F_j(M)|_m^{\dot{}} \ \leq \ C I_j\left ( |M|_m \right ) \ .
\eeq

We next estimate $x \cdot F_j(M)$. Using the commutation relation
\beq
\label{3.26e}
[x, P] = (n-1) \Delta^{-1} \nabla \ ,
\eeq

\noi easily proved in Fourier transformed variables, and estimating
$$| \sin ((\nu - 1) \omega)| \ \leq \ \nu \omega \vee 1 \ ,$$

\noi we obtain
\beq
\label{3.27e}
\parallel \omega^{m+1} x \cdot F_j(M) \parallel_2 \ \leq \ C \ 
I_{m-1+j} \left ( \parallel \omega^m (x \otimes M) \parallel_2 \ +
\ \parallel \omega^{m} M\parallel_{2}\right ) \ ,\eeq

\noi and in particular for $m \geq 1$,
\beq
\label{3.28e}
|x \cdot F_1(M)|_{m+1}^{\dot{}} \vee |\nabla x \cdot 
F_0(M)|_m^{\dot{}} \ \leq \ I_0 \left ( |x \otimes M|_m + |M|_m 
\right ) \ .
\eeq

The estimates (\ref{3.25e}) (\ref{3.28e}) will be the main tools used 
to estimate $B_a$ and $B_b$ given by (\ref{2.28e}) (\ref{2.29e}).

\mysection{Cauchy problem and asymptotics for the\break \noindent 
auxiliary system}
\hspace*{\parindent}
In this section, we study the Cauchy problem for the auxiliary system 
(\ref{2.34e})(\ref{2.35e}) for large but finite initial time, and
we derive asymptotic properties in time of its solutions. \par

The basic tool of this section consists of a priori estimates for 
suitably regular
solutions of the linearized system (\ref{2.37e}) (\ref{2.38e}). Those 
estimates can be proved by a
regularisation and limiting procedure and hold in the integrated form 
at the available
level of regularity. For brevity, we shall state them in differential 
form and we shall
restrict the proof to the formal computation. \par

We first estimate a single solution of the linearized system 
(\ref{2.37e}) (\ref{2.38e}) at the level of
regularity where we shall eventually solve the auxiliary system 
(\ref{2.34e}) (\ref{2.35e}). \\

\noi {\bf Lemma 4.1.} {\it Let $k > 3/2$, let $T \geq 1$, $I = 
[T,\infty )$ and let $(w, s, B_b) \in {\cal C} (I, X^k)$ with $|w|_k 
\vee
|xw|_k \in L^{\infty}(I)$. Let
\beq
\label{4.1e}
a = \ \parallel |w|_k \vee |xw|_k ; L^{\infty}(I)\parallel \ .
\eeq

\noi Let $I' \subset I$ be an interval, let $(w', s')$ be a solution 
of the system (\ref{2.37e}) with $(w', s', 0) \in {\cal C}(I',X^k)$
and let $B'_b$ be defined by (\ref{2.38e}). Let $0 \leq \theta \leq 
1$ and $k \leq \ell \leq k+2$. Then the following estimates hold~:
\beq
\label{4.2e}
|B_a|_{k+1}^{\dot{}} \leq \ C \ I_0 \left ( |w|_k \ |xw|_k \right ) \ 
\leq \ C\ a^2 \ ,
\eeq
\beq
\label{4.3e}
\parallel x\cdot B_a; \dot{H}^{k+1} \parallel \ \leq \ |\nabla (x 
\cdot B_a)|_k^{\dot{}} \ \leq \ C\ I_0 \Big ( |xw|_k (|xw|_k +
|w|_k)\Big ) \ \leq \ C\ a^2 \eeq

\noi for all $t \in I$.\par
\bea
\label{4.4e}
&&\left  | \partial_t| w'|_{k+\theta}  \right | \ \leq \ C\ t^{-2} 
\Big \{ |\nabla (s + B)|_k^{\dot{}} \ |w'|_{k+\theta} +
|\nabla \cdot s|_{k+\theta}^{\dot{}}\  |w'|_k \nn \\
&&+ \left ( |s|_{k+\theta}^{\dot{}} + |B|_{k+\theta}^{\dot{}}\right ) 
|B|_{k+\theta}^{\dot{}} \ |w'|_k \Big \} + C\ t^{-1-\beta (1 -
\theta)}  \parallel x \cdot B_a; \dot{H}^{k+1} \parallel \ |w'|_k \nn \\
&&+ C \ t^{-1} |x \cdot B_b|_{k+ \theta}^{\dot{}} \ |w'|_k \equiv M_4 
(\theta , w') \ ,
\eea
\beq
\label{4.5e}
\left  | \partial_t|x w'|_{k}  \right | \ \leq \ M_4(0,xw') + C \ 
t^{-2} \left ( |w'|_{k+1} + |s+ B|_{k}^{\dot{}}\  |w'|_k \right ) \ ,
\eeq
\bea
\label{4.6e}
&&\left  | \partial_t| s'|_{\ell}^{\dot{}}  \right | \ \leq \ C\ 
t^{-2} \left (  |s|_{k+1}^{\dot{}} \ |s'|_{\ell}^{\dot{}} + \chi (\ell
\geq k+1)|s|_{\ell}^{\dot{}}\  |s'|_{k+1}^{\dot{}} \right )  \nn \\
&&+ C \ t^{-1} |w|_{k} \ |w|_{\ell - 1} + C\ t^{-1+ \beta (\ell - k)} 
|\nabla (x \cdot B_a)|_k^{\dot{}}
\eea

\noi for all $t \in I'$,
\beq
\label{4.7e}
|B'_b|_{k+\theta}^{\dot{}} \leq \ C \ \ t^{-1}\ I_1 \left ( |w|_k ( 
|w|_{k+\theta} + |s + B|_{k}^{\dot{}} |w|_k ) \right )
\eeq
\beq
\label{4.8e}
|x \cdot B'_b|_{k+\theta}^{\dot{}} \leq \ C \ \ t^{-1}\ I_0 \left ( ( 
|w|_k + |xw|_k) ( |w|_{k+\theta} + |s + B|_{k}^{\dot{}} |w|_k )
\right )  \eeq

\noi for all $t \in I$.} \\

\noi {\bf Remark 4.1.} The statements on $B'_b$ are non empty only in 
so far as the integrals over $\nu$ in the RHS of (\ref{4.7e})
(\ref{4.8e}) are absolutely convergent. This requires suitable 
assumptions on the behaviour of $(w, s, B)$ at infinity in time. Such
assumptions will be made in due course. The only assumption of this 
type made so far is (\ref{4.1e}) which ensures (\ref{4.2e})
(\ref{4.3e}), thereby making (\ref{4.4e}) (\ref{4.5e}) (\ref{4.6e}) 
into non empty statements.\\

\noi {\bf Proof.} (\ref{4.2e}) and (\ref{4.3e}) follow immediately 
from (\ref{3.25e}) (\ref{3.28e}) and from (\ref{3.10e}) with $m =
\bar{m} = k$. \par

We next estimate $w'$ in $H^{k+\theta}$. It is clear from 
(\ref{2.37e}) that $\parallel w'\parallel_2 = const$. Let $m = k + 
\theta$.
We estimate by a standard energy method
\bea
\label{4.9e}
&&\left  | \partial_t\parallel \omega^mw'\parallel_2 \right | \ \leq 
\ t^{-2} \Big \{ \parallel [\omega^m, s+B] \cdot \nabla w'
\parallel_2 \ + \  \parallel [\omega^m, \nabla \cdot s] w'
\parallel_2\nn \\
&&+ \ \parallel [\omega^m, 2B\cdot s + B^2] w' \parallel_2 \Big \} + 
t^{-1} \parallel [\omega^m, (x\cdot B_a)_S] w' \parallel_2\nn \\
&&+\ t^{-1} \parallel [\omega^m, x\cdot B_b] w' \parallel_2  \ .
\eea

We estimate the first norm in the RHS by (\ref{3.14e}) with $\bar{m} 
= k$ and the other norms by (\ref{3.15e}) with $\bar{m} = k$.
Furthermore, by (\ref{3.19e})
  \beq
\label{4.10e}
\left |(x\cdot B_a)_S \right |_{k+\theta}^{\dot{}}  \ \leq \ 
t^{-\beta (1 - \theta)} \parallel  x \cdot B_a; \dot{H}^{k+1} 
\parallel
\eeq

\noi so that the RHS of (\ref{4.9e}) is estimated by that of 
(\ref{4.4e}), which together with $L^2$ norm conservation, yields
(\ref{4.4e}). \par

We next estimate $xw'$ in $H^k$, starting from (\ref{2.39e}). We obtain
\bea
\label{4.11e}
&&\left  | \partial_t\parallel \omega^k xw'\parallel_2 \right | \ 
\leq \ \hbox{terms containing $xw'$} \nn \\
&&+  t^{-2} \left ( \parallel \omega^{k+1}w' \parallel_2 \ + \ 
\parallel \omega^{k} ((s+B)w') \parallel_2 \right )
\eea

\noi where the terms containing $xw'$ are obtained from the RHS of 
(\ref{4.9e}) by replacing $w'$ by $xw'$ and taking $m = k$. Those
terms are estimated in the same way as before. Estimating the last 
norm in (\ref{4.11e}) by (\ref{3.15e}) and using the obvious $L^2$
estimate
\beq
\label{4.12e}
\left  | \partial_t\parallel xw'\parallel_2 \right | \ \leq \ t^{-2} 
\left ( \parallel \nabla w' \parallel_2 \ + \ \parallel
(s+B) \parallel_{\infty} \ \parallel w' \parallel_2  \right ) \eeq

\noi yields (\ref{4.5e}). \par

We next estimate $s'$. From (\ref{2.37e}) we obtain
\bea
\label{4.13e}
&&\left  | \partial_t\parallel \omega^{\ell} s'\parallel_2 \right | \ 
\leq \ t^{-2} \left (  \parallel [\omega^{\ell}, s] \nabla
s' \parallel_2 \ + \  \parallel (\nabla \cdot s) 
\omega^{\ell}s'\parallel_2 \right )\nn \\
&&+C\ t^{-1} \parallel [\omega^{\ell - 1}|w|^2 \parallel_2 \ + \ 
t^{-1} \parallel \omega^{\ell + 1}(x \cdot B_a)_L \parallel_2 \ .
\eea

We estimate the first norm in the RHS by (\ref{3.11e}) (\ref{3.12e}) 
with $m = \ell$ and $\bar{m} = k$ if $\ell \leq k+1$, while for
$\ell \geq k + 1$
\bea
\label{4.14e}
\parallel [\omega^{\ell}, s] \nabla s'\parallel_2 &\leq& C \left ( 
\parallel \nabla s \parallel_{\infty} \ \parallel \omega^{\ell} s'
\parallel_2 \ + \  \parallel \omega^{\ell}s\parallel_2 \ \parallel 
\nabla s' \parallel_{\infty}\right )\nn \\
&\leq& C \left ( |s|_{k+1}^{\dot{}}\ |s'|_{\ell}^{\dot{}} + 
|s|_{\ell}^{\dot{}}\ |s'|_{k+1}^{\dot{}} \right )
\eea

\noi by a direct application of Lemma 3.2. The last two norms in the 
RHS of (\ref{4.13e}) are estimated as
\beq
\label{4.15e}
\parallel \omega^{\ell-1} |w|^2 \parallel_2  \ \leq \ C \parallel w 
\parallel_{\infty} \ \parallel \omega^{\ell-1}w\parallel_2 \ \leq \
C |w|_k\ |w|_{\ell - 1}\eeq

\noi by Lemma 3.2, and
$$\parallel \omega^{\ell+1} (x\cdot B_a)_L \parallel_2 \ \leq \ 
t^{\beta (\ell - k)} \parallel \omega^{k+1} (x\cdot B_a) \parallel_2$$

\noi by (\ref{3.20e}). Together with the simpler estimate
\beq
\label{4.16e}
\left  | \partial_t\parallel \nabla s'\parallel_2 \right | \ \leq \ 
t^{-2} \parallel \nabla s\parallel_{\infty} \ \parallel \nabla s'
\parallel_2 \ + C\ t^{-1} \parallel w\parallel_4^2 \ + \ t^{-1} \parallel
\nabla^2 (x \cdot B_a) \parallel_2 \ , \eeq

\noi those estimates yield (\ref{4.6e}). \par

We finally estimate $B'_b$. From (\ref{3.25e}) with $j = 1$ and 
(\ref{3.28e}), we obtain
\beq \label{4.17e}
|B'_b|_{k+\theta}^{\dot{}} \ \leq \ C\ t^{-1}\ I_1 \left ( 
|M_b|_{k+\theta - 1} \right )
\eeq
\beq \label{4.18e}
|x \cdot B'_b|_{k+\theta}^{\dot{}}  \ \leq \ C\ t^{-1}\ I_0 \left ( 
|xM_b|_{k+\theta - 1} + |M_b|_{k+\theta - 1}\right )
\eeq

\noi from which (\ref{4.7e}) (\ref{4.8e}) follow by (\ref{3.9e}) 
(\ref{3.10e}). \par\nobreak \hfill $\sq$\par

We next estimate the difference of two solutions of the linearized 
system (\ref{2.37e}) (\ref{2.38e}) corresponding to two different
choices of $(w, s, B_b)$. We estimate that difference at a lower 
level of regularity than the solutions themselves. \\

\noi {\bf Lemma 4.2.} {\it Let $k > 3/2$, let $T \geq 1$, $I = 
[T,\infty )$ and let $(w_i, s_i, B_{bi}) \in {\cal C} (I, X^k)$, 
$i=1,2$
with $|w_i|_k \vee |xw_i|_k \in L^{\infty}(I)$. Let
\beq
\label{4.19e}
a = \ \mathrel{\mathop {\rm Max}_{i}}\ \parallel |w_i|_k \vee |xw_i|_k ; L^{\infty}(I)\parallel \ .\eeq\noi Let $I' \subset I$ be an inter 
val, let $(w'_i, s'_i)$ be solutions of the system (\ref{2.37e}) 
associated with $(w_i, s_i,
B_{bi})$ with $(w'_i, s'_i,0) \in {\cal C}(I',X^k)$ and let $B'_{bi}$ 
be defined by (\ref{2.38e}) in terms of $(w_i, s_i, B_{bi})$.
Define $(w_{\pm}, s_{\pm},  B_{b_{\pm}}) = (1/2) ((w_1, s_1, B_{b1}) 
\pm (w_2, s_2, B_{b2}))$ and similarly for the primed quantities
and for $B_a$, $B$, $B\cdot s$ and $B^2$. Let
\beq
\label{4.20e}
1/2 < k' \leq k-1 \quad , \quad 0 \leq \theta \leq 1 \quad , \quad 1 
\vee k' \leq \ell \leq k'+2 \ .
\eeq

\noi Then the following estimates hold~:
\beq
\label{4.21e}
  |B_{a_-}|_{k'+1}^{\dot{}} \ \leq \ C \ I_0 \left ( |xw_+|_k \ 
|w_-|_{k'} \right ) \ ,
\eeq
\beq
\label{4.22e}
\parallel x\cdot B_{a_-}; \dot{H}^{k'+1} \parallel \ \leq \ C\ 
I_{k'-1} \left ( (|xw_+|_k
+ (|w_+|_k)  |xw_-|_{k'}\right )\ , \eeq
\bea
\label{4.23e}
&&\left  | \partial_t| w'_-|_{k'+\theta}  \right | \ \leq \ C\Big \{ 
t^{-2} \left ( |\nabla s_+|_k^{\dot{}} + |\nabla B_+|_k^{\dot{}} +
|(Bs)_+|_k^{\dot{}} + |(B^2)_+|_k^{\dot{}} \right ) \nn \\
&&+ t^{-1-\beta (k-k')} \parallel x \cdot 
B_{a_+};\dot{H}^{k+1}\parallel \ + t^{-1} |x \cdot 
B_{b_+}|_k^{\dot{}} \Big \} |w'_-|_{k'+
\theta} \nn \\
&&+ C \ t^{-2} \Big\{ \left ( |s_-|_{k'+1}^{\dot{}} + 
|B_-|_{k'+1}^{\dot{}} \right ) |w'_+|_{k+\theta} + \Big ( |s_-|_{k'+1+
\theta}^{\dot{}}\nn \\
&&+ |B_+|_k^{\dot{}} \left ( |s_-|_{k'+1}^{\dot{}} + 
|B_-|_{k'+1}^{\dot{}} \right ) + |s_+|_k^{\dot{}}\ 
|B_-|_{k'+1}^{\dot{}} \Big )
|w'_+|_k \Big \}\nn \\
&&+ C \Big \{ t^{-1-\beta (1 - \theta)} \parallel x \cdot B_{a_-}; 
\dot{H}^{k'+1} \parallel \ + t^{-1} |x\cdot B_{b_-}|_{k'+1}^{\dot{}}
\Big\} |w'_+|_k \ , \eea
\bea
\label{4.24e}
&&\partial_t| xw'_-|_{k'}  \ \leq \ {\rm Idem} (\theta = 0, w' \to xw')\nn \\
&&+ C \ t^{-2}\left \{ |w'_-|_{k'+1} + |(s+B)_+|_k^{\dot{}}\ 
|w_-|_{k'} + |(s+B)_-|_{k'+1}^{\dot{}}\  |w'_+|_k \right \} \ ,
\eea
\bea
\label{4.25e}
&&\partial_t| s'_-|_{\ell}^{\dot{}}  \ \leq \ C \ t^{-2}\left \{ 
|\nabla s_+|_{k}^{\dot{}}\ |s'_-|_{\ell}^{\dot{}} +
|s_-|_{\ell}^{\dot{}} \ |\nabla s'_+|_{k}^{\dot{}} + \chi (\ell \geq 
k) |s_-|_k^{\dot{}} \ |\nabla s'_+|_{\ell}^{\dot{}} \right \}
\nn \\
&&+ C\ t^{-1} |w_+|_k \ |w_-|_{\ell-1} + C\ t^{-1+(\ell - k')\beta} 
\parallel \nabla (x \cdot B_{a_-});\dot{H}^{k'} \cap
\dot{H}^{k'\wedge 1} \parallel \ ,\eea
\beq
\label{4.26e}
  |B'_{b_-}|_{k'+1}^{\dot{}} \leq \ t^{-1} \ I_1 \left ( 
|M_{b_-}|_{k'} \right ) \ ,
\eeq
\beq
\label{4.27e}
  |x \cdot B'_{b_-}|_{k'+1}^{\dot{}} \leq \ t^{-1}  \ I_0 \left ( 
|xM_{b_-}|_{k'} + |M_{b_-}|_{k'} \right ) \ ,
\eeq

\noi where $M_{b_-} = (M_{b1} - M_{b2})/2$, $M_b$ is defined by 
(\ref{2.23e}), and}
\bea
\label{4.28e}
&&|M_{b_-}|_{k'} \vee |x\cdot M_{b_-}|_{k'} \ \leq \ C |<x>w_+|_k 
\Big \{ |w_-|_{k'+1} + |s_+ + B_+|_k^{\dot{}}\ |w_-|_{k'}\nn\\
&&+ |s_-+B_-|_{k'+1}^{\dot{}} \ |w_+|_k \Big \} + C |<x>w_-|_k\ 
|s_-+B_-|_{k'+1}^{\dot{}} \ |w_-|_k \ .
\eea
\vskip 3 truemm

\noi {\bf Proof.} The estimates (\ref{4.21e}) (\ref{4.22e}) follow 
immediately from (\ref{2.28e}), from (\ref{3.25e}) (\ref{3.28e}) and
from (\ref{3.10e}) with $m = k'$ and $\bar{m} = k$. \par

We next estimate $(w'_-, s'_-)$. Taking the difference of the systems 
(\ref{2.37e}) for $(w'_i, s'_i)$, we obtain the following system
for $(w'_-, s'_-)$~:
\bea
\label{4.29e}
&&\partial_t w'_- = i (2t^2)^{-1} \Delta w'_- + t^{-2} \left ( Q(s_+ 
+ B_+, w'_-) + Q(s_- + B_-, w'_+) \right )\nn \\
&&- i(2t^2)^{-1} \left ( (2(B \cdot s)_+ + (B^2)_+) w'_- + ( 2(B\cdot 
s)_- + (B^2)_-) w'_+ \right )\nn \\
&&+ i t^{-1} \Big ( ((x \cdot B_{a_+})_S + x \cdot B_{b_+} )w'_- + 
((x \cdot B_{a_-})_S + x \cdot B_{b_-}) w'_+ \Big )
\eea
\beq
\label{4.30e}
\partial_t s'_- = t^{-2}\left ( s_+ \cdot \nabla s'_- + s_- \cdot 
\nabla s'_+\right ) + t^{-1}\nabla g (w_+, w_-) - t^{-1}\nabla (x
\cdot B_{a_-})_L \ . \eeq

\noi Let $m = k' + \theta$, so that $1/2 < m \leq k$. We estimate 
$\omega^m w'_-$ by the same method as in Lemma 4.1, using
(\ref{3.14e}) for the term $(s+ B)_+\cdot \nabla w'_-$ with $\bar{m} = k$, and using (\ref{3.8e}) for all the
other terms, with $\bar{m} = 
k$ or $\bar{m} = k'+1$, depending on whether the estimated quantity 
is of $+$ or $-$ type. We obtain

$$\left  | \partial_t\parallel \omega^m w'_-\parallel_2 \right | \ 
\leq \ C \Big \{ t^{-2} \left ( | \nabla (s+ B)_+ |_k^{\dot{}} +
|\nabla \cdot s_+|_k^{\dot{}} + |(B\cdot s)_+ + (B^2)_+ |_k^{\dot{}} \right )$$
$$+ t^{-1} \left ( |(x\cdot B_{a_+})_S|_{k'+1}^{\dot{}} + |x \cdot 
B_{b_+}|_{k}^{\dot{}} \right ) \Big \} |w'_-|_m$$
$$+ C \Big \{ t^{-2} \left ( | (s+ B)_- |_{k'+1}^{\dot{}} \ 
|w'_+|_{m+1} + |\nabla \cdot s_-|_m \ |w'_+|_k + |(B\cdot s)_- +
(B^2)_-|_{k'+1}^{\dot{}} \ |w'_+|_m \right )$$
\beq
\label{4.31e}
+t^{-1} \left ( |(x \cdot B_{a_-})_S|_m\ |w'_+|_k + |x \cdot 
B_{b_-}|_{k'+1}^{\dot{}} \ |w'_+|_m \right ) \Big \} \ .
\eeq

We next estimate
$$|(x \cdot B_{a_+})_S|_{k'+1}^{\dot{}} \ \leq \ t^{-\beta (k-k')} 
\parallel x \cdot B_{a_+};\dot{H}^{k+1} \parallel \ ,$$
$$|(x \cdot B_{a_-})_S|_{m} \ \leq \ t^{-\beta (1-\theta)} \parallel 
x \cdot B_{a_-};\dot{H}^{k'+1} \parallel $$

\noi by (\ref{3.19e}). We estimate $\partial_t\parallel 
w'_-\parallel_2$ in a simpler way as above by taking $m = 0$ and 
omitting the
terms with $w'_-$ in (\ref{4.31e}). Collecting the previous estimates 
yields (\ref{4.23e}). \par

We now turn to the proof of (\ref{4.24e}), namely to the estimate of 
$xw'_-$. From (\ref{2.39e}) we obtain
\beq
\label{4.32e}
\partial_t xw'_- = \ {\rm Idem}(w'\to xw') - it^{-2} \nabla w'_- - 
t^{-2} \left ( (s+B)_+ w'_- + (s+B)_- w'_+ \right )
\eeq

\noi where Idem denotes the RHS of (\ref{4.29e}). The estimate 
(\ref{4.24e}) then follows from (\ref{4.23e}) with $\theta = 0$ and 
$w'$
replaced by $xw'$ and from an easy estimate of the additional terms 
using (\ref{3.10e}) with $m = k'$ and $\bar{m} = k$ or
$k'+1$.\par

We now turn to the proof of (\ref{4.25e}) namely to the estimate of 
$s'_-$. We estimate
\bea
\label{4.33e}
&& \partial_t \parallel \omega^{\ell}s'_- \parallel_2 \ \leq \ t^{-2} 
\Big \{ \parallel [\omega^{\ell}, s_+]\cdot \nabla s'_-
\parallel_{2} \ + \ \parallel \nabla \cdot s_+ \parallel_{\infty} \ 
\parallel \omega^{\ell}s'_-\parallel_2 \ \nn \\
&&+ \ \parallel
\omega^{\ell}(s_-\nabla s'_+)\parallel_2 \Big \} + C\ t^{-1} 
\parallel \omega^{\ell - 1} (\bar{w}_+ w_-) \parallel_2 \ + t^{-1}
\parallel \omega^{\ell + 1} (x\cdot B_{a_-})_L \parallel_2 \ . \nn \\ \eea

\noi We estimate the first norm in the RHS by (\ref{3.11e}) 
(\ref{3.12e}) with $m = \ell$ and $\bar{m} = k$. We estimate the third
norm by (\ref{3.8e}) with $m = \ell$ and $\bar{m} = k$ for $\ell \leq 
k$, while for $\ell \geq k$
\begin{eqnarray*}
&&\parallel \omega^{\ell} (s_- \nabla s'_+)\parallel_2 \ \leq \ 
C\left (\parallel \omega^{\ell} s_- \parallel_2 \ \parallel \nabla
s'_+\parallel_{\infty} \ + \ \parallel  s_-\parallel_{\infty} \ 
\parallel \omega^{\ell} \nabla s_+
\parallel_{2} \right )\\
&&\leq C \left ( |s_-|_{\ell}^{\dot{}}\ |\nabla s'_+|_k^{\dot{}} + 
|s_-|_k^{\dot{}}\ |\nabla s'_+|_{\ell}^{\dot{}} \right )
\end{eqnarray*}

\noi by a direct application of Lemma 3.2. We next estimate
$$\parallel \omega^{\ell-1} (\bar{w}_+ w_-)\parallel_2 \ \leq \ C 
|w_+|_k \ |w_-|_{\ell-1}$$

\noi by (\ref{3.10e}) with $m = \ell - 1$ and $\bar{m} = k$, and we 
estimate the last norm in (\ref{4.33e}) by (\ref{3.20e}). Together
with the special case $\ell = 1$, the previous estimates yield 
(\ref{4.25e}). \par

Finally (\ref{4.26e}) (\ref{4.27e}) are special cases of 
(\ref{3.25e}) (\ref{3.28e}), while (\ref{4.28e}) follows from 
repeated use of
(\ref{3.10e}) with $m = k'$ and $\bar{m} = k$ or $k' + 1$. 
\par\nobreak \hfill $\sq$\par

We now begin to study the auxiliary system (\ref{2.34e}) 
(\ref{2.35e}) and its linearized version (\ref{2.37e}) (\ref{2.38e}). 
The
first step is to solve the linear system (\ref{2.37e}) globally in time.\\
 
\noi {\bf Proposition 4.1.} {\it Let $k > 3/2$, let $T \geq 1$, $I = 
[T,\infty )$ and let $(w, s, B_{b}) \in {\cal C} (I, X^k)$
with $|w|_k \vee |xw|_k \in L^{\infty}(I)$. Let $t_0 \in I$ and 
$(w'_0, s'_0, 0) \in X^k$. Then the system (\ref{2.37e}) has a unique
solution $(w', s')$ in $I$ such that $(w', s', 0) \in {\cal C} (I, 
X^k)$ and $(w', s')(t_0) = (w'_0, s'_0)$. That solution satisfies the
estimates (\ref{4.4e}) (\ref{4.5e}) (\ref{4.6e}) for all $t \in I$. 
Two such solutions $(w'_i, s'_i)$ associated with $(w_i, s_i)$, $i =
1,2$, satisfy the estimates (\ref{4.23e}) (\ref{4.24e}) (\ref{4.25e}) 
for all $t \in I$.}\\

\noi {\bf Proof.}  We first prove the existence of a unique solution 
$(w', s') \in {\cal C} (I, Y^k)$ where $Y^k = H^{k+1} \oplus
K^{k+2}$. The proof proceeds in the same way as that of Proposition 
4.1 of \cite{6r}, through a parabolic regularization and a limiting
procedure. We define $U_1(t) = U(1/t)$ and $\widetilde{w}'(t) = 
U_1(t) w'(t)$. We first consider the case $t \geq t_0$. The system
(\ref{2.37e}) with a parabolic regularization added is rewritten in 
terms of the variables $(\widetilde{w}', s')$ as
\begin{eqnarray*}
&&\left \{ \begin{array}{l} \partial_t \widetilde{w}' = \eta \Delta 
\widetilde{w}' + U_1 \Big ( L - i (2t^2)^{-1} \Delta
\Big ) U_1^* \widetilde{w}' \equiv \eta \Delta \widetilde{w}' + 
G_1(\widetilde{w}')\\ \\ \partial_t s' = \eta \Delta s' +
t^{-2} s \cdot \nabla s' + t^{-1} \nabla g(w) - t^{-1} \nabla (x 
\cdot B_a)_L \\\end{array} \right . \\ && \\ &&\qquad \quad \equiv 
\eta
\Delta s' + G_2 (s') \end{eqnarray*}

\noi where $L$ is defined in (\ref{2.34e}) and where the parametric 
dependence of $L$, $G_1$, $G_2$ on $(w, s, B_b)$ has been omitted.
The Cauchy problem for that system can be recast into the integral form
\beq
\label{4.34e}
{\widetilde{w}' \choose s'} (t) = V_{\eta}(t-t_0) {\widetilde{w}'_0 
\choose s'_0} +
\int_{t_0}^t dt' \ V_{\eta} (t-t') {G_1(\widetilde{w}') \choose 
G_2(s')} (t') \eeq

\noi where $V_{\eta} (t) = \exp (\eta t \Delta)$. The operator 
$V_{\eta}(t)$ is a contraction
in $Y^{k}$ and satisfies the bound
$$\parallel \nabla V_{\eta}(t) ; {\cal L} (Y^{k}) \parallel \ \leq 
C(\eta t)^{-1/2} \ .$$

\noi From those facts and from estimates on $G_1$, $G_2$ similar to 
and mostly contained in those
of Lemma 4.1, it follows by a contraction argument that the system 
(\ref{4.34e}) has a unique
solution $(\widetilde{w}'_{\eta}, s'_{\eta}) \in {\cal C}([t_0, t_0 + 
T_1], Y^{k})$ for
some $T_1 > 0$ depending only on $|w'_0|_{k+1}$, 
$|s'_0|_{k+2}^{\dot{}}$ and $\eta$. That solution
satisfies the estimates (\ref{4.4e}) and (\ref{4.6e}) and can 
therefore be extended to $[t_0, \infty )$ by a standard globalisation
argument using Gronwall's inequality. \par

We next take the limit $\eta \to 0$. Let $\eta_1$, $\eta_2 > 0$ and 
let $(w'_i, s'_i) =
(w'_{\eta_i}, s'_{\eta_i})$, $i = 1,2$ be the corresponding 
solutions. Let $(w'_-, s'_-) =
(1/2)(w'_1 - w'_2, s'_1 - s'_2)$. By estimates similar to, but 
simpler than those of Lemma 4.2,
since in particular $(w_-, s_-, B_{b_-}) = 0$, we obtain
$$\left \{ \begin{array}{l} \partial_t \parallel w'_-\parallel_2^2 \ \leq
|\eta_1 - \eta_2 | \left ( \parallel \nabla w'_1 \parallel_2^2 + \parallel
\nabla w'_2 \parallel_2^2 \right ) \\
\\
\partial_t \parallel \nabla s'_-\parallel_2^2 \ \leq
|\eta_1 - \eta_2 | \left ( \parallel \nabla^2 s'_1 \parallel_2^2 + \parallel
\nabla^2 s'_2 \parallel_2^2\right ) + C\ t^{-2} \parallel \nabla s 
\parallel_{\infty} \parallel
\nabla s'_-\parallel_2^2 \ .\end{array}\right .$$
 
\noi Those estimates imply that $(w'_{\eta}, s'_{\eta})$ converges in 
$L^2 \oplus \dot{H}^1$ uniformly in
time in the compact subintervals of $[t_0, \infty)$, to a solution of 
the original system. It follows
then by a standard compactness argument using the estimates 
(\ref{4.4e}) (\ref{4.6e}) that the
limit belongs to ${\cal C}([t_0, \infty ), Y^{k})$. This completes 
the proof for $t \geq t_0$. The
case $t \leq t_0$ is treated similarly. \par

We next show that $xw' \in {\cal C}(I, H^k)$. For that purpose we 
choose a function $\psi \in {\cal C}_0^{\infty} ({I\hskip-1truemm
R}^n, {I\hskip-1truemm R}^+)$ with $0 \leq \psi \leq 1$, $\psi (x) = 
1$ for $|x| \leq 1$, $\psi (x) = 0$ for $x \geq 2$, and we
define $\psi_R$ by $\psi_R(x) = \psi (x/R)$. Clearly $\psi_R \ xw' 
\in {\cal C}(I, H^{k+1})$ and $\psi_R \ xw'$ satisfies the
equation      \bea \label{4.35e}
&&\partial_t \ \psi_R \ xw' = L \ \psi_R\ xw' - i t^{-2} \left ( 
\nabla (\psi_R x)\right ) \cdot \nabla w' - i(2t^2)^{-1} \left ( 
\Delta
(\psi_R x) \right ) w'\nn\\
&&\hskip 3 truecm - t^{-2} (s + B) \cdot \left ( \nabla (\psi_R x) 
\right ) w' \ . \eea

\noi Using Lemma 4.1, more precisely (\ref{4.5e}) and the fact that 
the operator of multiplication by a function $\varphi_R(x) = \varphi
(x/R)$ for $\varphi \in {\cal C}_0^{\infty}$ is a bounded operator in 
$H^m$ for all $m \geq 0$ uniformly in $R$ for $R \geq 1$, we
estimate
\bea \label{4.36e}
&&\left | \partial_t | \psi_R \ xw' |_k \right | \leq M_4 (0, \psi_R\ 
xw') + C \ t^{-2} \left ( |w'|_{k+1} + R^{-1} |w'|_k \right )\nn\\
&&\hskip 3 truecm + C\  t^{-2} |s + B|_k^{\dot{}} \ |w'|_k  \ . \eea

\noi Integrating (\ref{4.36e}) between $t_0$ and $t$ and using 
Gronwall's inequality, we obtain
$$| \psi_R \ xw' |_k  \leq  C(t) \left ( 1 + |\psi_R \ x w'_0|_k 
\right ) \leq  C(t) \left ( 1 + |xw'_0|_{k}\right )$$

\noi so that $\psi_R \ xw'$ is bounded in $H^k$ uniformly in $R$. 
This implies that $x w' \in L_{loc}^{\infty}(I, L^2)$ and that $\psi_R
\ xw'$ tends to $xw'$ strongly in $L^2$ pointwise in $t$ when $R \to 
\infty$. Moreover, it follows from (\ref{4.35e}) that $xw'$ is weakly
continuous in $L^2$ as a function of $t$. Together with 
(\ref{4.36e}), this implies that $xw' \in  {\cal C}(I, H^{k})$ by 
standard
compactness arguments.\par \nobreak \hfill
$\sq$\par

We now turn to the study of the auxiliary system (\ref{2.34e}) 
(\ref{2.35e}). The main results will be the existence and uniqueness 
of
solutions of that system, defined in a neighborhood of infinity in 
time and with suitable bounds at infinity, and some asymptotic
behaviour of those solutions. The bounds on the solutions at infinity 
will be essentially dictated by the existence result (Proposition
4.3 below), and for simplicity we shall mostly restrict our attention 
to solutions satisfying those bounds, although more general
solutions could be considered in the uniqueness and asymptotic 
behaviour results. We shall thus consider solutions $(w, s, B_b) \in 
{\cal
C}(I, X^{k})$ for some interval $I = [T, \infty )$ such that  \bea
\label{4.37e}
\Lambda &\equiv& \Lambda (t) \equiv |w|_k \vee |xw|_k \vee (\ell n \ 
t)^{-1} |w|_{k+1} \vee (\ell n \ t)^{-1} |s|_k^{\dot{}} \vee
t^{-\beta} \ |s|_{k+1}^{\dot{}} \nn \\
&&\vee t^{-2 \beta} \ |s|_{k+2}^{\dot{}} \vee |B_b|_k^{\dot{}} \in 
L^{\infty}(I) \ . \eea

\noi We first state the uniqueness result.\\

\noi {\bf Proposition 4.2.} {\it Let $k > 3/2$, $0 < \beta < 1/2$ and 
$I = [T,\infty )$. \par

(1) Let $t_0 \in I$ and $(w_0, s_0, 0) \in X^k$. Then for $t_0$ 
sufficiently large, the auxiliary system (\ref{2.34e}) (\ref{2.35e}) 
has
at most one solution $(w, s, B_b) \in {\cal C}(I, X^{k})$ satisfying 
(\ref{4.37e}) and $(w,s)(t_0) = (w_0,s_0)$.\par

(2) Let $(w_i, s_i, B_{bi})$ $i = 1,2$, be two solutions of the 
auxiliary system (\ref{2.34e}) (\ref{2.35e}) in ${\cal C}(I, X^{k})$
satisfying (\ref{4.37e})  and such that for some $\varepsilon > 0$
\beq
\label{4.38e}
|s_-|_{k+1}^{\dot{}} \vee t^{2\beta + \varepsilon} \left ( |w_-|_k 
\vee |xw_-|_{k-1} \right ) \to 0 \ {\it when} \ t \to \infty \ .
\eeq

\noi Then $(w_1, s_1, B_{b1}) = (w_2, s_2, B_{b2})$}.\\

\noi {\bf Remark 4.2.} The condition $t_0$ sufficiently large in Part 
(1) takes the form
\beq \label{4.39e}
t_0 \geq \ \parallel 1 + \Lambda ; L^{\infty}([t_0, \infty )) \parallel^N
\eeq

\noi for some $N >0$. For a given solution $(w, s, B_b)$, the RHS of 
(\ref{4.39e}) is decreasing in $t_0$ while the LHS is increasing,
and Part (1) supplemented by (\ref{4.39e}) gives a lower bound of the 
initial time for which the given solution is uniquely determined
as a solution of the Cauchy problem with that initial time.\\

\noi {\bf Proof.} The proof relies on Lemma 4.2, and we first recast 
the estimates of that lemma in a simplified form for solutions
satisfying (\ref{4.37e}). We consider two solutions $(w_i, s_i, 
B_{bi})$ $(\equiv (w'_i, s'_i, B'_{bi}))$, $i = 1,2$ of the system
(\ref{2.34e}) (\ref{2.35e}) satisfying (\ref{4.37e}) and we define
\beq \label{4.40e}
A(t_1) = \ \mathrel{\mathop {\rm Max}_{i=1,2}} \ \parallel 
\Lambda_i(t) ; L^{\infty}([t_1, \infty )) \parallel
\eeq

\noi for $t_1 \geq T$, where $\Lambda_i$ are the quantities defined 
by  (\ref{4.37e}) for the two
solutions. We define $(w_-, s_-, B_{b_-})$ as in Lemma 4.2, and
  \beq
\label{4.41e}
\left \{ \begin{array}{l} y = |w_-|_{k-1} \vee |xw_-|_{k-1} \quad , 
\qquad y_1 = |w_-|_k \ , \\ \\ z_j = |s_-|_{k+j-1}^{\dot{}} \quad ,
\qquad j = 1,2 \ .\end{array} \right . \eeq

We rewrite the estimates (\ref{4.21e})-(\ref{4.28e}) for general $t 
\in [t_1, \infty)$ for some $t_1 \geq T$. The $+$ quantities are
estimated by (\ref{4.37e}) (\ref{4.40e}) supplemented by (\ref{4.2e}) 
(\ref{4.3e}) (\ref{4.7e}) (\ref{4.8e}) as regards $B_a$ and
$B_b$. This produces overall constants depending polynomially on 
$A(t_1)$, which we omit for brevity, but the occurrence of which 
should
be kept in mind for subsequent arguments. In terms having the same 
dependence on the dynamical variables, we keep only the terms with
the leading behaviour in $t$, namely we use
$$t^{1-\beta} \geq \ell n \ t \geq 1 \ .$$

\noi We take $k' = k-1$ and rewrite (\ref{4.21e})-(\ref{4.28e})  as follows
\beq
\label{4.42e}
|B_{a_-}|_k^{\dot{}} \leq I_0 (y) \ ,
\eeq
\beq
\label{4.43e}
\parallel x \cdot B_{a_-}; \dot{H}^k \parallel \ \leq I_{k-2}(y) \ ,
\eeq
\beq
\label{4.44e}
|\partial_t y| \leq t^{-1-\beta} y + t^{-2} \left ( z_1 + 
|B_-|_k^{\dot{}} \ \ell n \ t\right ) + t^{-1-\beta} \ I_{k-2}(y) \nn 
\\
  + t^{-1} |x\cdot B_{b_-}|_k^{\dot{}} + t^{-2} y_1 \ ,
\eeq
\beq
\label{4.45e}
|\partial_t y_1| \leq t^{-1-\beta} y_1 + t^{-2} \left ( z_1\  \ell n 
\ t + z_2 + |B_-|_k^{\dot{}} \ \ell n \ t\right ) + t^{-1}
I_{k-2}(y) + t^{-1} |x\cdot B_{b_-}|_k^{\dot{}} \ ,
\eeq
\beq
\label{4.46e}
|\partial_t z_1| \leq t^{-2+\beta} z_1 + t^{-1} y + t^{-1+\beta} \ I_m(y) \ ,
\eeq
\beq
\label{4.47e}
|\partial_t z_2| \leq t^{-2+\beta} \left ( z_2 + t^{\beta} z_1 \right 
) + t^{-1} y_1  + t^{-1+2\beta} \ I_m(y)
\eeq

\noi with $m = (k-2) \wedge 0$,
\beq
\label{4.48e}
|B_{b_-}|_k^{\dot{}}  \leq t^{-1}\ I_1 \left ( y_1 + y \ \ell n \ t + 
z_1 + |B_-|_k^{\dot{}}  \right ) \ ,
\eeq
\beq
\label{4.49e}
|x \cdot B_{b_-}|_k^{\dot{}}  \leq t^{-1}\ I_0 \left ( y_1 + y \ \ell 
n \ t + z_1 + |B_-|_k^{\dot{}}  \right ) \ .
\eeq

The system (\ref{4.42e})-(\ref{4.49e}) will be the starting point for 
the proof of Proposition 4.2.\\

\noi \underbar{Part (1)}. We first prove uniqueness for $t \geq t_0$. 
Note that this region is autonomous in the
sense that the equations in this region involve the dynamical 
variables in this region only, since the integrals (\ref{2.25e}) 
occurring
in (\ref{2.28e})-(\ref{2.29e}) are taken for $\nu \geq 1$. We define
  \beq \label{4.50e}
Y = \ \parallel y; L^{\infty}([t_0, \infty )) \parallel \quad , \quad 
Y_1 = \ \parallel y_1(\ell n \ t)^{-1}; L^{\infty}([t_0, \infty )) 
\parallel
\eeq

\noi which are finite by (\ref{4.37e}) and we prove that those 
quantities are zero by integrating (\ref{4.44e})-(\ref{4.47e}) with
initial condition $(y,y_1,z_1,z_2)(t_0) = 0$.\par

In keeping with the previous simplification, we perform that 
computation up to constants (depending on $A(t_0))$ and under 
conditions
that $t_0$ is sufficiently large in the sense of (\ref{4.39e}) when 
needed. We furthermore eliminate the diagonal terms in
(\ref{4.46e})-(\ref{4.47e}) by exponentiation, namely by using the fact that
\beq
\label{4.51e}
\partial_t y \leq fy + g \Rightarrow y(t) \leq \exp \left ( 
\int_{t_0}^{\infty} f\right ) \left ( y(t_0) + \int_{t_0}^t g \right )
\eeq

\noi for integrable $f$. Using (\ref{4.50e}) and integrating 
(\ref{4.46e}) (\ref{4.47e}) (and using $z_1 \leq z_2$ and $2\beta < 
1$), we
obtain
  \beq
\label{4.52e}
z_1 \leq t^{\beta} Y \qquad , \qquad z_2 \leq (\ell n \ t)^2 Y_1 + 
t^{2\beta } Y \ .
\eeq

\noi Substituting (\ref{4.50e}) (\ref{4.52e}) into (\ref{4.42e}) 
(\ref{4.48e}) (\ref{4.49e}) yields
\beq
\label{4.53e}
|B_{a_-}|_k^{\dot{}}  \ \vee \parallel x \cdot B_{a_-} ; \dot{H}^k 
\parallel \ \leq Y \ ,
\eeq
\bea \label{4.54e}
|B_{b_-}|_k^{\dot{}} \vee |x \cdot B_{b_-}|_k^{\dot{}} &\leq& 
t^{-1}\ell n \ t \ (Y_1 + Y) + t^{-1+\beta} \ Y + t^{-1} \parallel
|B_{-}|_k^{\dot{}} ; L^{\infty}([t_0, \infty )) \parallel \nn \\
&\leq& t^{-1} \ell n \ t \ Y_1 + t^{-1+\beta} \ Y
  \eea

\noi for $t_0$ sufficiently large, so that one can assume
\beq
\label{4.55e}
|B_{-}|_k^{\dot{}} \leq Y \ .
\eeq

\noi Substituting (\ref{4.50e}) (\ref{4.52e})-(\ref{4.55e}) into 
(\ref{4.44e}) (\ref{4.45e}) yields
\bea
\label{4.56e}
\partial_t y &\leq& \left ( t^{-1-\beta} + t^{-2+ \beta} + t^{-2} 
\ell n \ t \right ) Y + t^{-2} \ell n\ t \ Y_1 \nn \\
&\leq & t^{-1-\beta} \ Y + t^{-2} \ell n\ t \ Y_1  \ ,
\eea
\bea
\label{4.57e}
\partial_t y_1 &\leq& \left ( t^{-2+\beta} \ell n \ t + t^{-2} \ell n 
\ t + t^{-2+2\beta} + t^{-1}\right ) Y + t^{-2} (\ell n\ t)^2 \ Y_1
\nn \\ &\leq & t^{-1} \ Y + t^{-2} (\ell n\ t)^2 \ Y_1  \ .
\eea

\noi Integrating (\ref{4.56e}) (\ref{4.57e}) in $[t_0, t]$ and 
comparing with (\ref{4.50e}) yields
\beq
\label{4.58e}
\left \{ \begin{array}{l} Y \leq t_0^{-\beta} \ Y + t_0^{-1} \ell n\ 
t_0\ Y_1\\ \\ Y_1 \leq Y + t_0^{-1} \ell n\ t_0\ Y_1 \\
\end{array} \right . \eeq

\noi which implies $Y = Y_1 = 0$ for $t_0$ sufficiently large. \par

We next prove uniqueness for $t \leq t_0$. Since we have already 
proved uniqueness for $t \geq t_0$, $(y, y_1, z_1, z_2)$ vanish for $t
\geq t_0$, and the region $t \leq t_0$ also becomes autonomous in the 
previous sense, which it would not have been if treated first. Now
however we are in a standard situation where the variables $\{y_i\} = 
\{y, y_1, z_1, z_2, I_0 (|B_{-}|_k^{\dot{}}\}$ satisfy a system of
inequalities of the type
\beq
\label{4.59e}
|\partial_t y_i| \leq \sum_j f_{ij}\ y_j + g_i \int_t^{t_0} dt' 
\sum_j h_{ij}(t')\ y_j(t')
\eeq

\noi for $t \leq t_0$. This can be reduced to the case of a single 
function $\bar{y} = \sum y_i$ namely
\beq
\label{4.60e}
|\partial_t \bar{y}| \leq f\bar{y} + g \int_t^{t_0} dt' \ h(t') \ \bar{y}(t')
\eeq

\noi with
$$f = \ \mathrel{\mathop {\rm Max}_{j}}\ \sum_i f_{ij} \quad , \quad 
g = \ \mathrel{\mathop {\rm Max}_{i}}\ g_i \quad , \quad h = \ 
\mathrel{\mathop {\rm
Max}_{j}} \  \sum_i h_{ij} \ .$$

\noi One can then eliminate $f$ by exponentiation, in very much the 
same way as in (\ref{4.51e}). Integrating (\ref{4.60e}) with $f = 0$
and $\bar{y}(t_0) = 0$ yields
$$\bar{y}(t) \leq G(t) \int_t^{t_0} dt'\ h(t')\ \bar{y}(t') \equiv 
G(t) \ \bar{Y}(t)$$

\noi with $G(t) = \int_t^{t_0} dt'\ g(t')$ and therefore
$$|\partial_t \bar{Y}| \leq h\ G\ \bar{Y}$$

\noi with $\bar{Y}(t_0) = 0$, which implies $\bar{Y} = 0$ and 
therefore $\bar{y} = 0$ for all $t \leq t_0$.\\

\noi \underbar{Part (2)}. We proceed in the same way as in Part (1) 
for $t \geq t_0$. Let $\lambda = 2 \beta +
\varepsilon$, so that $2\beta < \lambda < 1$ and define
\beq
\label{4.61e}
\left \{ \begin{array}{l} Y(t) = \parallel \cdot^{\lambda} y(\cdot ); 
L^{\infty}([t, \infty )) \parallel \\ \\ Y_1(t) = \parallel
\cdot^{\lambda} y_1(\cdot ); L^{\infty}([t, \infty )) \parallel \ .\\ 
\end{array} \right . \eeq

\noi We take $t_1$ and $t_0$ such that $T \leq t_1 \leq t_0 < 
\infty$, with $t_1$ sufficiently large, and we estimate $(y, y_1,
z_1,z_2)$ for $t \in [t_1, t_0]$ by integrating 
(\ref{4.44e})-(\ref{4.47e}) between $t$ and $t_0$ with final data 
$(y_0, y_{10}, z_{10},
z_{20})$ at $t_0$. Let $Y = Y(t_1)$, $Y_1 = Y_1(t_1)$. From 
(\ref{4.42e}) (\ref{4.43e}) (\ref{4.61e}) we obtain
\beq
\label{4.62e}
|B_{a_-}|_k^{\dot{}}  \vee \parallel x \cdot B_{a_-} ; \dot{H}^k 
\parallel \ \leq t^{-\lambda}\ Y
\eeq

\noi for $t \geq t_1$. Integrating (\ref{4.46e}) (\ref{4.47e}) as 
before, we now obtain
\beq
\label{4.63e}
\left \{ \begin{array}{l}  z_1 \leq z_{10} + t_-^{\beta - \lambda} \ 
Y\\ \\ z_2 \leq z_{20} + t_-^{- \lambda} \ Y_1 + \varepsilon^{-1}
\ t_-^{-\varepsilon}\ Y\\ \end{array} \right . \eeq

\noi for $t \geq t_1$, with $t_- = t \wedge t_0$. (Note that we need 
an estimate of $z_1$ for $t \geq t_0$ for substitution in
(\ref{4.48e}) (\ref{4.49e})). Substituting (\ref{4.61e}) 
(\ref{4.62e}) (\ref{4.63e}) into (\ref{4.48e}) (\ref{4.49e}) yields
\bea \label{4.64e}
|B_{b_-}|_k^{\dot{}} \vee |x \cdot B_{b_-}|_k^{\dot{}} &\leq& 
t^{-1}\left \{ t^{-\lambda}(Y_1 + Y \ell n \ t) + z_{10} + t_-^{\beta 
-
\lambda} Y + t^{-\lambda} Y + I_0\left (  |B_{b_-}|_k^{\dot{}}\right 
) \right \}  \nn \\
&\leq& t^{-1} \left ( t^{-\lambda} Y_1 + z_{10} + t^{\beta - \lambda} 
Y \right )
  \eea

\noi for $t_1 \leq t \leq t_0$ and $t_1$ sufficiently large. 
Substituting (\ref{4.61e})-(\ref{4.64e}) into (\ref{4.44e}) 
(\ref{4.45e})
yields
\bea
\label{4.65e}
|\partial_t y| &\leq& \left ( t^{-1-\beta- \lambda} + t^{-2+\beta - 
\lambda} \right ) Y + t^{-2} z_{10} + t^{-2-\lambda} \ Y_1\nn \\
&\leq& t^{-1-\beta - \lambda} Y + t^{-2} z_{10} + t^{-2-\lambda}  \ Y_1   \eea

\noi since $2\beta < 1$,
\bea
\label{4.66e}
|\partial_t y_1| &\leq& \left ( t^{-1-\beta- \lambda} + t^{-2- 
\lambda} \right ) Y_1 + t^{-2} \left ( z_{10} \ \ell n \ t + z_{20} 
\right
)\nn \\
&&+ \left ( t^{-2+\beta- \lambda} + \varepsilon^{-1} \ t^{-2- 
\varepsilon} + t^{-1- \lambda} \right ) Y  \nn \\
&\leq& t^{-1-\beta - \lambda} \ Y_1 + t^{-2} \left ( z_{10} \ \ell n 
\ t + z_{20} \right ) + t^{-1-\lambda}  \ Y \ .   \eea

\noi Integrating (\ref{4.65e}) (\ref{4.66e}) between $t$ and $t_0$ yields
\beq
\label{4.67e}
\left \{ \begin{array}{l} y \leq y_{0} + t^{-\beta - \lambda} \ Y + 
t^{-1} z_{10} + t^{-1 - \lambda} \ Y_1\\ \\ y_1 \leq y_{10} +
t^{-\beta - \lambda} \ Y_1 + t^{-1} \left ( z_{10} \ \ell n\ t + 
z_{20}\right ) + t^{-\lambda} \ Y \ .\\ \end{array} \right . \eeq

\noi Substituting the result into the definition (\ref{4.61e}) and 
using the fact that $t_0^{\lambda} \ y_0 \leq Y(t_0)$, $t_0^{\lambda}
\ y_{10} \leq Y_1(t_0)$ by definition, we obtain
\beq
\label{4.68e}
\left \{ \begin{array}{l} Y \leq Y(t_{0}) + t_1^{-1 + \lambda} z_{10} 
+ t_1^{-\beta} \ Y + t_1^{-1} \ Y_1\\ \\ Y_1 \leq Y_1(t_{0})
+ t_1^{-1+ \lambda} \left ( z_{10} \ \ell n\ t_1 + z_{20}\right ) + Y 
+ t_1^{-\beta} \ Y \ .\\ \end{array} \right . \eeq

\noi Taking $t_1$ sufficiently large to eliminate the diagonal terms, we obtain
$$\left \{ \begin{array}{l} Y \leq m_{0} + t_1^{-1}  \ Y_1\\ \\ Y_1 
\leq m_{1} +Y \\ \end{array} \right . $$

\noi with
\beq
\label{4.69e}
\left \{ \begin{array}{l} m_0 = Y(t_0) + t_1^{-1+ \lambda} z_{10} \\ 
\\ m_1 = Y_{1}(t_0)  + t_1^{-1+ \lambda} \left ( z_{10} \ \ell n \
t_1 + z_{20}\right )\\ \end{array} \right . \eeq

\noi which implies
$$Y \leq m_0 + t_1^{-1} m_1\qquad , \qquad Y_1 \leq m_0 + m_1$$

\noi for $t_1$ sufficiently large. We now let $t_0 \to \infty$ for 
fixed $t_1$. Then $m_0$ and $m_1$ tend to zero, and therefore $Y = Y_1
= 0$, which together with (\ref{4.63e}) (\ref{4.64e}) and another 
limit $t_0 \to \infty$ for fixed $t_1$, implies that $(w_1, s_1,
B_{b1}) = (w_2, s_2, B_{b2})$.\par\nobreak \hfill $\sq$\par

We now turn to the main result of this section, namely the existence 
of solutions of the auxiliary system (\ref{2.34e}) (\ref{2.35e})
with sufficiently large initial time, defined for sufficiently large 
time and satisfying (\ref{4.37e}). \\
 
\noi {\bf Proposition 4.3.} {\it Let $k > 3/2$ and $0 < \beta < 1/2$. 
Let $(w_0,s_0, 0) \in X^k$. Then \par
(1) There exists $T_0 < \infty$, depending on $(w_0, s_0)$, such that 
for all $t_0 > T_0$, there exists $T < t_0$ such that the
auxiliary system (\ref{2.34e}) (\ref{2.35e}) has a unique solution 
$(w, s, B_b) \in {\cal C}(I, X^{k})$, where $I = [T,\infty )$,
satisfying (\ref{4.37e}) and  $(w,s)(t_0) = (w_0,s_0)$. \par

The dependence of $T_0$ (resp. $T$) on $(w_0,s_0)$ (resp. and on 
$t_0$) can be formulated more precisely as follows. For any set ${\cal
A} = \{a, a_1, b_0, b_1, b_2\}$ of five positive numbers, there 
exists $T_0 = T_0 ({\cal A})$ and for any $t_0 > T_0$, there exists 
$T =
T(t_0,{\cal A}) < t_0$, increasing in $t_0$ and such that 
$T(T_0({\cal A}), {\cal A}) = T_0 ({\cal A})$, such that the previous 
statement
holds for such $T_0$, $t_0$, $T$ for all $(w_0,s_0,0) \in X^k$ such that
\beq
\label{4.70e}
\left \{ \begin{array}{l} |w_0|_k \vee |xw_0|_k \leq a \quad , \quad |w_0|_{k+1} \leq a_1 \ \ell n \ t_0  \\ \\ |s_0|_{k+j}^{\dot{}}\leq b 
_j \left ( \ell n \ t_0 + t_0^{j\beta} \right ) \quad , \quad j = 
0,1,2 \ .\\ \end{array} \right . \eeq

\noi The solution $(w,s,B_b)$ satisfies the estimates
\beq
\label{4.71e}
|w|_{k} \vee |xw|_k \leq Ca \quad , \quad |w|_{k+1} \leq C(a_1 + a^3) 
\ell n \ t_+ \ ,
\eeq
\beq
\label{4.72e}
|s|_{k+j}^{\dot{}} \leq C\left ( b_j + a^2\right ) \left ( \ell n \ 
t_+ + t_+^{j\beta}\right ) \quad , \quad j = 0, 1, 2 \ ,
\eeq
\beq
\label{4.73e}
|B_b|_{k+1}^{\dot{}} \vee |x\cdot B_b|_{k+1}^{\dot{}} \leq C\left ( 
aa_1 + a^2b_0 + a^4\right ) t^{-1} \ell n \ t_+
\eeq

\noi for all $t \geq T$, with $t_+ = t \vee t_0$. \par

In addition $w \in L^{\infty}(I, H^{k+\theta})$ for $0 \leq \theta < 1$.\par

(2) The map $(w_0, s_0) \to (w,s, B_b)$ is continuous for fixed 
$t_0$, on the bounded sets of $X^k$, from the norm of $(w_0, s_0, 0) 
\in
X^{k-1}$ to the norm of $(w, s, B_b)$ in $L^{\infty}(J, X^{k-1})$ for 
any interval $J \subset \subset I$, and in the weak-$*$ sense to 
$L^{\infty}(J,
X^{k})$.}\\

\noi {\bf Proof.} \underbar{Part (1)}. The proof consists in 
exploiting the estimates of Lemmas 4.1 and 4.2 in order to show that 
the
map $\Gamma : (w, s, B_b) \to (w', s', B'_b)$ defined by Proposition 
4.1 with $(w', s')(t_0) = (w,s)(t_0)$ and by (\ref{2.38e}) is a
contraction of a suitable set ${\cal R}$ of ${\cal C}(I, X^k)$ for a 
suitably time rescaled norm of $L^{\infty}(I, X^{k-1})$. More
precisely, let $I = [T,\infty )$ and $t_0 \in I$. For $(w, s, 0) \in 
{\cal C}(I, X^k)$, we define
\beq
\label{4.74e}
\left \{ \begin{array}{l} y = |w|_k \vee |xw|_k \quad , \quad y_1 = 
|w|_{k+1}\\ \\ z_j = |s|_{k+j}^{\dot{}}
\quad , \quad j = 0,1,2 \ ,\\ \end{array} \right . \eeq

\noi and we define ${\cal R}$ by
\bea
\label{4.75e}
&&{\cal R} = \Big \{ (w, s, B_b) \in {\cal C}(I, X^k) :(w,s)(t_0) = 
(w_0,s_0), y\leq Y, y_1 \leq Y_1 \ell n\ t_+, \nn \\
&&z_j \leq Z_j \left ( \ell n \ t_+ + t_+^{j\beta}\right ), j= 0,1,2, 
|B_b|_{k+1}^{\dot{}} \vee |x \cdot B_b|_{k+1}^{\dot{}} \leq N
t^{-1} \ell n \ t_+ \Big \} \eea

\noi for some positive constants $(Y, Y_1, Z_j, N)$ to be chosen 
later. Actually, those constants will turn out to take the form that
appears in (\ref{4.71e})-(\ref{4.73e}). The proof will require 
various lower bounds on $T$ and $t_0$, depending on $(Y,Y_1,Z_j,N)$.
Since those constants will take the form that appear in 
(\ref{4.71e})-(\ref{4.73e}), the lower bounds on $T$ and $t_0$ will 
eventually
be expressed in terms of $(a,a_1,b_j)$, thereby taking the form 
stated in the proposition. \par

We first show that the set ${\cal R}$ is mapped into itself by 
$\Gamma$. Let $(w, s, B_b) \in {\cal R}$. From 
(\ref{4.2e})-(\ref{4.3e}) it
follows that
\beq
\label{4.76e}
|B_a|_{k+1}^{\dot{}} \vee |\nabla (x\cdot B_a)|_{k}^{\dot{}} \leq C\ 
I_0(y^2) \leq C\ Y^2 \ .
\eeq

\noi We now get rid of the variable $B_b$. From  (\ref{4.7e}) 
(\ref{4.8e})  (\ref{4.76e}) it follows that
\bea \label{4.77e}
&&|B'_{b}|_{k+1}^{\dot{}} \vee |x \cdot B'_{b}|_{k+1}^{\dot{}} \leq C 
\ t^{-1}\ I_0 \left (yy_1 + y^2\left (z_0 + I_0(y^2) +
|B_{b}|_k^{\dot{}}\right ) \right )  \nn \\
&&\leq C\  t^{-1} \ell n \ t_+ \left ( YY_1 + Y^2 \left (Z_0 + Y^2 + 
N t^{-1}\right ) \right ) \leq N  t^{-1} \ell n \ t_+
  \eea

\noi with
\beq \label{4.78e}
N = C \left ( YY_1 + Y^2 (Z_0 + Y^2) \right )
\eeq

\noi provided $T$ is sufficiently large, namely $T \geq CY^2$. 
Therefore under that condition, if we choose $N$ in (\ref{4.75e}) as 
given
by (\ref{4.78e}), the condition on $B_b$ in (\ref{4.75e}) is 
automatically reproduced by $\Gamma$, and it remains only to show that
$\Gamma$ reproduces the conditions on $(w,s)$. \par

In order to proceed, we furthermore impose the condition
\beq
|B_{b}|_{k+1}^{\dot{}} \leq Y^2\label{4.79e}
\eeq

\noi which together with (\ref{4.76e}) ensures that
\beq
|B|_{k+1}^{\dot{}} \leq C\ Y^2 \ . \label{4.80e}
\eeq

\noi That condition is ensured by taking $T$ and $t_0$ sufficiently 
large so that
\beq
\label{4.81e}
t\left ( \ell n \ t_+\right )^{-1} \geq N/Y^2 = C \left ( Y_1/Y + Z_0 + Y^2 \right )\eeq\noi for all $t \geq T$, a condition which can be rewri 
tten as
\beq
\label{4.82e}
t_0 \geq T \geq N Y^{-2} \ell n \ t_0  = C \left ( Y_1/Y + Z_0 + Y^2 
\right ) \ell n \ t_0 \ ,
\eeq

\noi where we have included the condition $t_0 \geq T$ for 
completeness. Let now $(w', s')$ be the solution of the system 
(\ref{2.37e})
with initial data $(w',s')(t_0) = (w_0,s_0)$ obtained by using 
Proposition 4.1. We define
\beq
\label{4.83e}
\left \{ \begin{array}{l} y' = |w'|_k \vee |xw'|_k \quad , \quad y'_1 
= |w'|_{k+1}  \\ \\ z'_j = |s'|_{k+j}^{\dot{}}
\quad , \quad j = 0,1,2 \ .\\ \end{array} \right . \eeq

 From Lemma 4.1, more precisely from (\ref{4.4e})-(\ref{4.6e}), we obtain
\beq
\label{4.84e}
|\partial_t y'| \leq E\ y' + t^{-2} \ y'_1 \ ,
\eeq
\beq
\label{4.85e}
|\partial_t y'_1| \leq E\ y'_1 + C\ t^{-2} z_2 \ y' + C\ t^{-1} 
|\nabla (x \cdot B_a)|_{k}^{\dot{}} \ y' \ ,
\eeq
\beq
\label{4.86e}
|\partial_t z'_j| \leq C\ t^{-2} \left (z_1\ z'_j + z_j\ z'_1\right ) 
+ C\ t^{-1} y (y+ \delta_{j_2} \ y_1) + C\ t^{-1+j \beta }|\nabla (x
\cdot B_a)|_{k}^{\dot{}} \ , \eeq

\noi where $\delta_{j_2}$ is the Kronecker symbol and where
\beq
\label{4.87e}
E = C\ t^{-2} \left (z_1 + |B|_{k+1}^{\dot{}} \right ) \left (1 + 
|B|_{k+1}^{\dot{}} \right ) + C\ t^{-1-\beta } |\nabla (x \cdot
B_a)|_{k}^{\dot{}} + C\ t^{-1} |x \cdot B_b|_{k+1}^{\dot{}}\ . \eeq

We want to integrate (\ref{4.84e})-(\ref{4.86e}) between $t_0$ and 
$t$ with appropriate initial data at $t_0$. For that purpose, we
first eliminate the diagonal terms in  (\ref{4.84e}) (\ref{4.85e}) by 
exponentiation according to  (\ref{4.51e}). By
(\ref{4.75e}), (\ref{4.76e}) (\ref{4.80e}) we estimate
\beq
\label{4.88e}
E \leq \bar{E} = C\ t^{-2} \left (Z_1 \ t_+^{\beta} + Y^2 \right ) (1 
+ Y^2) + C\ t^{-1-\beta}\ Y^2 +  C\ N\ t^{-2} \ell n \ t_+
\eeq

\noi and therefore by integration
\bea
\label{4.89e}
\left | \int_{t_0}^t dt'\ E(t') \right |&\leq& C\ (t \wedge t_0)^{-1} 
\left (t_0^{\beta} \ Z_1 + Y^2 \right ) (1 + Y^2)
+  C(t \wedge t_0)^{-\beta } Y^2\nn \\
&&+ C(t \wedge t_0)^{-1 } \ell n\ t_0 \ N \leq C   \eea

\noi where we have used the estimates
\beq
\label{4.90e}
\left | \int_{t_0}^t dt'\ t'^{-2} \ t'^{\beta}_+ \right | \leq (1 - 
\beta )^{-1} (t \wedge t_0)^{-1}\ t_0^{\beta} \ ,
\eeq

\beq
\label{4.91e}
\left | \int_{t_0}^t dt'\ t'^{-2} \ell n \ t'_+\right | \leq (t 
\wedge t_0)^{-1} (1 + \ell n \ t_0)
\eeq

\noi and where the last inequality is achieved for $T$ and $t_0$ 
sufficiently large, namely
\beq
\label{4.92e}
t_0 \geq T \geq t_0^{\beta}\ Z_1(1 + Y^2) + Y^2 + Y^{2/\beta} + N\ 
\ell n\ t_0 \ .
\eeq

We next estimate $y'$ and $y'_1$ by integrating (\ref{4.84e}) 
(\ref{4.85e}) between $t_0$ and $t$ with initial condtions $y'(t_0) 
\leq
a$ and $y'_1(t_0) \leq a_1 \ \ell n \ t_0$ according to (\ref{4.70e}) 
and with $E$ exponentiated to a constant under the
condition (\ref{4.92e}). We obtain
\beq
\label{4.93e}
\left \{ \begin{array}{l} y' \leq C\ a + C \displaystyle{\int} t^{-2} 
\ y'_1   \\ \\ y'_1 \leq C\ a_1 \ \ell n\ t_0 + C\ Z_2
\displaystyle{\int} t^{-2} \ t_+^{2\beta} y' + C \ Y^2 
\displaystyle{\int} t^{-1} \ y' \\ \end{array} \right . \eeq

\noi by the use of (\ref{4.75e}) (\ref{4.76e}) and with the short hand notation
$$\int f(t) = \left | \int_{t_0}^t dt'\ f(t') \right | \ .$$

\noi Let now
\beq \label{4.94e}
Y' = \ \parallel y';L^{\infty}([T, \infty )) \parallel \quad, \quad 
Y'_1 = \ \parallel (\ell n \ t_+)^{-1} y'_1;L^{\infty}([T, \infty ))
\parallel \ .\eeq

\noi (Strictly speaking, we should use a bounded interval $[T, T_1]$ 
with $T_1$ large instead of $[T, \infty )$, check that the
subsequent estimates are uniform in $T_1$, and take the limit $T_1 
\to \infty$ at the end. We omit that step for simplicity). From
(\ref{4.93e}) we obtain
\beq
\label{4.95e}
\left \{ \begin{array}{l} y' \leq C\ a + C \ Y'_1 (t \wedge t_0)^{-1} 
\ \ell n \ t_0\\ \\ y'_1 \leq C\ a_1 \ \ell n\ t_0 + C\ Z_2 \
Y'(t\wedge t_0)^{-1} \ t_0^{2\beta} + C \ Y^2  \ Y'\ \ell n \ t_+  \\ 
\end{array} \right . \eeq

\noi by integration and by the use of (\ref{4.90e}) (\ref{4.91e}), and therefore \beq\label{4.96e}\left \{ \begin{array}{l} Y' \leq C\ a + C \ 
Y'_1 \ T^{-1} \ \ell n \ t_0\\ \\ Y'_1 \leq C\ a_1 + C\ Z_2\ Y' \ 
T^{-1} \
t_0^{2\beta} + C \ Y^2  \ Y' \\ \end{array} \right . \eeq

\noi which implies $Y' \leq Y$ and $Y'_1 \leq Y_1$ provided we take
\beq \label{4.97e}
Y = C\ a \qquad , \qquad Y_1 = C(a_1 + a^3)
\eeq

\noi and provided we take $T$ and $t_0$ sufficiently large so that
\beq
\label{4.98e}
t_0 \geq T \geq C(a_1/a + a^2) \ell n\ t_0 \vee C\ a \ Z_2 \ 
t_0^{2\beta}/(a_1 + a^3)\ .
\eeq

\noi This shows that the conditions on $(y, y_1)$ in (\ref{4.75e}) 
are preserved by $\Gamma$. \par

We finally estimate $z'_j$, $j = 0,1,2$, by integrating (\ref{4.86e}) 
between $t_0$ and $t$ with initial condition $z'_j(t_0) \leq b_j
(\ell n \ t_0 + t_0^{j\beta})$ according to (\ref{4.70e}). Using 
(\ref{4.75e}) and exponentiating the diagonal term with $z_1z'_j$ to
a constant under the condition
\beq \label{4.99e}
t_0 \geq T \geq C\ Z_1\ t_0^{\beta} \ ,
\eeq

\noi we obtain
\bea
\label{4.100e}
z'_j &\leq& C\ b_j \left (\ell n \ t_0 + t_0^{j\beta}\right ) + C\ 
Z_j \int t^{-2} \left ( \ell n \ t_+ + t_+^{j\beta} \right ) z'_1 \nn
\\ &&+ C\ Y^2 \left ( \ell n \ t_+ + t_+^{j\beta} \right ) + C\ Y \ 
Y_1 \ \delta_{j_2} ( \ell n\ t_+)^2\ .\eea

\noi We let now
\beq
\label{4.101e}
Z'_j = \ \parallel \left ( \ell n\ t_+ + t_+^{j\beta } \right )^{-1} 
z'_j ; L^{\infty}([T, \infty )) \parallel
\eeq

\noi and obtain from  (\ref{4.100e})
$$Z'_j \leq C\ b_j  + C\ Z_j \ Z'_1 \ T^{-1} \ t_0^{\beta} + C\ Y^2 + 
C\ Y\ Y_1 ( \ell n \ t_0)^2 \ t_0^{-2\beta}$$

\noi which implies $Z'_j \leq Z_j$, $j = 0, 1,2$ provided we take
\beq
\label{4.102e}
Z_j = C \left ( b_j + a^2 \right )
\eeq

\noi and provided $T$ and $t_0$ are sufficiently large so that 
(\ref{4.99e}) holds and in addition
\beq
\label{4.103e}
t_0^{2\beta} \geq C\left ( a_1/a + a^2 \right ) (\ell n \ t_0)^2 \ .
\eeq

We have proved that $\Gamma$ maps ${\cal R}$ into itself provided 
$(Y, Y_1,Z_j, N)$ are chosen according to (\ref{4.97e})
(\ref{4.102e}) (\ref{4.78e}) and provided $t_0$ and $T$ are taken 
sufficiently large according to (\ref{4.82e})
(\ref{4.92e}) (\ref{4.98e}) (\ref{4.99e}) (\ref{4.103e}). The latter 
conditions clearly take the form stated in the proposition.\par

We next show that the map $\Gamma$ is a contraction in ${\cal R}$ for 
a suitably time rescaled norm of $L^{\infty}(I, X^{k-1})$. Let
$(w_i, s_i, B_{bi}) \in {\cal R}$ and let $(w'_i, w'_i, B'_{bi}) = 
\Gamma (w_i, s_i, B_{bi})$, $i = 1,2$. As in Lemma 4.2, we define
$(w_-, w_-, B_{b_-}) = (1/2) ((w_1, s_1, B_{b1}) - (w_2, s_2, 
B_{b2}))$ and similarly for the primed quantities. Furthermore, we
define
\beq
\label{4.104e}
\left \{ \begin{array}{l} y_- = |w_-|_{k-1} \vee |xw_-|_{k-1} \quad , 
\quad y_{1_-} = |w_-|_k \\ \\ z_{j_-} = |s_-|_{k+j-1}^{\dot{}} \ ,
\ j = 1, 2 \quad , \quad n_- = |B_{b_-}|_k^{\dot{}} \vee |x \cdot 
B_{b_-}|_k^{\dot{}}
   \\ \end{array} \right . \eeq

\noi and similarly for the primed quantities and for $B_a$ and $B$. \par

 From Lemma 4.2 with $k' = k-1$ and from the fact that the $+$ 
quantities are estimated by the definition (\ref{4.75e}) of ${\cal 
R}$ and
by (\ref{4.76e}), we obtain
\beq
|B_{a_-}|_{k}^{\dot{}} \leq C\ Y\ I_0(y_-)\ ,\label{4.105e}
\eeq
\beq
\label{4.106e}
\parallel x \cdot B_{a_-}; \dot{H}^k \parallel \ \leq C\ Y\ I_{k-2}(y_-) \ ,
\eeq
\bea
\label{4.107e}
|\partial_t y'_-| &\leq& \bar{E} y'_- + C\ t^{-2} \ Y \left \{ (1 + 
Y^2) z_{1_-} + \left ( 1 + Z_0 \ \ell n\ t_+ + Y^2 \right ) \left (
YI_0(y_-) + n_-\right ) \right \}\nn \\
&&+ C \ t^{-1-\beta} \ Y^2\ I_{k-2}(y_-) + C \ t^{-1} \ Y\ n_- + 
t^{-2} \ y'_{1_-} \ ,
\eea

$$|\partial_t y'_{1_-}| \leq \bar{E} y'_{1_-} + C\ t^{-2}  \Big \{ 
\left (Y_1 \ell n\ t_+ + Y^3\right ) z_{1_-} + Y z_{2_-} + \left (
Y_1 \ell n\ t_+ + YZ_0 \ell n\ t_+ + Y^3  \right )$$
\beq
\label{4.108e}
\left ( Y\ I_0(y_-) + n_- \right ) \Big \} + C \ t^{-1} \ Y^2\ 
I_{k-2}(y_-) + C \ t^{-1}\  Y\ n_-   \eeq
\noi where $\bar{E}$ is defined by (\ref{4.88e}),
\bea
\label{4.109e}
|\partial_t z'_{j_-}| &\leq& C\ t^{-2}  \left (Z_1 \ t_+^{\beta} 
(z_{j_-}+ z'_{j_-}) + \delta_{j_2} \ Z_2 \ t_+^{2\beta}\  z_{1_-}
\right) \nn \\ &+& C \ t^{-1}  Y\left ( y_- + \delta_{j_2} \ 
y_{1_-}\right ) + C \ t^{-1+j\beta} \ Y\ I_m(y_-) \ {\rm for}\ j = 1,
2\ ,  \eea

\noi where $m = (k-2) \wedge 0$,
\beq
\label{4.110e}
n'_{-} \leq C\ t^{-1} \ Y\ I_0 \left \{ y_{1_-} + \left ( Z_0 \ \ell 
n\ t_+ + Y^2\right ) y_- + Y \left ( z_{1_-} + Y \ I_0 (y_-) +
n_- \right ) \right \} \ . \eeq

\noi For brevity, we continue the argument with a simplified version 
of the system (\ref{4.107e})-(\ref{4.110e}) where we eliminate the
diagonal terms and in particular $\bar{E}$ by exponentiation 
according to (\ref{4.51e}), and where we eliminate the constants and 
the
factors $(Y, Y_1,Z_j, N)$. As a consequence we shall not be able to 
follow the dependence of the additional lower bounds of $t_0$ and
$T$ on those factors. That dependence is of the same type as that 
encountered in the proof of stability of ${\cal R}$ under $\Gamma$.
Thus we rewrite (\ref{4.107e})-(\ref{4.110e})  as
\beq
\label{4.111e}
|\partial_t y'_-| \leq t^{-2} \left \{ z_{1_-} + \ell n\ t_+ \left 
(I_0 (y_-) + n_- \right ) \right \} + t^{-1-\beta} \
I_{k-2}(y_-) + t^{-1}\ n_- + t^{-2} y'_{1_-} \ ,
\eeq
\beq
\label{4.112e}
|\partial_t y'_{1_-}| \leq t^{-2} \left \{ z_{2_-} + \ell n\ t_+ 
\left (z_{1_-} + I_0 (y_-) + n_- \right ) \right \} + t^{-1} \
I_{k-2}(y_-) + t^{-1}\ n_- \ ,
\eeq
\beq
\label{4.113e}
|\partial_t z'_{j_-} | \leq t^{-2} \left (t_+^{\beta}z_{j_-} + 
\delta_{j_2} t_+^{2\beta} z_{1_-}\right ) + t^{-1} \left (y_- + 
\delta_{j_2}
y_{1_-} \right ) + t^{1 + j\beta}\ I_m(y_-) \ ,
\eeq
\beq
\label{4.114e}
n'_{-} \leq t^{-1} \ I_0 \left ( y_{1_-} + \ell n\ t_+ \ y_- + 
z_{1_-} + I_0 (y_-) +
n_- \right )  \ . \eeq

We now define
\beq
\label{4.115e}
\left \{ \begin{array}{l} Y_- = \ \parallel y_-;L^{\infty}([T, \infty 
)) \parallel \ ,  \quad Y_{1_-} = \ \parallel (\ell n \
t_+)^{-1}y_{1_-} ; L^{\infty}([T, \infty ))\parallel \ ,\\ \\ Z_{j_-} 
= \ \parallel t_+^{-j\beta} z_{j_-};L^{\infty}([T, \infty ))
\parallel  \quad , \ j = 1,2 \ , \\ \\ N_- = \parallel t \
t_+^{-\beta}n_- ; L^{\infty}([T, \infty ))\parallel \\
   \\ \end{array} \right . \eeq

\noi and similarly for the primed quantities. Using those definitions 
and omitting the $-$ indices in the remaining part of the
contraction proof, we obtain from (\ref{4.111e})-(\ref{4.114e})
\beq
\label{4.116e}
\partial_t y' \leq t^{-2} \left \{ Z_1 \ t_+^{\beta} + Y\ell n\ t_+ + 
N t^{-1} t_+^{\beta} \ell n \ t_+ \right \} + t^{-1-\beta} Y +
t^{-2} t_+^{\beta} N + t^{-2} \ell n \ t_+ Y'_1 \ ,  \eeq
\beq
\label{4.117e}
\partial_t y'_1 \leq t^{-2} \left \{ Z_2 \ t_+^{2\beta} + Z_1 
t_+^{\beta} \ell n\ t_+ + Y \ell n \ t_+ + N t^{-1} t_+^{\beta} \ell 
n \ t_+
\right \} + t^{-1} Y + t^{-2} \ t_+^{\beta} N \ ,  \eeq
\beq
\label{4.118e}
\partial_t z'_j \leq t^{-2} \ t_+^{(j+1)\beta} (  Z_j + Z_1 ) + 
t^{-1} \left ( Y + \delta_{j_2} \ell n\ t_+  \ Y_1 \right ) +
t^{-1+j\beta} \ Y
  \ ,  \eeq
\beq
\label{4.119e}
n' \leq t^{-1} \left \{ (Y+Y_1) \ell n\ t_+  +  Z_1 \ t_+^{\beta} + N 
t^{-1} \ t_+^{\beta} \right \}  \ . \eeq

\noi Integrating (\ref{4.116e})-(\ref{4.118e}) between $t_0$ and $t$ 
with initial condition $(y', y'_1, z'_j)(t_0) = 0$, using
(\ref{4.90e}) (\ref{4.91e}) and similar estimates, and omitting again 
some absolute constants, we obtain
\bea
\label{4.120e}
y' &\leq& (t \wedge t_0)^{-1} \left \{ Z_1 \ t_0^{\beta} + Y\ \ell n\ 
t_0 + N(t\wedge t_0)^{-1} \ t_0^{\beta} \ \ell n\ t_0 \right \} + (t
\wedge t_0)^{-\beta} \ Y \nn \\
&&+ (t \wedge t_0)^{-1} \ t_0^{\beta} \ N + (t \wedge t_0)^{-1} \ell 
n\ t_0 \ Y'_1 \ ,\eea
\bea
\label{4.121e}
y'_1 &\leq& (t \wedge t_0)^{-1} \left \{ Z_2 \ t_0^{2\beta} + Z_1 \ 
t_0^{\beta} \ \ell n\ t_0 + Y\ \ell n\ t_0 + N(t\wedge t_0)^{-1}
\ t_0^{\beta} \ \ell n\ t_0 \right \} \nn \\
&&+ \ell n \ t_+ \ Y + (t \wedge t_0)^{-1} \ t_0^{\beta} N \ ,\eea
\beq
\label{4.122e}
z'_{j}  \leq (t \wedge t_0)^{-1}  \ t_0^{\beta} \ t_+^{j\beta} (Z_j + 
Z_1) + Y \ell n \ t_+  + \delta_{j_2} Y_1 (\ell n\ t_+)^2  +
t_+^{j\beta} \ Y \ .  \eeq

\noi Substituting (\ref{4.120e})-(\ref{4.122e}) and (\ref{4.119e}) 
into the primed analog of the definition (\ref{4.115e}) (and with the
$-$ indices omitted), we obtain
\beq
\label{4.123e}
\left \{ \begin{array}{l} Y' \leq T^{-1} \left \{ Z_1 \ t_0^{\beta} + 
Y \ell n \ t_0 + N t_0^{\beta} + Y'_1 \ell n \ t_0 \right \} +
T^{-\beta} Y \\ \\ Y'_1 \leq T^{-1} \left \{ Z_2 \ t_0^{2\beta} + Z_1 
\ t_0^{\beta} + Y  + N t_0^{\beta} \right \} +
Y \\ \\ Z'_j \leq T^{-1}\  t_0^{\beta} (Z_j +  Z_1) + Y + 
\delta_{j_2} \ Y_1 \ t_0^{-2\beta} (\ell n\ t_0)^2 \\ \\ N' \leq (Y + 
Y_1)
\ t_0^{-\beta} \ell n\ t_0 + Z_1 + T^{-1} N\ . \\
   \end{array} \right . \eeq

\noi Substituting $Y'_1$ from the second inequality into the first 
one, we recast (\ref{4.123e}) into the form
\beq
\label{4.124e}
\left \{ \begin{array}{l} Y' \leq \varepsilon (Y + Z_1 + Z_2 + N) \\ 
\\ Y'_1 \leq Y + \varepsilon (Z_1 + Z_2 + N) \\ \\ Z'_j \leq Y +
\varepsilon (Y_1 + Z_1 + Z_j) \\ \\ N' \leq Z_1 + \varepsilon (Y + Y_1 + N)
  \\
   \end{array} \right . \eeq

\noi where $\varepsilon$ can be made arbitrarily small by taking 
$t_0$ and $T$ sufficiently large. We then define
\beq
\label{4.125e}
X = Y + Z_1/4 + (Y_1 + Z_2 + N)/8
\eeq

\noi and similarly for the primed quantities. It then follows from 
(\ref{4.124e}) that
$$X' \leq (1/2 + O(\varepsilon )) X$$

\noi and therefore $\Gamma$ is a contraction of ${\cal R}$ in the 
norms defined by (\ref{4.104e}) (\ref{4.115e}) for $T$ and $t_0$
sufficiently large. By a standard compactness argument, ${\cal R}$ is 
easily shown to be closed for the latter norms. Therefore
$\Gamma$ has a unique fixed point in ${\cal R}$. \par

The uniqueness of the solution in ${\cal C} (I, X^k)$ under the 
assumption (\ref{4.37e}) follows from Proposition 4.2, part (1). The
estimates (\ref{4.71e})-(\ref{4.73e}) follow from the definition 
(\ref{4.75e}) of ${\cal R}$ and from the choices (\ref{4.97e})
(\ref{4.102e}) (\ref{4.78e}) of $(Y, Y_1, Z_j, N)$. The dependence of 
$t_0$ and $T$ on $(w_0,s_0)$ stated in the proposition
follows from that choice and from the fact that the lower bounds on 
$t_0$ and $T$ are expressed in terms of those quantities, as
explained above. \par

Finally, the fact that $w \in L^{\infty}(I, H^{k+\theta})$ for $0 
\leq \theta < 1$ follows immediately from (\ref{4.4e}) by
substituting the estimates contained in the definition of ${\cal R}$ 
into the RHS, and integrating in time after exponentiation of
the diagonal term. The crucial point is the fact that the 
contribution of the term in $x\cdot B_a$ is integrable in time for
$\theta < 1$. \\

\noi \underbar{Part (2)}. Let $(w_i, s_i, B_{bi})$, $i = 1,2$, be two 
solutions of the system (\ref{2.34e}) (\ref{2.35e}) with
initial conditions $(w_i, s_i)(t_0) = (w_{i_0}, s_{i_0})$ as obtained 
in Part (1). In particular those solutions satisfy the
estimates (\ref{4.71e})-(\ref{4.73e}). We define $(y_-, y_{1_-}, 
z_{j_-}, n_-)$ by (\ref{4.104e}). By the same estimates as in the
contraction proof, we obtain (\ref{4.111e})-(\ref{4.114e}) with the 
primes omitted. We next define $(Y_-, Y_{1_-}, Z_{j_-},
N_-)$ by (\ref{4.115e}). Omitting again the $-$ indices as in the 
contraction proof, we obtain (\ref{4.116e})-(\ref{4.119e}) with
the primes omitted. Integrating the first three equations thereof 
between $t_0$ and $t$ with initial condition $(y,y_1,z_j)(t_0)
= (y_0, y_{10}, z_{j0})$, we obtain
$$\left \{ \begin{array}{l} y \leq y_0 + \ {\rm RHS\ of}\ 
(\ref{4.120e}) \\ \\ y_1 \leq y_{10} + \ {\rm RHS\ of}\ 
(\ref{4.121e})  \\ \\
z_j \leq z_{j0} + \ {\rm RHS\ of}\  (\ref{4.122e})
  \\
   \end{array} \right .$$

\noi and therefore
$$\left \{ \begin{array}{l} Y \leq y_0 \\ \\ Y_1 \leq y_{10} (\ell n 
\ t_0)^{-1}  \\ \\
Z_j \leq z_{j0} \ t_0^{-j\beta}\\ \\ N \leq \\ \end{array} \right \} 
+ \ {\rm RHS\ of}\  (\ref{4.123e})$$

\noi which by the same argument as in the contraction proof, implies 
that $X$ defined by (\ref{4.125e}) satisfies
\beq \label{126e}
X \leq C \left ( y_0 + y_{10} (\ell n\ t_0)^{-1} + z_{10} \ 
t_0^{-\beta} + z_{20} \ t_0^{-2\beta } \right ) \ .
\eeq

\noi This proves the continuity of the map $(w_0, s_0) \to (w, s, 
B_b)$ from the norm of $(w_0, s_0, 0)$ in $X^{k-1}$ on the bounded
sets of $X^k$ to the norm $(w, s, B_b)$ in ${\cal C}(I, X^{k-1})$ in 
the norms defined by (\ref{4.115e}), and a fortiori to the norm of
$(w, s, B_b)$ in $L^{\infty} (J, X^{k-1})$ for $J \subset \subset I$. 
The last continuity follows by a standard compactness
argument.\par\nobreak \hfill $\sq$\par

We conclude this section by deriving asymptotic properties in time of 
the solutions of the auxiliary system (\ref{2.34e})
(\ref{2.35e}) obtained in the previous proposition. We prove in 
particular the existence of asymptotic states $(w_+, \sigma_+)$ for 
those
solutions.\\

\noi {\bf Proposition 4.4.} {\it Let $k > 3/2$, $0 < \beta < 1/2$, $T 
\geq 1$, $I = [T,\infty )$ and let $(w,s, B_b) \in {\cal C}(I,
X^{k})$ be a solution of the auxiliary system (\ref{2.34e}) 
(\ref{2.35e}) satisfying (\ref{4.37e}). Then\par

(1) There exists $w_+ \in H^k$ such that $x w_+ \in H^k$, $w(t)$ 
tends to $w_+$ strongly in $H^k$ and $x w(t)$ tends to $xw_+$
strongly in $H^{k'}$ for $0 \leq k' < k$ and weakly in $H^k$ when $t 
\to \infty$. Furthermore the following estimates hold
\beq
\label{4.127e}
|w_+|_k \vee |xw_+|_k \leq a_{\infty} = \lim_{t \to \infty} \sup 
(|w(t)|_k \vee |xw(t)|_k )  \ ,
\eeq
\beq
\label{4.128e}
|w(t) - U^*(1/t) w_+|_k \vee |w(t) - w_+|_k \leq C\ t^{-\beta} \ ,
\eeq
\beq
\label{4.129e}
|xw(t) - U^*(1/t) xw_+|_{k-1} \vee |xw(t) - xU^*(1/t) w_+|_{k-1} 
\leq C\ t^{-2\beta} \ ,
\eeq
\beq
\label{4.130e}
|xw(t) - xw_+|_{k-1} \leq C\left ( t^{-2\beta} + t^{-1/2} \right ) \ .
\eeq

\noi The constants $C$ in (\ref{4.128e})-(\ref{4.130e}) depend on 
$\parallel \Lambda ; L^{\infty}(I)\parallel$, where $\Lambda$ is
defined by (\ref{4.37e}).\par

Assume in addition that $w \in L^{\infty} (I, H^{k+\theta})$ for all 
$\theta$, $0 \leq \theta < 1$, and define $(W,S)$ by (\ref{2.40e}).
Then\par

(2) For all $\theta$, $0 \leq \theta < 1$, $w_+ \in H^{k+\theta}$ and 
$w(t)$ tends to $w_+$ strongly in $H^{k+\theta}$ when $t \to
\infty$. Furthermore $xW(t) \in L^{\infty} (I, H^{k+\theta-1})$ and 
$x W(t)$ tends to $xw_+$ strongly in $H^{k+\theta-1}$ when $t \to
\infty$.\par

(3) For all $j$, $0 \leq j < 2$, $S \in {\cal C}([1, \infty ), 
K^{k+j})$ and $S$ satisfies the estimate
\beq
\label{4.131e}
|S|_{k+j}^{\dot{}} \leq C_{\varepsilon} \ t^{\beta (j + \varepsilon)}
\eeq
\noi for any $\varepsilon > 0$. \par

Furthermore there exists $\sigma_+$ such that for all $\theta$, $0 
\leq \theta < 1$, $\sigma_+ \in K^{k+ \theta}$, $s - S$ tends to
$\sigma_+$ strongly in $K^{k+ \theta}$, and the following estimate holds~: }
\beq
\label{4.132e}
|s- S- \sigma_+|_{k+\theta }^{\dot{}} \leq C \ t^{-\beta (1 - \theta)} \ .
\eeq
\vskip 3 truemm

\noi {\bf Proof.} \underbar{Part (1)}. Let $\widetilde{w}(t) = U(1/t) 
w(t)$ and $\widetilde{w}_0 = \widetilde{w}(t_0)$ for some $t_0 \in
I$. From (\ref{2.34e}) we obtain
\bea
\label{4.133e}
\partial_t (\widetilde{w} - \widetilde{w}_0) &=& U(1/t) \Big \{ 
t^{-2} Q(s + B, w) - i(2t^2)^{-1} (2B\cdot s + B^2) w \nn \\
&&+ i \ t^{-1}
\left ( (x \cdot B_a)_S + x \cdot B_b\right ) w \Big \} \ . \eea

\noi By estimates similar to, but simpler than, those of Lemma 4.1, we estimate
\bea
\label{4.134e}
&&|\partial_t (\widetilde{w} - \widetilde{w}_0)|_k \leq C \Big \{ 
t^{-2} \Big ( (|s|_k^{\dot{}} + |B|_k^{\dot{}}) (|w|_{k+1} +
|B|_k^{\dot{}}\ |w|_k) + |s|_{k+1}^{\dot{}} \ |w|_k \Big ) \nn \\
&&+ \left ( t^{-1- \beta } \parallel x \cdot B_a;\dot{H}^{k+1}
\parallel \ + \ t^{-1} |x \cdot B_b|_k^{\dot{}} ) |w|_k \right )\Big \}\nn \\
&&\leq C\ t^{-1-\beta}   \eea

\noi by (\ref{4.2e}) (\ref{4.3e}) (\ref{4.7e}) (\ref{4.8e}) and 
(\ref{4.37e}), so that by integration
\beq
\label{4.135e}
|\widetilde{w}(t) - \widetilde{w}(t_0)|_k \leq C \left ( t \wedge 
t_0\right )^{-\beta} \ .
\eeq

\noi This implies the existence of $w_+ \in H^k$ such that $w(t) \to 
w_+$ strongly in $H^k$, and the first estimate of (\ref{4.128e}).
The second estimate follows from the first one and from the fact that
\beq
\label{4.136e}
|(U(1/t) - \1 )w|_k \ \leq t^{-1/2} |w|_{k+1} \ \leq C\ t^{-1/2} \ell n \ t
\eeq

\noi by (\ref{4.37e}), and that $\beta < 1/2$. \par

Let similarly $\widetilde{xw}(t) = U(1/t) xw(t)$ and 
$(\widetilde{xw})_0 = \widetilde{xw}(t_0)$. From (\ref{2.36e}) we 
obtain
\bea
\label{4.137e}
&&\partial_t (\widetilde{xw} - (\widetilde{xw})_0)= U(1/t) \Big \{ 
t^{-2} Q(s + B, xw) - t^{-2} (s+B)w - it^{-2} \nabla w\nn\\
&&- i(2t^2)^{-1} (2B\cdot s + B^2) xw + i t^{-1} \left ( (x \cdot 
B_a)_S + x \cdot B_b\right ) xw \Big \}\eea

\noi and by similar estimates as previously
\bea
\label{4.138e}
&&|\partial_t (\widetilde{xw} - (\widetilde{xw})_0)|_{k-1} \leq C 
\Big \{ t^{-2} \Big [ |w|_k + \left ( |s|_k^{\dot{}} +
|B|_k^{\dot{}} \right ) \left ( |xw|_k + |w|_{k-1}  \right .\nn \\
&&\left . + |B|_k^{\dot{}}\ |xw|_{k-1} \right ) \Big ] + 
t^{-1-2\beta} \parallel x \cdot B_a; \dot{H}^{k+1}\parallel \ |xw|_k 
+ t^{-1} |x\cdot B_b|_k^{\dot{}} \ |xw|_{k-1} \Big \}\nn \\
&&\leq C \left ( t^{-2} \ell n \ t + t^{-1-2\beta} \right ) \eea

\noi so that by integration
\beq
\label{4.139e}
|\widetilde{xw}(t) - \widetilde{xw}(t_0)|_{k-1} \leq C \left ( t 
\wedge t_0\right )^{-2\beta}
\eeq

\noi which together with (\ref{4.135e}) implies that $xw_+ \in 
H^{k-1}$ and that the first estimate of (\ref{4.129e}) holds. The fact
that $xw_+ \in H^k$ and that $xw_+$ satisfies (\ref{4.127e}) and the 
additional convergences of $xw$ to $xw_+$ follow therefrom by
standard compactness and interpolation arguments. \par

The second estimate of (\ref{4.129e}) follows from the first one and 
from the identity
\beq
\label{4.140e}
U^*(1/t) x = \left ( x + i \ t^{-1} \nabla \right ) U^*(1/t)
\eeq

\noi which implies
\beq
\label{4.141e}
|U^*(1/t)xw_+ - xU^*(1/t)w_+|_{k-1} \leq t^{-1} |w_+|_k \ .
\eeq

\noi Finally (\ref{4.130e}) follows from (\ref{4.129e}) and from
\beq
\label{4.142e}
|(U^*(1/t)- \1) xw_+|_{k-1} \leq t^{-1/2} |xw_+|_k \leq a_{\infty} \ 
t^{-1/2}\ .
\eeq
\vskip 3 truemm

\noi \underbar{Part (2)}. The first statement follows from Part (1) 
by standard compactness and interpolation arguments. The second
statement follows from the first one, and from (\ref{4.140e}) which implies
  \beq
\label{4.143e}
|U^*(1/t)xw_+ - xW(t)|_{k+\theta -1} \leq t^{-1} |w_+|_{k+\theta} \ .
\eeq
\vskip 3 truemm

\noi \underbar{Part (3)}. We estimate in the same way as in Lemma 4.1
\beq
\label{4.144e}
\partial_t |S|_{k+j}^{\dot{}} \leq C\  t^{-1} |W|_{k+j-1} \ |W|_k + C 
\ t^{-1+ \beta (j + \varepsilon)} I_0 \left ( |xW|_{k-
\varepsilon}^2\right ) \eeq

\noi for $0 < \varepsilon < k - 3/2$. The first statement and the 
estimate (\ref{4.131e}) then follow from Part (2) and from
(\ref{4.144e}) by integration. \par

Let now $q = w - W$ and $\sigma = s - S$. From (\ref{2.34e}), we obtain
\beq
\label{4.145e}
\partial_t \sigma = t^{-2} s \cdot \nabla s + t^{-1} \nabla g (q, w + 
W) - t^{-1} \nabla (x \cdot B_a)_L \ (q, w + W)
\eeq

\noi and by estimates similar to those in Lemma 4.1

\begin{eqnarray*}
\partial_t |\sigma|_{k+\theta}^{\dot{}} &\leq& C\  t^{-2} \left \{ 
|s|_k^{\dot{}} \ |s|_{k+2}^{\dot{}} + |s|_{k+1}^{\dot{}2}\right \}  +
C\ t^{-1}\ a |q|_k  \\
&&+ C\ t^{-1+ \beta (1 + \theta)} a\ I_m \left ( |xq|_{k-1} \right ) 
\end{eqnarray*}

\noi where $a$ is defined by (\ref{4.1e}) and $m = (k-2) \wedge 0$,
$$ \cdots \leq C \left ( t^{-2+2\beta} \ell n\ t + t^{-1-\beta} + 
t^{-1-\beta (1 - \theta )} \right )$$

\noi by (\ref{4.128e}) (\ref{4.129e}), and therefore by integration
$$|\sigma (t) - \sigma (t_0)|_{k+ \theta}^{\dot{}} \leq C (t \wedge 
t_0)^{-\beta (1 - \theta )}$$

\noi from which the result follows.\par\nobreak \hfill $\sq$\par

\noi {\bf Remark 4.3}. In Proposition 4.4, we have stated the 
asymptotic properties of $(w, s)$ that follow readily from the bounds 
on
the solutions obtained in Proposition 4.3, expressed in terms of 
$(W,S)$ defined by (\ref{2.40e}). However part of the results hold
under more general assumptions. For instance if we drop the 
assumptions on the higher norms $|w|_{k+1}$ and $|s|_{k+2}^{\dot{}}$, 
we
still get the existence of a limit $w_+$ of $w(t)$ with $w_+ \in H^k$ 
and $xw_+ \in H^k$, with the estimate (\ref{4.129e}) for$xw$ and
a similar estimate for $w$. On the other hand, we could have 
expressed the asymptotic properties of $(w,s)$ in terms of the simpler
$(W,S)$ defined by (\ref{2.41e}). However the convergence properties 
of $w$ would be weaker (compare (\ref{4.130e}) with
(\ref{4.129e})), thereby yielding weaker convergence properties of 
$s- S$ in Part (3).

\mysection{Cauchy problem at infinity for the auxiliary system}
\hspace*{\parindent}
In this section, we construct the wave operators for the auxiliary 
system (\ref{2.34e}) (\ref{2.35e}) in the difference form
(\ref{2.49e}) (\ref{2.50e}) for infinite initial time, for the choice 
of $(W,S)$ given by (\ref{2.40e}). In the same spirit as in
Section 4, we first solve the linearized version (\ref{2.54e}) of the 
system (\ref{2.49e}), which together with (\ref{2.55e}) defines a
map $\Gamma : (q, \sigma , B_b) \to (q', \sigma ', B'_b)$. We then 
show that this map is a contraction on a suitable set in suitable
norms. \par

The basic tool of this section again consists of a priori estimates 
for suitably regular solutions of the linearized system
(\ref{2.54e}) (\ref{2.55e}). We first estimate a solution of that 
system at the level of regularity where we shall eventually solve the
auxiliary system (\ref{2.49e}) (\ref{2.50e}). \\

\noi {\bf Lemma 5.1.} {\it Let $k > 3/2$, $0 < \beta < 1/2$, $T \geq 
1$, $I = [T,\infty )$. Let $(W,S,0) \in {\cal C}(I,
X^{k+1})$ with $(U(1/t) W,S,0) \in {\cal C}^1(I,X^{k})$ and let $(q, 
\sigma, B_b) \in {\cal C}(I,X^{k})$, with $W,xW$, $q,xq \in
L^{\infty}(I,H^k)$. Let $I' \subset I$ be an interval, let $(q', 
\sigma ')$ be a solution of the system (\ref{2.54e}) with $(q', \sigma
', 0) \in {\cal C}(I',X^{k})$ and define $B'_b$ by (\ref{2.55e}). Let 
$0 \leq \theta \leq 1$ and $k \leq \ell \leq k+2$. Then the
following estimates hold~:
\beq \label{5.1e}
|G|_{k+1}^{\dot{}} \leq C I_0 \left ( |q|_k \ |x(2W + q)|_k \right ) \ ,
\eeq
\beq
\label{5.2e}
\parallel x \cdot G; \dot{H}^{k+1}\parallel \ \leq |\nabla (x\cdot 
G)|_k^{\dot{}} \leq C\ I_0 \left ( ( |xq|_k + |q|_k) (x
|2W+q|_k)\right ) \ , \eeq
\bea
\label{5.3e}
&&\left |\partial_t| q'|_{k+ \theta } \right | \leq M_4(\theta, q') + 
C \Big \{ t^{-2} \left ( |\sigma|_{k+ \theta + 1}^{\dot{}} + 
|\sigma|_{k+
\theta}^{\dot{}}\ |B|_{k+ \theta}^{\dot{}} \right . \nn \\
&&+ |G + B_b|_{k+ \theta}^{\dot{}} \left ( 1 + |2S+B + B_*|_{k+ 
\theta}^{\dot{}} \right ) + t^{-1-\beta (1 - \theta )} \parallel x
\cdot G; \dot{H}^{k+1}\parallel \nn \\
&&+ t^{-1} |x\cdot B_b|_{k+ \theta}^{\dot{}} \Big \} |W|_{k+ \theta + 
1} + |R_1(W,S,0)|_{k+ \theta} \nn \\
&&\equiv M_4 (\theta , q') + M_5(\theta ) |W|_{k + \theta + 1} + |R_1 
(W,S,0)|_{k+ \theta}\eea

\noi where $M_4(\theta, \cdot )$ is defined by the RHS of (\ref{4.4e}),}
\bea
\label{5.4e}
&&\left |\partial_t|x q'|_{k} \right | \leq M_4(0, xq') + M_5(0) 
|xW|_{k+1} \nn \\
&&+ C \ t^{-2} \left ( |q'|_{k+ 1} + |s+B|_{k}^{\dot{}}\ |q'|_{k} + 
|\sigma + G + B_b|_{k}^{\dot{}} |W|_k \right ) + |xR_1(W,S,0)|_k
\ ,\nn \\\eea
\bea
\label{5.5e}
&&\partial_t| \sigma '|_{\ell}^{\dot{}} \leq  C \  t^{-2}\Big \{ 
|s|_{k+ 1}^{\dot{}}
|\sigma '|_{\ell}^{\dot{}} + \chi (\ell \geq k+1) |s|_{\ell}^{\dot{}} 
\  |\sigma '|_{k+ 1}^{\dot{}} + |\sigma |_{\ell}^{\dot{}}\  |S|_{k+
1}^{\dot{}} + |\sigma|_{k}^{\dot{}} |S|_{\ell + 1}^{\dot{}} \Big \} \nn \\
&& + C \ t^{-1} \left ( |W|_{k+1} + |q|_k \right ) |q|_{\ell - 1} + 
C\ t^{-1+\beta (\ell - k)} |\nabla (x \cdot G)|_{k}^{\dot{}} +
|R_2(W,S)|_{\ell}^{\dot{}}\ , \eea \bea
\label{5.6e}
&&|B'_b|_{k+1}^{\dot{}} \leq C \ t^{-1} \ I_1 \left ( |q|_{k+1} 
(|W|_{k+1} + |q|_k ) + |s+B|_k^{\dot{}}\  |q|_k (|W|_k + |q|_k) \right
.\nn \\ &&\left . +  |\sigma + B_b + G|_{k}^{\dot{}} \ |W|_{k}^2 
\right ) + |R_3(W,S,0)|_{k+1}^{\dot{}}\ ,
\eea
\bea
\label{5.7e}
&&|xB'_b|_{k+1}^{\dot{}} \leq C \ t^{-1} \ I_0 \left ( (|q|_{k+1} + 
|xq|_k) (|W|_{k+1} + |xw|_k ) + |s+B|_k^{\dot{}}( |xq|_k+ |q|_k)
  \right .\nn \\
&&\left . (|W|_k + |q|_k) + |\sigma + B_b + G|_{k}^{\dot{}} ( 
|xW|_{k} + |W|_k ) |W|_k \right ) + |x \cdot 
R_3(W,S,0)|_{k+1}^{\dot{}} .
\eea
\vskip 3 truemm

\noi {\bf Remark 5.1.} The boundedness assumptions in time of $W$ and 
$q$ ensure that the integrals $I_0$ occurring in (\ref{5.1e})
(\ref{5.2e}) are convergent. Furthermore, by estimates similar to but 
simpler than those of Lemma 4.1, one sees easily that the norms of
the remainders $R_1$ and $R_2$ that occur in (\ref{5.3e}) 
(\ref{5.4e}) (\ref{5.5e}) are finite under the assumptions made on 
$(W,S)$.
On the other hand (see Remark 4.1), the statements on $B'_b$ are non 
empty only in so far as the integrals over $\nu$ in the RHS
of (\ref{5.6e}) (\ref{5.7e}) are convergent. This requires additional 
assumptions on the behaviour of $(W,S)$ and of $(q, \sigma,
B_b)$ at infinity in time, which will be made in due course.\\

\noi {\bf Proof}. The proof is very similar to that of Lemma 4.1. The 
estimates (\ref{5.1e}) (\ref{5.2e}) follow immediately from
(\ref{3.25e}) (\ref{3.28e}) and from (\ref{3.10e}) with $m = \bar{m} = k$.\par

We next estimate $q'$ in $H^{k+\theta}$, starting from (\ref{2.54e}). 
Let $m = k + \theta$. We estimate $\partial_t \parallel
\omega^mq'\parallel_2$ by an energy method in the same way as in 
(\ref{4.9e}). The terms containing $q'$ are estimated in the same way
as in (\ref{4.9e}), thereby yielding $M_4(\theta, q')$, while the 
remaining terms are estimated with the help of (\ref{3.10e}) with
$\bar{m} = m$, supplemented by (\ref{3.19e}) for the term containing 
$x\cdot G$. Together with an elementary estimate of $\partial_t
\parallel q'\parallel_2$ (to which the terms containing $q'$ make no 
contribution), this proves (\ref{5.3e}). \par

We next estimate $xq'$ in $H^k$, starting from (\ref{2.56e}). The 
terms containing $xq'$ or $xW$ explicitly yield the first two terms in
the RHS of (\ref{5.4e}) by the special case $\theta = 0$ of the proof 
of (\ref{5.3e}), while the remaining terms are estimated by
(\ref{3.10e}) with $m = \bar{m} = k$. This proves (\ref{5.4e}). \par

We next estimate $\sigma '$, starting with $\partial_t 
\parallel\omega^{\ell} \sigma '\parallel_2$. The term $s\cdot \nabla 
\sigma '$
is estimated exactly as in the proof of (\ref{4.6e}). The term 
$\sigma \cdot \nabla S$ is estimated directly by Lemma 3.2 as
\bea \label{5.8e}
\parallel\omega^{\ell} (\sigma \cdot \nabla S)\parallel_2\ &\leq& C 
\left ( \parallel\omega^{\ell} \sigma \parallel_2\ \parallel\nabla S
\parallel_{\infty} + \parallel \sigma \parallel_{\infty} \ 
\parallel\omega^{\ell+1} S\parallel_2 \right ) \nn \\
&\leq& C \left ( |\sigma|_{\ell}^{\dot{}} \ |S|_{k+1}^{\dot{}} + 
|\sigma|_{k}^{\dot{}}\ |S|_{\ell+1}^{\dot{}} \right ) \ .\eea

\noi The term containing $g$ is estimated by (\ref{3.10e}) with $(m, 
\bar{m}) = (k \wedge (\ell - 1), k \vee (\ell - 1))$ or $(\ell -
1, k+1)$ as
\beq
\label{5.9e}
\parallel\omega^{\ell-1}(q(q+ 2W))\parallel_2\ \leq C |q|_{\ell - 1} 
\left ( |q|_k + |W|_{k+1} \right ) \ .
\eeq

\noi The term containing $G$ is estimated by (\ref{3.20e}). Together 
with a simpler estimate of $\partial_t \parallel\nabla \sigma
'\parallel_2$, the previous estimates yield (\ref{5.5e}).\par

Finally (\ref{5.6e}) (\ref{5.7e}) follow from (\ref{2.55e}) 
(\ref{3.25e}) (\ref{3.28e}) and from repeated use of (\ref{3.10e}) 
with $m =
\bar{m} = k$. \par\nobreak \hfill $\sq$\par

We shall also need estimates for the difference of two solutions of 
the linearized system (\ref{2.54e}) (\ref{2.55e}) corresponding to
two different choices of $(q, \sigma, B_b)$ but to the same choice of 
$(W, S)$. Those estimates will be provided by Lemma 4.2, since
for such solutions $(q_-, \sigma_-) = (w_-, s_-)$ and $(q'_-, \sigma 
'_-) = (w'_-, s'_-)$ in the notation of that Lemma extended in an
obvious way.\par

We now begin the study of the Cauchy problem for the auxiliary system 
(\ref{2.49e}) (\ref{2.50e}) and for that purpose we first study
that problem for the linearized system (\ref{2.54e}). For finite 
initial time $t_0$, that problem is solved by Proposition 4.1. The
following is a special case of that proposition and of Lemma 5.1\\

\noi {\bf Proposition 5.1.} {\it Let $k > 3/2$, $0 < \beta < 1/2$, $T 
\geq 1$ and $I = [T,\infty )$. Let $(W,S,0) \in {\cal C}(I,
X^{k+1})$ with $(U(1/t) W,S,0) \in {\cal C}^1(I,X^{k})$ and let $(q, 
\sigma, B_b) \in {\cal C}(I,X^{k})$, with $W,xW,q,xq \in
L^{\infty}(I,H^k)$. Let $t_0 \in I$ and $(q'_0, \sigma '_0, 0) \in 
X^k$. Then the system (\ref{2.54e}) has a unique solution
$(q', \sigma ')$ in $I$ such that $(q', \sigma ', 0) \in {\cal 
C}(I,X^{k})$ and $(q', \sigma ') (t_0) = (q'_0, \sigma '_0)$. That
solution satisfies the estimates (\ref{5.3e}) (\ref{5.4e}) 
(\ref{5.5e}) for all $t \in I$, with $G$ estimated by (\ref{5.1e})
(\ref{5.2e}).}\\

In order to study the Cauchy problem with infinite initial time, both 
for the linearized system (\ref{2.54e}) and for the nonlinear
system (\ref{2.49e}) (\ref{2.50e}), we shall need stronger 
assumptions on the asymptotic behaviour in time of $(W,S)$. For 
simplicity,
from now on we make the final choice of $(W,S)$ that will turn out to 
satisfy those assumptions. Thus we choose $(W,S)$ as explained in
Section 2, namely
$$\left \{ \begin{array}{l}
W(t) = U^*(1/t) w_+ \\
  \\
S(t) = \int_1^t dt' \ t'^{-1} \left ( \nabla g (W) - \nabla (x \cdot 
B_a)_L (W)\right ) \\
\end{array} \right . \eqno(2.40)\equiv(5.10)$$

\noi for some fixed $w_+ \in H^{k+\alpha + 1}$ with $x w_+ \in H^{k+ 
\alpha}$ for some $\alpha \geq 1$ (we shall eventually need $\alpha
> 1$), and we define
$$a_+ = |w_+|_{k+ \alpha + 1} \vee |xw_+|_{k+\alpha} \ . \eqno(5.11)$$

\noi Using the fact that
$$\partial_t + i(2t^2)^{-1} \Delta = U^*(1/t) \partial_t \ U(1/t) \eqno(5.12)$$

\noi we recast the remainders that occur in the system (\ref{2.49e}) 
(\ref{2.50e}) into the form
$$R_1(W,S,0) = t^{-2} Q(S + B_*, W) - i(2t^2)^{-1} (2B_* \cdot S + 
B_*^2)W + i t^{-1} (x \cdot B_*)_S W \eqno(5.13)$$
$$R_2(W,S) = t^{-2} S \cdot \nabla S \eqno(5.14)$$
$$R_3(W,S,0) = t^{-1} F_1 \left ( {\rm Im} \bar{W} \nabla W - (S + 
B_*) |W|^2\right ) \eqno(5.15)$$

\noi where $B_* = B_a (W)$ (see Section 2). \par

We shall need the following estimates.\\

\noi {\bf Lemma 5.2.} {\it Let $k > 3/2$ and $0 < \beta < 1/2$. Let 
$w_+ \in H^{k+ \alpha + 1}$ with $xw_+ \in H^{k+\alpha}$ for some
$\alpha \geq 1$. Define $(W,S)$ and $a_+$ by (5.10) (5.11). Then 
$(W,S,0) \in {\cal C}([1, \infty ), X^{k+\alpha})$ and
the following estimates hold}
$$|xW(t)|_{k+\alpha} \leq |xw_+|_{k+ \alpha} + t^{-1} |\nabla 
w_+|_{k+ \alpha} \leq 2a_+  \ , \eqno(5.16)$$
$$|B_*|_{k+\alpha +1}^{\dot{}} \vee |\nabla (x\cdot 
B_*)|_{k+\alpha}^{\dot{}} \leq I_0 \left ( |xW|_{k+\alpha} 
(|xW|_{k+\alpha} +
|W|_{k + \alpha} ) \right ) \leq C\ a_+^2 \ , \eqno(5.17)$$
$$|S|_{k+j}^{\dot{}} \leq C\ a_+^2 \left ( \ell n\ t + t^{\beta (j- 
\alpha)} \right ) \quad {\it for}\ 0 \leq j \leq 2 + \alpha \ ,
\eqno(5.18)$$
$$|R_1(W,S,0)|_{k+\theta }\leq C(a_+) \left ( t^{-2} \ell n\ t + 
t^{-1- \beta (\alpha + 1 - \theta)} \right ) \quad {\it for}\ 0 \leq
\theta \leq 1 \ ,  \eqno(5.19)$$
$$|xR_1(W,S,0)|_{k}\leq C(a_+) \left ( t^{-2} \ell n\ t + t^{-1- 
\beta (\alpha + 1 )} \right )  \ ,  \eqno(5.20)$$
$$|R_2(W,S)|_{k+j}^{\dot{}}\leq C\ a_+^4 \ t^{-2} \ell n\ t \left ( 
\ell n \ t + t^{\beta (j + 1 - \alpha)} \right )
\quad {\it for}\ 0 \leq j \leq 1 + \alpha \ ,  \eqno(5.21)$$
$$|R_3(W,S,0)|_{k+1}^{\dot{}} \vee |x\cdot 
R_3(W,S,0)|_{k+1}^{\dot{}}\leq C\ a_+^2 \ t^{-1} (1 + a_+^2 \ell n\ 
t) \ .
\eqno(5.22)$$ \vskip 3 truemm

\noi {\bf Proof.} We first estimate $x\cdot W$. From the commutation 
relation (\ref{4.140e}), it follows that
$$U(1/t) \ x\ W(t) = x w_+ + i t^{-1} \nabla w_+ \eqno(5.23)$$

\noi which implies (5.16).\par

The estimate (5.17) follows immediately form (\ref{3.25e}) 
(\ref{3.28e}) (\ref{3.10e}) and (5.16). \par

We next estimate $S$. Let $\ell = k + j \geq k$. From (5.10) 
(\ref{3.20e}) we obtain
$$\parallel \omega^{\ell} S \parallel_2 \ \leq C \int_1^t dt' \left 
\{ t'^{-1} \parallel \omega^{\ell - 1} |W|^2 \parallel_2 \ + t'^{-1
+ \beta (\ell - k- \alpha)_+} \ |\nabla (x \cdot B_*)|_{k + 
\alpha}^{\dot{}} \right \} \eqno(5.24)$$

\noi from which (5.18) follows by (\ref{3.10e}), (5.16), (5.17) and 
integration on time, provided $\ell - 1 \leq k + 1 + \alpha$ or
equivalently $j \leq 2 + \alpha$. \par

We next estimate $R_1$. By repeated use of (\ref{3.10e}) and by 
(\ref{3.19e}), we obtain from (5.13)
$$|R_1(W,S,0)|_{k+\theta} \leq C\ t^{-2} \left \{ |S|_{k+\theta 
+1}^{\dot{}} + |B_*|_{k+\theta}^{\dot{}} \left ( 1 +
|S|_{k+\theta}^{\dot{}} + |B_*|_{k+\theta}^{\dot{}} \right ) \right 
\} |W|_{k+ \theta + 1}$$
$$+ C\ t^{-1-\beta (\alpha + 1 - \theta )} \ \parallel x \cdot B_*; 
\dot{H}^{k+1+ \alpha} \parallel \ |w_+|_{k + \theta}
\eqno(5.25)$$

\noi and therefore by (5.17) (5.18) and with $\alpha \geq 1 \geq \theta$
$$|R_1(W,S,0)|_{k+\theta} \leq C\ t^{-2} \ a_+^3 \left ( \ell n \ t + 
t^{\beta (1+ \theta - \alpha)} + a_+^2 (1 + \ell n\ t)  \right ) +
C\ t^{-1-\beta(\alpha + 1 - \theta)}\ a_+^3 \eqno(5.26)$$

\noi which proves (5.19) since $\beta < 1/2$. \par

The proof of (5.20) follows from the fact that
$$x R_1(W,S,0) = L_0 \ x W - t^{-2} (S + B_*) W \eqno(5.27)$$

\noi where the linear operator $L_0$ is defined by rewriting the RHS 
of (5.13) as $L_0 \ W$. The term $L_0\ xW$ is estimated in the same
way as in the proof of (5.19) with $\theta = 0$, and by using in 
addition (5.16), while the last term in (5.27) is estimated by using
(5.17) (5.18) and (\ref{3.10e}). \par

We next estimate $R_2$. By a direct application of Lemma 3.2, we 
obtain from (5.14)
$$\begin{array}{ll}
\parallel \omega^{\ell} R_2(W,S) \parallel_2 \ &\leq C \ t^{-2} \left 
( \parallel S \parallel_{\infty}\ \parallel \omega^{\ell + 1} S
\parallel_2 \ + \ \parallel \nabla S\parallel_3\ \parallel \omega^{\ell} S
\parallel_6 \right )\\
&\\
&\leq C \ t^{-2} |S|_{k}^{\dot{}} \  |S|_{\ell + 1}^{\dot{}}\\ 
\end{array} \eqno(5.28)$$

\noi for $\ell \geq 1$, which yields (5.22) by the use of (5.18). \par

Finally the estimate (5.22) follows readily from (\ref{3.25e}) 
(\ref{3.28e}) (\ref{3.10e}) and from (5.16) (5.17) (5.18). \\

\noi {\bf Remark 5.2.} If one makes the simpler choice (\ref{2.41e}) 
for $(W,S)$, Lemma 5.2 and its proof remain essentially unchanged,
the only difference being that the proof of (5.19) (5.20) now 
requires $\alpha \geq 2$ in order to estimate the contribution of 
$\Delta
w_+$ to $R_1$.\\

We can now solve the Cauchy problem with infinite initial time for 
the linearized system (\ref{2.54e}) (\ref{2.55e}) for the previous
choice of $(W,S)$. \\

\noi {\bf Proposition 5.2.} {\it Let $k > 3/2$, $0 < \beta < 1/2$, $T 
\geq 1$ and $I = [T,\infty )$. Let $w_+ \in H^{k+\alpha + 1}$ with
$xw_+ \in H^{k+\alpha}$ for some $\alpha > 1$ with $\beta (\alpha + 
1) \geq 1$. Define $(W,S)$ and $a_+$ by (5.10) (5.11). Let $(q,
\sigma, B_b) \in {\cal C}(I, X^{k})$ satisfy
$$|q|_k \vee |xq|_k \leq Y\ t^{-1} \ell n\ t \ , \eqno(5.29)$$
$$|q|_{k+1} \leq Y_1\left ( t^{-1} \ell n\ t  + t^{-\alpha\beta} 
\right ) \ , \eqno(5.30)$$
$$|\sigma|_{k+j}^{\dot{}}  \leq Z_j\ t^{-1} \ell n\ t \left ( \ell n 
\ t  + t^{j\beta} \right ) \quad {\rm for}\ j = 0, 1, 2\ ,
\eqno(5.31)$$
$$|B_b|_{k+1}^{\dot{}} \vee |x \cdot B_b|_{k+1}^{\dot{}} \leq N\ 
t^{-1}\ell n\ t \ ,\eqno(5.32)$$
 
\noi for some constants $(Y, Y_1, Z_j, N)$ and for all $t \in I$.\par

Then the linearized system (\ref{2.54e}) (\ref{2.55e}) has a unique 
solution $(q', \sigma ', B'_b) \in {\cal C}(I,X^{k})$ satisfying
$$|q'|_k \vee |xq'|_k \leq Y'\ t^{-1}\ell n \ t \ , \eqno(5.33)$$
$$|q'|_{k+1} \leq Y'_1\left ( t^{-1}\ell n \ t + t^{-\alpha \beta} 
\right )  \ , \eqno(5.34)$$
$$|\sigma|_{k+j}^{\dot{}}  \leq Z'_j\ t^{-1} \ell n\ t \left ( \ell n 
\ t  + t^{j\beta} \right ) \quad {\it for}\ j = 0, 1, 2\ ,
\eqno(5.35)$$
$$|B'_b|_{k+1}^{\dot{}}  \vee |x \cdot B'_b|_{k+1}^{\dot{}} \leq N'\ 
t^{-1} \ell n\ t \ , \eqno(5.36)$$

\noi for some constants $(Y',Y'_1, Z'_j, N')$ depending on 
$(Y,Y_1,Z_j,N,a_+,T)$. The solution is actually unique in ${\cal C}(I,
X^{k})$ under the condition that $(q', \sigma ')$ tends to zero in 
$L^2 \oplus \dot{H}^1$ when $t \to \infty$.} \\

\noi {\bf Remark 5.3.} Whereas the conditions $0 < \beta < 1/2$ and 
$\alpha > 1$ are used in an essential way in the proof of
Proposition 5.2, the condition $\beta (\alpha + 1) \geq 1$ has been 
imposed for convenience only, in order to obtain rather simple and
optimal decay properties for $(q', \sigma ')$, and could be relaxed 
at the expense of a weakening of those properties. For given
$\alpha > 1$, it can be achieved by taking $\beta$ sufficiently close 
to 1/2. Its meaning is that for a given regularity of $w_+$, one
should put a sufficiently large part of the interaction $B_1$ into 
the long range part $(x \cdot B_a)_L$ so as to obtain a sufficiently
good decay of the short range part $(x \cdot B_a)_S$.\\

\noi {\bf Proof.} With Proposition 5.1 available, it is sufficient to 
prove Proposition 5.2 for $T$ sufficiently large, depending
possibly on $(Y,Y_1,Z_j,N,a_+)$. Furthermore it is sufficient to 
solve the system (\ref{2.54e}) for $(q', \sigma ')$, since $B'_b$ is
given by an explicit formula, namely (\ref{2.55e}). The proof 
consists in showing that the solution $(q'_{t_0}, \sigma '_{t_0})$ of 
the
linearized system (\ref{2.54e}) with initial data $(q'_{t_0}, \sigma 
'_{t_0})(t_0) = 0$ for some finite $t_0 \geq T$, obtained from
Proposition 5.1, satisfies the estimates (5.33)-(5.35) uniformly in 
$t_0$ for $T \leq t \leq t_0$, namely with $(Y',Y'_1,Z'_j)$
independent of $t_0$, and that when $t_0 \to \infty$, that solution 
converges on the compact intervals of $I$ uniformly in suitable
norms. \par

Let therefore $T$ be sufficiently large, in a sense to be made clear 
below. We define
$$\left \{ \begin{array}{l} y = |q|_k \vee |xq|_k \qquad, \quad y_1 = 
|q|_{k+1} \ , \\ \\ z_j = |\sigma|_{k+j}^{\dot{}} \qquad , \quad j =
0, 1, 2 \ , \\ \end{array}\right .  \eqno(5.37)$$

\noi and we first take $T$ large enough so that (5.29)-(5.32) imply
$$y \vee y_1 \leq a_+ \eqno(5.38)$$
$$z_j \le a_+^2\ \ell n \ t \quad {\rm for}\ j = 0, 1, 2, \qquad 
|B_b|_{k+1}^{\dot{}} \leq a_+^2 \eqno(5.39)$$

\noi for $t \geq T$. It follows from (\ref{5.1e}) (\ref{5.2e}) (5.29) 
(5.38) that
$$|G|_{k+1}^{\dot{}} \vee |\nabla (x \cdot G)|_{k}^{\dot{}} \leq C\ 
a_+\ I_0(y) \leq C\ a_+\left ( Y t^{-1}\ell n \ t \vee a_+ \right
) \ . \eqno(5.40)$$

\noi Let $t_0 > T$ and let $(q'_{t_0}, \sigma '_{t_0})$ be defined as 
above. We want to estimate $(q'_{t_0}, \sigma '_{t_0})$ for $T
\leq t \leq t_0$. We define
$$\left \{ \begin{array}{l} y' = |q'_{t_0}|_k \vee |xq'_{t_0}|_k 
\qquad, \quad y'_1 = |q'_{t_0}|_{k+1} \ , \\ \\  z'_j = |\sigma
'_{t_0}|_{k+j}^{\dot{}} \qquad , \quad j = 0, 1, 2 \ , \\ \end{array} 
\right . \eqno(5.41)$$
$$\left \{ \begin{array}{l} Y' = \ \parallel t(\ell n\ t)^{-1} 
y';L^{\infty}([T, t_0])\parallel \  , \\ \\ Y'_1 = \ \parallel \left (
t^{-1}\ell n\ t + t^{-\alpha \beta} \right )^{-1} y'_1;L^{\infty}([T, 
t_0])\parallel \ , \\ \end{array} \right . \eqno(5.42)$$
$$ Z'_j = \ \parallel \left ( t^{-1}\ell n\ t (\ell n \ t + t^{j 
\beta}) \right )^{-1} z'_j;L^{\infty}([T, t_0])\parallel \ .
\eqno(5.43)$$

We first estimate $q'_{t_0}$ and more precisely $y'$ and $y'_1$, 
starting from (\ref{5.3e}) (\ref{5.4e}). We estimate $(W,S,B_*)$ by
(5.16)-(5.18), $(\sigma , B_b)$ by (5.39), $x\cdot B_b$ by (5.32), 
$G$ by (5.40) and $R_1$ by (5.19) (5.20), with $\beta (\alpha + 1)
\geq 1$. We obtain
$$|\partial_t y'| \leq C \left \{ \left ( a_+^2 (1 + a_+^2) + N 
\right ) t^{-2} \ell n\ t + a_+^2 \ t^{-1- \beta} \right \}
y'$$
$$+ C \left ( a_+^3 (1 + a_+^2)+ N a_+\right ) t^{-2} \ell n\ t + C\ 
a_+^2 \ Y\ t^{-2-\beta}\ell n\ t + t^{-2} \ y'_1
  + C (a_+) t^{-2} \ell n \ t \ ,\eqno(5.44)$$
$$|\partial_t y'_1| \leq C \left \{ \left ( a_+^2 (1 + a_+^2) + N 
\right ) t^{-2} \ell n\ t \right \}
y'_1 + C \left ( t^{-2+\beta (2 - \alpha )} + t^{-1}\right ) a_+^2\ y'$$
$$+ C \left ( a_+^3 (1 + a_+^2)+ N a_+\right ) t^{-2} \ell n\ t + C\ 
a_+^2 \ Y\ t^{-2}\ell n\ t + C(a_+) \left (  t^{-2} \ell n\ t
  + t^{-1-\alpha\beta} \right ) \ .\eqno(5.45)$$

\noi We integrate (5.44) (5.45) between $t$ and $t_0$, with $y'(t_0) 
= y'_1(t_0) = 0$, exponentiating the diagonal terms to a constant
for $T$ sufficiently large depending on $(a_+, N)$ and substituting 
(5.42) into the remaining terms. We obtain
$$\left \{ \begin{array}{l} y' \leq C\ a_+^2 Y t^{-1-\beta} \ell n \ 
t + CY'_1 \left ( t^{-2} \ell n\ t + t^{-1-\alpha \beta}\right ) +
\left ( CNa_+ + C(a_+)\right ) t^{-1} \ell n \ t  \\ \\  y'_1  \leq 
C\ a_+^2 (Y' + Y)  t^{-1} \ell n \ t + C Na_+ + C(a_+)  \left (
t^{-1} \ell n\ t + t^{-\alpha \beta}\right ) \\ \end{array} \right . 
\eqno(5.46)$$

\noi and therefore by substituting (5.46) into (5.42)
$$\left \{ \begin{array}{l} Y' \leq C\ a_+^2 Y T^{-\beta} + CY'_1 \ 
T^{-\theta} + CNa_+ + C(a_+) \equiv C Y'_1 \ T^{-\theta} + A\\
\\ Y'_1  \leq C\ a_+^2 (Y' + Y) + C N a_+ + C(a_+) \equiv C\ a_+^2\ 
Y' + A_1 \\ \end{array}
\right . \eqno(5.47)$$

\noi with $\theta = 1 \wedge \alpha \beta > \beta$. It follows from (5.47) that
$$\left \{ \begin{array}{l} Y' \leq C\ a_+^2 T^{-\theta} Y' + A + C\ 
T^{-\theta} A_1 \leq  C a_+^2 \ Y\ T^{-\beta} + C N a_+ + C(a_+) \\
\\ Y'_1 \leq C\ a_+^2 \ T^{-\theta} Y'_1 + C \ a_+^2 A + A_1 \leq C\ 
a_+^2 Y + C N a_+(1 + a_+^2) + C(a_+) \\ \end{array} \right .
\eqno(5.48)$$

\noi for $T^{\theta} \geq T^{\beta} \geq C a_+^2$. Finally
$$\left \{ \begin{array}{l} Y' \leq C\ a_+^2 \ T^{-\beta} \ Y + (1+N) 
C(a_+) \\
\\ Y'_1  \leq C\ a_+^2 \  Y + (1+ N) C(a_+) \\ \end{array}
\right . \eqno(5.49)$$

\noi for $T$ sufficiently large depending on $(Y, Y_1, Z_j, N, a_+)$. 
This proves that $q'_{t_0}$ satisfies the estimates (5.33) (5.34)
for $T \leq t \leq t_0$. \par

We next estimate $\sigma '_{t_0}$, and more precisely $z'_j$, 
starting from (\ref{5.5e}), which we rewrite with the help of the
definitions (5.37) (5.41) and of the estimates (5.18) (5.21) and 
(5.38)-(5.40) as
$$|\partial_t z'_j| \leq C \ a_+^2 \ t^{-2} \left \{ \ell n\ t (z_j + 
z'_j) + \delta_{j_2}\  t^{\beta (2 - \alpha )} (z_1 + z'_1) \right
\}$$
$$+ C \ t^{-1} \ a_+ \left ( y + \delta_{j_2}\  y_1 \right ) + C\ 
t^{-1+j\beta } \ a_+\ I_0(y)$$
$$+ C\ a_+^4 \ t^{-2} \ell n\ t \left (\ell n\ t + t^{\beta (j+1 - 
\alpha)} \right ) \eqno(5.50)$$

\noi for $j = 0, 1, 2$ with $\delta_{j_2}$ the Kronecker symbol in 
order to include additional terms for $j = 2$. Using (5.29)-(5.31)
and (5.43), we obtain from (5.50)
$$|\partial_t z'_j| \leq C \ a_+^2 \ t^{-2} \ell n\ t \ z'_j  + C \ 
a_+^2 \ t^{-3} \ell n\ t \Big \{ \ell n\ t(\ell n\ t +
t^{j\beta})Z_j$$
$$+ \delta_{j_2}\  t^{\beta (3- \alpha ) }(Z_1 + Z'_1) \Big \} + C \ 
a_+ \ t^{-2} \Big \{ Y \ell n\ t + \delta_{j_2}\ Y_1 (\ell n\ t
+ t^{1 - \alpha \beta } ) \Big \}$$
$$+ C \ a_+ \ t^{-2+j\beta} \ell n\ t \ Y + C \ a_+^4 \ t^{-2} \ell n 
\ t (\ell n\ t +
t^{\beta (j+1 - \alpha)})\ .\eqno(5.51)$$

\noi We integrate (5.51) between $t$ and $t_0$ with $z'_j(t_0) = 0$, 
exponentiating the diagonal terms to a constant for $T$ sufficiently
large in the sense that
$$T(\ell n \ T)^{-1} \geq C\ a_+^2 \eqno(5.52)$$

\noi and we substitute the result into the definition (5.43), thereby obtaining
$$\begin{array}{ll}
Z'_j &\leq C\ a_+^2 \left ( T^{-1} \ell n \ T \ Z_j + \delta_{j_2} \ 
T^{-1+ \beta (1 - \alpha )} (Z_1 + Z'_1) \right )\\ &\\
&+ C\ a_+ \left ( Y + \delta_{j_2} (T^{-2\beta} + T^{1-(\alpha + 
2)\beta} )Y_1\right ) + C\ a_+^4\\ \end{array}
$$

\noi so that for $\alpha \geq 1$ and $\beta (\alpha + 1) \geq 1$ and 
under the condition (5.52)
$$\left \{ \begin{array}{l} Z'_j \leq C\ a_+^2 \ T^{-1} \ell n \ T\ 
Z_j + C\ a_+\ Y + C\ a_+^4  \qquad, \quad {\rm for}\ j = 0, 1 \ , \\
\\
  Z'_2 \leq C\ a_+^2 \ T^{-1} \ell n \ T(Z_2 + Z_1) + C\ a_+( Y 
+t^{-\beta}Y_1) + C\ a_+^4 \\ \end{array}
\right . \eqno(5.53)$$

\noi for $T$ sufficiently large. This proves that $\sigma '_{t_0}$ 
satisfies the estimate (5.35) for $T \leq t\leq t_0$. \par

We now prove that $(q'_{t_0}, \sigma '_{t_0})$ tends to a limit when 
$t_0 \to \infty$. For that purpose we consider two solutions
$(q'_i, \sigma '_i ) = (q'_{t_i}, \sigma '_{t_i})$, $i=1,2$, of the 
system (\ref{2.54e}) corresponding to the same choice of $(q,
\sigma, B_b)$ and to $t_0 = t_i$, $i = 1,2$, for $T \leq t_1 \leq 
t_2$. Let $(q'_-, \sigma '_- ) = 1/2 (q'_1 - q'_2, \sigma '_1 - \sigma
'_2)$. For fixed $(q, \sigma, B_b)$, the inhomogeneous term in $q'$ 
in the equation for $q'$ is the same and therefore $q'_-$ satisfies
an homogeneous linear equation which preserves the $L^2$ norm, so that
$$\parallel q'_-(t) \parallel_2\ = \ \parallel q'_-(t_1) \parallel_2\ 
= 1/2\parallel q'_2(t_1) \parallel_2\ \leq Y'\ t_1^{-1} \ell n\
t_1 \eqno(5.54)$$

\noi for $T \leq t \leq t_1$, by (5.33) applied to $q'_2$ at $t = t_1 
\in [T, t_2]$. Similarly, $\sigma '_-$ satisfies the equation
$$\partial_t \ \sigma '_- = t^{-2} (S + \sigma ) \cdot \nabla \sigma 
'_- \eqno(5.55)$$

\noi so that by an elementary subestimate of Lemma 4.2 and by (5.18) (5.39)
$$\partial_t \parallel \nabla \sigma '_-\parallel_2\ \leq  C\ a_+^2\ 
t^{-2} \ell n\ t  \parallel \nabla \sigma '_-\parallel_2
\eqno(5.56)$$

\noi and therefore under the condition (5.52)
$$\parallel \nabla \sigma '_-(t) \parallel_2\ \leq \ C\parallel 
\nabla \sigma '_-(t_1) \parallel_2\ \leq C\ Z'_0\  t_1^{-1} \ell n\ 
t_1
\eqno(5.57)$$

\noi for $T \leq t \leq t_1$ by (5.35) applied to $\sigma '_2$ at $t 
= t_1$. \par

 From (5.54) (5.57), it follows that $(q'_{t_0}, \sigma '_{t_0})$ 
converges to a limit $(q', \sigma ' ) \in {\cal C}(I,L^2 \oplus
\dot{H}^1)$ uniformly on the compact subintervals of $I$. From the 
uniform estimates (5.33)-(5.35) and from Lemma 5.1, it then follows
by a standard compactness argument that $(q', \sigma ',0 ) \in {\cal 
C}(I,X^k)$ and that $(q', \sigma ' )$ also satisfies the estimates
(5.33)-(5.35). Clearly $(q', \sigma ' )$ satisfies the system 
(\ref{2.54e}). This completes the existence part of the proof. \par

The uniqueness statement follows immediately from the $L^2$ norm 
conservation for the difference of the $q'$ components and from (5.55)
(5.57) for the difference of the $\sigma '$ components of two solutions.\par

As mentioned at the beginning of the proof, the existence and 
uniqueness of $B'_b$ follow from the fact that it is given by an 
explicit
formula (\ref{2.55e}) and the estimate (5.36) follows immediately 
from (\ref{5.6e}) (\ref{5.7e}) (5.22) (5.38) (5.39) (5.40) with
$$N' = C\ a_+^2(1 + a_+^2) \eqno(5.58)$$

\noi for $T$ large enough to ensure (5.38) (5.39).\par\nobreak \hfill $\sq$\par

We now turn to the main result of this section, namely the fact that 
for $T$ sufficiently large, depending on $a_+$, the auxiliary
system (\ref{2.49e}) (\ref{2.50e}) with $(W,S)$ defined by (5.10) has 
a unique solution $(q, \sigma, B_b)$ defined for all $t \geq T$
and decaying at infinity in the sense of (5.29)-(5.32). In the same 
spirit as for Proposition 4.3, this will be done by showing that the
map $\Gamma : (q, \sigma , B_b) \to (q', \sigma ', B'_b)$ defined by 
Proposition 5.2 is a contraction in suitable norms.\\

\noi {\bf Proposition 5.3.} {\it Let $k > 3/2$ and $0 < \beta < 1/2$. 
Let $w_+ \in H^{k+\alpha +
1}$ with $xw_+ \in H^{k+\alpha}$ for some $\alpha > 1$ with $\beta 
(\alpha + 1) \geq 1$. Define $(W,S)$ and $a_+$ by (5.10) (5.11). Then
\par

(1) There exists $T = T(k, \beta , \alpha, a_+)$, $1 \leq T < \infty$ 
such that the system (\ref{2.49e}) (\ref{2.50e})  has a unique
solution $(q, \sigma, B_b) \in {\cal C}(I, X^{k})$, with $I = [T, 
\infty)$, satisfying the estimates (5.29)-(5.32) for some constants
$(Y, Y_1, Z_j, N)$ depending on $(k, \beta , \alpha , a_+)$. Furthermore
$$|G|_{k+1}^{\dot{}} \vee |\nabla (x \cdot G)|_{k}^{\dot{}} \leq C\ 
a_+\ Y \ t^{-1}\ell n \ t  \eqno(5.40)\equiv(5.59)$$

\noi where $G$ is defined by (\ref{2.46e}). The solution is actually 
unique in ${\cal C}(I, X^{k})$ under the conditions that
$$|xq|_{k} \vee |q|_{k+1 } \vee |\sigma|_{k+2}^{\dot{}} \vee 
|B_b|_{k}^{\dot{}} \in L^{\infty}(I) \ , \eqno(5.60)$$
$$|\sigma|_{k+1}^{\dot{}} \vee t^{2\beta + \varepsilon} \left ( |q|_k 
\vee |xq|_{k-1} \right ) \to 0 \quad {\it when} \ t\to \infty
\eqno(5.61)$$

\noi for some $\varepsilon > 0$. \par
(2) The map $w_+ \to (W + q, S + \sigma , B_b) \equiv (w, s, B_b)$ is 
continuous on the bounded sets of the norm (5.11) from the norm
$|w_+|_k \vee |xw_+|_{k-1}$ for $w_+$ to the norm of $(w, s, B_b)$ in 
$L^{\infty}(J, X^{k-1})$ and in the weak-$*$ sense to
$L^{\infty}(J, X^k)$ for any interval $J \subset \subset I$.}\\

\noi {\bf Proof.} \underbar{Part (1)}. The proof consists in showing 
that the map $\Gamma : (q, \sigma , B_b) \to (q',\sigma ', B'_b)$ 
defined by
Proposition 5.2 is a contraction of a suitable set ${\cal R}$ of 
${\cal C}(I, X^{k})$, with $I = [T, \infty )$, for $T$ sufficiently
large and for a suitably time rescaled norm of $L^{\infty}(I, 
X^{k-1})$. We define $(y, y_1, z_j)$ by (5.37) and we define ${\cal 
R}$ by
$${\cal R} = \Big \{ (q, \sigma ,B_b) \in {\cal C}(I, X^{k}) : y \leq 
Y\ t^{-1}\ell n \ t, y_1 \leq Y_1 \left ( t^{-1} \ell n\ t +
t^{-\alpha \beta} \right ) \ ,$$
$$z_j \leq Z_j\ t^{-1} \ell n \ t \left ( \ell n \ t + t^{j\beta} 
\right )\ , \ j = 0, 1, 2, \ |B_b|_{k+1}^{\dot{}} \vee
|x \cdot B_b|_{k+1}^{\dot{}} \leq N\ t^{-1}\ell n \ t \Big \} \eqno(5.62)$$

\noi for some constants $(Y, Y_1, Z_j, N)$ depending on $a_+$, to be 
chosen later, and we take $T$ large enough so that (5.38) (5.39)
and therefore (5.40) hold for $(q, \sigma , B_b) \in {\cal R}$. \par

We first show that the set ${\cal R}$ is mapped into itself by 
$\Gamma$ for suitable $(Y, Y_1, Z_j, N)$ and sufficiently large $T$.
Let $(q', \sigma ', B'_b) = \Gamma  (q, \sigma ,B_b)$ as defined by 
Proposition 5.2. As mentioned in the proof of that proposition, it
follows from (5.6) (5.7) (5.22) (5.38)-(5.40) that $B'_b$ satisfies 
(5.36) with $N'$ defined by (5.58) so that if we define $N$ by
$$N = C\ a_+^2(1 + a_+^2) \eqno(5.63)$$

\noi then the condition on $B_b$ contained in (5.62) is reproduced by 
$\Gamma$. It remains to estimate $(q' , \sigma ')$. Now from the
proof of Proposition 5.2, it follows that $(q' , \sigma ')$ satisfies 
the estimates (5.33) (5.34) (5.35) with $(Y', Y'_1, Z'_j)$
satisfying (5.49) (5.53), for $T$ sufficiently large depending on 
$(Y, Y_1, Z_j, a_+)$. With $N$ given by (5.63), it follows immediately
from (5.49) that one can choose $Y$ and $Y_1$ depending on $a_+$ such 
that (5.49) implies $Y' \leq Y$ and $Y'_1 \leq Y_1$ for $T$
sufficiently large, namely $T^{\beta} \geq C a_+^2$. Therefore the 
conditions on $q$ in (5.62) are also reproduced by $\Gamma$ for that
choice. Finally, it follows from (5.53) that one can choose $Z_j$, 
actually in the form
$$Z_j = C\ a_+\ Y + C\ a_+^4 \qquad , \quad j = 0, 1, 2 \eqno(5.64)$$

\noi so that (5.53) implies $Z'_j \leq Z_j$ for $T$ sufficiently 
large, in fact for
$$T(\ell n\ T)^{-1} \geq C\ a_+^2 \qquad , \qquad T^{\beta} \geq 
Y_1/Y \ . \eqno(5.65)$$

This completes the proof of the fact that $\Gamma$ maps ${\cal R}$ 
into itself for $(Y, Y_1, Z_j, N)$ chosen as above, depending on
$a_+$, and for $T$ sufficiently large, depending on $a_+$. \par

We next show that the map $\Gamma$ is a contraction in ${\cal R}$ for 
a suitably weighted norm of $L^{\infty}(I, X^{k-1})$. Let $(q_i,
\sigma_i ,B_{bi}) \in {\cal R}$ and let $(q'_i,\sigma '_i ,B'_{bi}) = 
\Gamma (q_i, \sigma_i ,B_{bi})$, $i = 1,2$. We define
$(q_{\pm}, \sigma_{\pm} ,B_{b_{\pm}}) = 1/2 ((q_1, \sigma_1 ,B_{b1}) 
\pm (q_2, \sigma_2 ,B_{b2}))$ and similarly for the primed
quantities. Furthermore we define
$$y_- = |q_-|_{k-1} \vee |xq_-|_{k-1} \qquad , \qquad y_{1_-} = |q_-|_k$$
$$z_{j_-} = |\sigma_-|_{k+j-1}^{\dot{}} \ , \ j = 1, 2\ ,\ n_- = 
|B_{b_-}|_{k}^{\dot{}} \vee |x \cdot B_{b_-}|_{k}^{\dot{}}
\eqno(5.66)$$

\noi and similarly for the primed quantities and for $B_a$ and $B$. \par

We shall estimate $(q_-, \sigma_- ,B_{b_-})$ by Lemma 4.2 applied to 
the solutions $(w'_i , s'_i, B'_{bi})$ $= (W + q'_i , S + \sigma '_i,
B'_{bi})$ of the system (\ref{2.37e}) (\ref{2.38e}) associated with 
$(w_i, s_i, B_{bi}) = (W + q_i , S + \sigma_i,
B_{bi})$. As a consequence we shall apply Lemma 4.2 with $(w_-, s_-) 
= (q_-, \sigma_-), (w'_-, s'_-) = (q'_-, \sigma '_-)$ and $(w_+,
s_+) = (W + q_+, S + \sigma_+), (w'_+, s'_+) = (W + q'_+, S + \sigma 
'_+)$, and for that purpose we shall use the estimates (5.16)
(5.18) for $(W,S)$ and the fact that $(q_i, \sigma_i)$ and therefore 
$(q_+, \sigma_+)$ are not larger than $(W,S)$ in the sense of
(5.38) (5.39). From the first part of the proof, it follows that 
$(q'_i, \sigma '_i)$ and therefore $(q'_+, \sigma '_+)$ and $(w'_+,
s'_+)$ satisfy the same properties. Together with (5.17) (5.18) 
(5.30) (5.63) (5.40), this implies that the available estimates on
$(w_+, s_+, B_{b_+})$ and on $B_{a_+}$, $B_+$ are
$$\left \{ \begin{array}{l} |xw_+|_k \vee |w_+|_{k+1} \leq C \ a_+ \\ \\
|s_+|_{k+j}^{\dot{}} \leq C\ a_+^2 \left ( \ell n\ t + t^{(j-\alpha 
)\beta} \right ) \\ \\
|B_{b_+}|_{k+j}^{\dot{}} \vee |x \cdot B_{b_+}|_{k+j}^{\dot{}} \leq 
C\ a_+^2 (1 + a_+^2) t^{-1}\ell n\ t \\ \\
|B_{+}|_{k+j}^{\dot{}} \vee |\nabla (x \cdot B_{a_+})|_{k}^{\dot{}} 
\leq C\ a_+^2 \\
\end{array} \right . \eqno(5.67)$$

\noi and that the same estimates hold for the primed quantities. \par

 From Lemma 4.2 with $k' = k -1$ and from the estimates (5.67) of the 
$+$ quantities, we obtain (compare with
(\ref{4.105e})-(\ref{4.110e}))
$$|B_{a_-}|_{k}^{\dot{}}  \leq C\ a_+ \ I_0(y_-) \ , \eqno(5.68)$$
$$\parallel x\cdot B_{a_-}; \dot{H}^k\parallel \ \leq C\ a_+ \ 
I_{k-2}(y_-) \ , \eqno(5.69)$$
$$|\partial_t y'_-| \leq E_1\ y'_- + C \ t^{-2} \ a_+ \left \{ (1 + 
a_+^2) z_{1_-} + (1 + a_+^2  \ell n\ t) (a_+ \ I_0(y_-) + n_- )
\right \}$$
$$+ C \ t^{-1-\beta} \ a_+^2\ I_{k-2}(y_-) + C \ t^{-1}\ a_+ \ n_- + 
t^{-2} \ y'_{1_-} \ ,\eqno(5.70)$$
$$|\partial_t y'_{1_-}| \leq E_1\ y'_{1_-} + C \ t^{-2} \ a_+ \left 
\{ a_+^2  \ z_{1_-} + z_{2_-} + (1 + a_+^2\ \ell n\ t) (a_+\
I_0 (y_-) + n_-)  \right \}$$
$$+ C \ t^{-1} \ a_+^2\ I_{k-2}(y_-) + C \ t^{-1}\ a_+ \ n_- \ ,\eqno(5.71)$$

\noi where
$$E_1 = C\ a_+^2 \left \{ (1 + a_+^2) t^{-2} \ell n \ t + 
t^{-1-\beta} \right \}$$
  $$|\partial_t z'_{j_-}| \leq C \ t^{-2} \ a_+^2 \left ( \ell n\ 
t(z_{j_-} + z'_{j_-}) + \delta_{j_2} \ t^{(2- \alpha )\beta}\ z_{1_-}
\right )$$
$$+ C \ t^{-1} \ a_+\ (y_- +  \delta_{j_2} \ y_{1_-}) + C\ t^{-1 + j 
\beta} \ a_+\ I_m (y_-)  \quad {\rm for} j = 1, 2\ , \eqno(5.72)$$

\noi where $m = (k-2) \vee 0$,
$$n'_- \leq C \ t^{-1} \ a_+ \ I_{0}\left \{ y_{1_-} + a_+^2 \ \ell n 
\ t \ y_- + a_+ (z_{1_-} + a_+ \ I_0(y_-) + n_-) \right \} \ .
\eqno(5.73)$$

\noi In the same way as in the proof of Proposition 4.3, we continue 
the argument with a simplified version of the system (5.70)-(5.73)
where we exponentiate the diagonal terms and in particular $E_1$ to a 
constant according to (\ref{4.51e}) and where we eliminate the
constants and the factors containing $a_+$, with the consequence that 
we lose detailed control of the dependence of the lower bounds for
$T$ on $a_+$. Thus we rewrite (5.70)-(5.73) as (compare with 
(\ref{4.111e})-(\ref{4.114e}))
$$|\partial_t y'_-| \leq t^{-2} \left \{ z_{1_-} + \ell n\ t(I_0(y_-) 
+ n_- )\right \} + t^{-1-\beta} \ I_{k-2}(y_-) + t^{-1}\ n_- +
t^{-2}\ y'_{1_-} \ , \eqno(5.74)$$
$$|\partial_t y'_{1_-}| \leq t^{-2}\left \{ z_{2_-} + \ell n \ t 
(I_0(y_-) + n_-) \right \} + t^{-1} \ I_{k-2}(y_-) + t^{-1} \
n_- \ , \eqno(5.75)$$
$$|\partial_t z'_{j_-}| \leq t^{-2} \left \{ \ell n\ t\ z_{j_-} + 
\delta_{j_2} \ t^{(2- \alpha )\beta}\ z_{1_-}
\right \} + t^{-1} (y_- + \delta_{j_2} \ y_{1_-} ) + t^{-1 + j \beta} 
\ I_m (y_-)  \ , \eqno(5.76)$$
$$n'_- \leq t^{-1} \ I_{0}\left ( y_{1_-} + \ell n \ t \ y_- + 
z_{1_-} + I_0(y_-) + n_- \right ) \ .
\eqno(5.77)$$

\noi We now define
$$\left \{ \begin{array}{l} Y_- = \ \parallel t(\ell n\ t)^{-1} 
y_-;L^{\infty}([T, \infty))\parallel \  , \ Y_{1_-} = \ \parallel t(
\ell n\ t)^{-1} \ y_{1_-}; L^{\infty}([T, \infty))\parallel \ , \\ \\
Z_{j_-} = \ \parallel t^{1- j \beta} (\ell n\ t)^{-1} 
z_{j_-};L^{\infty}([T, \infty))\parallel \  , \ j = 1,2 \ ,\\ \\
N_{-} = \ \parallel t^{2- \beta} (\ell n\ t)^{-1} n_- ;L^{\infty}([T, 
\infty))\parallel \  \\ \end{array} \right . \eqno(5.78)$$

\noi and similarly for the primed quantities. It follows from 
(5.29)-(5.36) that all those quantities are finite. Using those
definitions and omitting the $-$ indices in the remaining part of the 
contraction proof, we obtain from (5.74)-(5.77)
$$|\partial_t y'| \leq t^{-2} \ell n \ t \left \{ Z_{1} \ t^{-1 + 
\beta} + Y\ t^{-\beta} + N\  t^{-1+\beta} + Y'_1 \ t^{-1} \right \}
  \ , \eqno(5.79)$$
$$|\partial_t y'_1| \leq t^{-2} \ell n \ t \left \{ Z_{2} \ t^{-1 + 
2\beta} + Y + N\  t^{-1+\beta} \right \}
  \ , \eqno(5.80)$$
$$|\partial_t z'_{j}| \leq t^{-2} \ell n \ t \left \{ Z_j \ t^{-1 + 
j\beta} \ \ell n \ t + \delta_{j_2} \ Z_1\ t^{-1+(3- \alpha
)\beta} + \delta_{j_2} \ Y_{1}  + Y\ t^{j \beta} \right \}  \ , \eqno(5.81)$$
$$n' \leq t^{-2} \ell n \ t \left \{ Y_1 + Y \  \ell n \ t + Z_1\ 
t^{\beta} + N\ t^{-1+ \beta} \right \}  \ . \eqno(5.82)$$

Integrating (5.79)-(5.81) between $t$ and $\infty$ with initial 
condition $(y',y'_1,z'_j)(\infty ) = 0$, and substituting the result 
and
(5.82) into the primed analog of the definition (5.78) (with the $-$ 
indices omitted), we obtain
$$\left \{ \begin{array}{l} Y' \leq  T^{-1+ \beta} \ Z_1 + 
T^{-\beta}\ Y + T^{-1+\beta} \ N + T^{-1} \ Y'_1 \ ,\\ \\
Y'_1 \leq  T^{-1+ 2\beta} \ Z_2 + Y + T^{-1+\beta}\ N  \ ,\\ \\
Z'_j \leq  T^{-1} \ell n \ T \ Z_j + \delta_{j_2} \ T^{-1-(\alpha 
-1)\beta}\ Z_1 + \delta_{j_2} \ T^{-2\beta} \ Y_1 + Y \ ,\\ \\
N' \leq  T^{-\beta} \ Y_1 + T^{-\beta} \ell n \ T \ Y + Z_1 + T^{-1}\ 
N \ .   \\ \end{array} \right . \eqno(5.83)$$

\noi Substituting $Y'_1$ from the second inequality into the first 
one, we recast the system (5.83) into the form (\ref{4.124e}), from
which the contraction property follows exactly as in the proof of 
Proposition 4.3. This proves that $\Gamma$ has a unique fixed point in
${\cal R}$.\par

The uniqueness of the solution in ${\cal C}(I, X^{k})$ under the 
conditions (5.60) (5.61) follow from Proposition 4.2, part (2). Note
that the time decay (5.60) (5.61) required for uniqueness is weaker 
than that contained in (5.29)-(5.32).\\

\noi \underbar{Part (2)}. Let $w_{+i}$, $i = 1,2$, satisfy the 
assumptions of the proposition with
$$|xw_{+i}|_{k + \alpha} \vee |w_{+i}|_{k+ \alpha + 1} \leq a_+ \eqno(5.84)$$

\noi and define $(W_i, S_i)$ by (5.10). Let $(w_i, s_i, B_{bi}) = 
(W_i + q_i, S_i + \sigma_i, B_{bi})$ be the two solutions of the
system (\ref{2.34e}) (\ref{2.35e}) obtained in Part (1). Let $(W_-, 
S_-) = (1/2) (W_1 - W_2, S_1 - S_2)$, define $(w_-, s_-, B_{b_-})$
as in Lemma 4.2 and define
$$\left \{ \begin{array}{l} y_- = |w_-|_{k-1} \vee |xw_-|_{k-1} \quad 
, \quad y_{1_-} = |w_-|_k \\ \\
  z_{j_-} = |s_-|_{k+j-1}^{\dot{}} \quad ,  \ j = 1,2, \quad n_- = 
|B_{b_-}|_{k}^{\dot{}} \vee  |x \cdot B_{b_-}|_{k}^{\dot{}} \ .\\
\end{array} \right . \eqno(5.85)$$

We assume that $w_{+1} - w_{+2}$ is small in the sense that
$$|x(w_{+1} - w_{+2})|_{k-1} \vee |w_{+1} - w_{+2}| \leq \eta \eqno(5.86)$$

\noi and we want to show that $(w_-, s_-, B_{b_-})$ is small by 
estimating $(y_-, y_{1_-}, z_{j_-}, n_-)$ in terms of $\eta$. We first
estimate $(W_-, S_-)$. From (\ref{4.140e}) it follows that
$$|xW_-|_{k-1} \vee |W_-|_k \leq \eta + t^{-1} \ a_+ \leq 2 \eta \eqno(5.87)$$

\noi for $t \geq a_+/\eta$. Furthermore, in the same way as in Lemma 4.2
$$|S_-|_{k+j-1}^{\dot{}} \leq C\ a_+\ \eta\ t^{j\beta} \qquad {\rm 
for}\ j = 1,2 \ . \eqno(5.88)$$

\noi We now take $t_0$ large in a sense to be specified below, and we 
estimate $(y_-, y_{1_-}, z_{j_-}, n_-)$ separately for $t \geq
t_0$ and for $t\leq t_0$. For $t \geq t_0$, using (5.87) (5.88) and 
(5.29)-(5.31), we obtain
$$\left \{ \begin{array}{l} y_- \vee y_{1_-} \leq  2 \eta + Y\ t^{-1} 
\ell n\ t \leq C\ \eta  \ ,\\ \\
z_{j_-} \leq  C\ a_+\ \eta \ t^{j\beta} + Z_j \ t^{-1} \ell n \ t 
(\ell n \ t + t^{(j-1)\beta}) \leq C\ a_+\ \eta \ t^{j\beta}   \\
\end{array} \right . \eqno(5.89)$$

\noi for $t_0$ large in the sense that
$$\eta\ t_0(\ell n\ t_0)^{-1} = C(a_+) \geq Y \vee \left ( (Z_1 \vee 
Z_2) a_+^{-1} \ t_0^{-\beta} \ell n\ t_0 \right ) \ .\eqno(5.90)$$

\noi (Remember that in this proposition, $Y$, $Z_1$, $Z_2$ are 
functions of $a_+$). Using the fact that for $t \geq t_0$, $I_m(f)$
depends only on $f$ restricted to $t \geq t_0$, we obtain in addition
$$n_- \leq C\ t^{-1} \left ( \eta (1 + \ell n \ t + a_+ \ t^{\beta} ) 
+ n_-\right )$$

\noi and therefore
$$n_- \leq C(1 + a_+) \eta\ t^{-1+\beta} \eqno(5.91)$$

\noi for $t \geq t_0$. \\

We next estimate $(y_-, y_{1_-}, z_{j_-}, n_-)$ for $t \leq t_0$ and 
for that purpose, we use the system (5.74)-(5.77) with the primes
omitted, since $(w_i, s_i, B_{bi})$ are solutions of the system 
(\ref{2.34e}) (\ref{2.35e}). We choose $\lambda$ such that $2\beta <
\lambda < 1$ (in the same way as in the proof of Proposition 4.2, 
part (2)) and we define
$$\left \{ \begin{array}{l} Y_- = \ \parallel t^{\lambda} 
y_-;L^{\infty}([T, t_0])\parallel \  , \ Y_{1_-} = \ \parallel 
t^{\lambda} \
y_{1_-}; L^{\infty}([T, t_0])\parallel \ , \\ \\  Z_{j_-} = \ 
\parallel t^{\lambda - j\beta} \ z_{j_-};L^{\infty}([T,
t_0])\parallel \  , \ j = 1,2 \ ,\\ \\ N_{-} = \ \parallel t^{\lambda 
+ 1- \beta} \  n_- ;L^{\infty}([T, t_0])\parallel \ .
\\ \end{array} \right . \eqno(5.92)$$

\noi Substituting those definitions into (5.74)-(5.77), using the fact that
$$I_m(y_-) \leq I_m \left (t^{-\lambda} \ Y_- + C \eta \right ) \leq C\left
( t^{-\lambda} \ Y_- + \eta \right ) \eqno(5.93)$$

\noi and similar relations for $y_{1_-}$, $z_{j_-}$, $n_-$, and 
omitting the $-$ indices, we obtain

$$\left \{ \begin{array}{l} |\partial_t y| \leq \eta\ t^{-1-\beta} + 
t^{-1-\lambda} \{ \cdot \} \\ \\
|\partial_t y_1| \leq \eta\ t^{-1} + t^{-1-\lambda} \{ \cdot \} \\ \\
|\partial_t z_j| \leq \eta\ t^{-1+j\beta} + t^{-1-\lambda} \{ \cdot \} \\ \\
n_- \leq \eta\ t^{-1+\beta} + t^{-1-\lambda} \{ \cdot \} \\ 
\end{array} \right . \eqno(5.94)$$
\noi where the brackets in the RHS are the same as in (5.79)-(5.82). 
Integrating the first three inequalities of (5.94) between $t$ and
$t_0$ for $t \leq t_0$ with initial condition at $t_0$ estimated by 
(5.89), and omitting again absolute constants, we obtain
$$\left \{ \begin{array}{l}  y \leq \eta + t^{-\lambda} \{ \cdot \} \\ \\
y_1 \leq \eta\ \ell n \ t_0 + t^{-\lambda} \{ \cdot \} \\ \\
z_j \leq \eta\ t_0^{j\beta} + t^{-\lambda} \{ \cdot \} \\ \\
n_- \leq \eta\ t^{-1+\beta} + t^{-1-\lambda} \{ \cdot \} \\ 
\end{array} \right . \eqno(5.95)$$

\noi where the brackets in the RHS are the same as in (5.94). We 
substitute (5.95) into the definitions (5.92) and obtain a system
similar to (5.83), with however the primes omitted, and with an 
additional term bounded by $\eta t_0^{\lambda} \ell nt_0$ in each of 
the
RHS. Proceeding therefrom as in the contraction proofs of Proposition 
4.3 and of Part (1) of this proposition, we obtain
$$X \leq \eta \ t_0^{\lambda}\ \ell n \ t_0 \eqno(5.96)$$

\noi where $X$ is defined by (\ref{4.125e}), so that by (5.90), $X$ 
tends to zero when $\eta$ tends to zero (actually as a power of
$\eta$). This proves the norm continuity of the map $w_+ \to (w, s, 
B_b)$ from the norm $|w_+|_k \vee |xw_+|_{k-1}$ (see (5.86)) to the
norm in $L^{\infty} (J, X^{k-1})$ for compact $J$. The last 
continuity follows from a standard compactness argument. \par\nobreak 
\hfill
$\sq$\par

\noi {\bf Remark 5.4.} In part (2) of Proposition 5.3, we prove 
actually a stronger continuity than stated, namely a suitably weighted
$L^{\infty}$ continuity in the whole interval $[T, \infty)$, as 
follows from (5.89) (5.91) for $t \geq t_0$ and (5.92) (5.96) for $t 
\leq
t_0$ defined in terms of $\eta$ by (5.90).

\mysection{Wave operators and asymptotics for (u, A)}
\hspace*{\parindent} In this section we complete the construction of 
the wave operators for the system (\ref{2.6e}) (\ref{2.7e}) in the
special case of vanishing asymptotic magnetic field, and we derive 
asymptotic properties of solutions in their range. The construction
relies in an essential way on Proposition 5.3. So far we have worked 
with the system (\ref{2.34e}) for $(w,s)$ and the first task is to
reconstruct the phase $\varphi$. Corresponding to $S$ defined by 
(\ref{2.40e}), we define
$$\phi = \int_1^t dt'\ t'^{-1}\left ( g(W) - (x \cdot B_a)_L (W) 
\right ) \eqno(6.1)$$

\noi so that $S = \nabla \phi$. Let now $(q, \sigma, B_b)$ be the 
solution of the system (\ref{2.49e}) (\ref{2.50e}) obtained in
Proposition 5.3 and let $(w, s) = (W + q, S + \sigma )$. We define 
$\psi$ by $\psi (\infty ) = 0$ and
$$\partial_t \psi = (2t^2)^{-1} \ |s|^2 + t^{-1} \left \{ g(w) - g(W) 
- (x \cdot B_a)_L (w) + (x \cdot B_a)_L (W) \right \}
\eqno(6.2)$$

\noi or equivalently
$$\psi = - \int_t^{\infty} dt' \left \{  (2t'^2)^{-1} \ |s(t')|^2 + 
t'^{-1} \ g(q, q + 2W) - t'^{-1}(x \cdot B_a)_L (q,q + 2W)  \right \}
\eqno(6.3)$$

\noi which is taylored to ensure that $\nabla \psi = \sigma$, given 
the fact that $S$ and $\sigma$ are gradients. The integral
converges in $\dot{H}^1$, as follows from (5.16) (5.18), from (5.27) 
(5.29) and from the estimate
$$\begin{array}{l}
\partial_t \parallel \sigma \parallel_2 \ \leq t^{-2} \parallel s 
\cdot \nabla s  \parallel_2 \ + t^{-1} \parallel \nabla g (q, q +
2W)\parallel_2 \\ \\
+ t^{-1} \parallel \nabla (x \cdot B_a)(q,q + 2W)\parallel_2 \\ \\
\leq t^{-2} \parallel s\parallel_{\infty} \ \parallel \nabla s
\parallel_2 \ + C \ t^{-1} a_+ ( \parallel q\parallel_2 + I_{-1} 
(\parallel xq \parallel_2 ) \\ \\
\leq C(a_+) \ t^{-2} (\ell n\ t)^2 \\ \end{array}
\eqno(6.4)$$

\noi so that
$$\parallel \nabla \psi \parallel_2 \ = \ \parallel \sigma 
\parallel_2 \ \leq C(a_+) \ t^{-1} (\ell n\ t)^2 \ . \eqno(6.5)$$

\noi Finally we define $\varphi = \phi + \psi$ so that $\nabla 
\varphi = s$. \par

We can now define the modified wave operators for the MS system in 
the form (\ref{2.6e}) (\ref{2.7e}) in the special case of vanishing
asymptotic magnetic field. We start from the asymptotic state $u_+$ 
for $u$ and we  define $w_+ = Fu_+$. The asymptotic state $(A_+,
\dot{A}_+)$ for $A$ is taken to be zero. We define $(W,S)$ by 
(\ref{2.40e}). We solve the system (\ref{2.49e}) (\ref{2.50e}) for 
$(q,
\sigma, B_b)$ by Proposition 5.3. Through (\ref{2.43e}), this yields 
a solution $(w, s, B_b)$ of the auxiliary system (\ref{2.34e})
(\ref{2.35e}). We reconstruct the phase $\varphi = \phi + \psi$ with 
$\phi$ and $\psi$ defined by (6.1) (6.3). We finally
substitute $(w, \varphi , B_b)$ into (\ref{2.17e}) and (\ref{2.18e}) 
with $B = B_a + B_b$ and $B_a$ defined by (\ref{2.28e}). This
yields a solution $(u, A)$ of the system (\ref{2.6e}) (\ref{2.7e}) 
defined for large time. The modified wave operator is the map
$\Omega  : u_+ \to (u,A)$ thereby obtained. \par

In order to state the regularity properties of $u$ that follow in a 
natural way from the
previous construction, we introduce appropriate function spaces. In 
addition to the operators
$M = M(t)$ and $D = D(t)$ defined by (\ref{2.14e}) (\ref{2.15e}), we 
introduce the operator
$$J = J(t) = x + it \ \nabla \ , \eqno(6.6)$$
\noi the generator of Galilei transformations. The operators $M$, 
$D$, $J$ satisfy the
commutation relation
$$i\ M\ D \ \nabla = J\ M\ D \ . \eqno(6.7)$$
\noi For any interval $I \subset [1, \infty )$ and any $k \geq 0$, we 
define the space
$${\cal X}^k(I) = \Big \{ u : D^*M^* u \in {\cal C} (I, H^{k+1}), 
D^*M^* x u \in {\cal C} (I, H^k) \Big \}$$
$$ = \Big \{ u:<J(t)>^{k+1} u \ {\rm and} \ <J(t)>^k xu \in {\cal 
C}(I,L^2)\Big \} \eqno(6.8)$$
 
\noi where $< \lambda > = (1 + \lambda^2)^{1/2}$ for any real number 
or self-adjoint operator
$\lambda$ and where the second equality follows from (6.7). \par

We now collect the information obtained for the solutions of the 
system (\ref{2.6e}) (\ref{2.7e}) in the range of the modified wave
operators and state the main result of this paper as follows.\\

\noi {\bf Proposition 6.1.} {\it Let $k > 3/2$, $0 < \beta < 1/2$ and 
let $\alpha > 1$ be such that $\beta (\alpha + 1) \geq 1$. Let $u_+$
be such that $w_+ = Fu_+ \in H^{k+ \alpha + 1}$ and $x w_+ \in 
H^{k+\alpha}$. Define $(W,S)$ by (\ref{2.40e}) and $a_+$ by (5.11). 
Then
\par

(1) There exists $T = T(a_+)$, $1 \leq T < \infty$, such that the 
auxiliary system (\ref{2.34e}) (\ref{2.35e}) has a unique solution
$(w,s,B_b) \in {\cal C}(I, X^{k})$ where $I = [T,\infty )$, satisfying
$$|w - W|_k \vee |x(w-W)|_k \leq C\ t^{-1} \ell n\ t  \eqno(6.9)$$
$$|w - W|_{k+1} \leq C\left ( t^{-1} \ell n\ t  + t^{-\alpha\beta} 
\right )  \eqno(6.10)$$
$$|s-S|_{k+j}  \leq C\ t^{-1} \ell n\ t \left ( \ell n \ t  + 
t^{j\beta} \right ) \quad {for}\ j = 0, 1, 2\ .
\eqno(6.11)$$
$$|B_b|_{k+1}^{\dot{}} \vee |x \cdot B_b|_{k+1}^{\dot{}} \leq C\ 
t^{-1}\ell n\ t \ .\eqno(6.12)$$
 
(2) Let $\phi$ and $\psi$ be defined by (6.1) and (6.3) with $q = w - 
W$, and let $\varphi = \phi + \psi$. Let
$$u = M D \exp (-i \varphi ) w \ , \eqno(2.17)\equiv (6.13)$$
$$A = t^{-1} \ D_0 B \eqno(2.18)\equiv (6.14)$$

\noi with $B = B_a + B_b$ and $B_a$ defined by (\ref{2.28e}). Then $u 
\in {\cal X}^k(I)$, $(A, \partial_t A) \in {\cal C}(I,K^{k+1}
\oplus H^k)$, $(u,A)$ solves the sytem (\ref{2.6e}) (\ref{2.7e}) and 
$u$ behaves asymptotically in time as $MD \exp (- i \phi ) w_+$ in
the sense that $u$ satisfies the following estimates.
$$ \parallel <J(t)>^k <|x|/t> ( \exp (i \phi (t, x/t)) u(t) - M(t) \ 
D(t) \ W(t)) \parallel_2 \ \leq C \ t^{-1} (\ell n\ t)^2 \ ,
\eqno(6.15)$$
$$ \parallel <J(t)>^{k+1} ( \exp (i \phi (t, x/t)) u(t) - M(t) \ D(t) 
\ W(t)) \parallel_2 \ \leq C \left ( t^{-1} (\ell n\ t)^2 +
t^{-\alpha \beta} \right ) \ , \eqno(6.16)$$
$$ \parallel <|x|/t> ( u(t) - M(t) \ D(t) \exp (-i \phi (t)\ 
W(t))\parallel_r \ \leq C \ t^{-1-\delta (r) } (\ell n\ t)^2
\eqno(6.17)$$

\noi for $2 \leq r \leq \infty$, with $\delta (r) = 3/2 - 3/r$. \par

Furthermore $A$ behaves asymptotically in time as $t^{-1}D_0 B_a(W)$ 
in the sense that the following estimates hold
$$|B - B_a(W)|_{k+1}^{\dot{}} \vee |\nabla x \cdot (B - 
B_a(W))|_{k}^{\dot{}} \ \leq C\ t^{-1}\ell n\ t \eqno(6.18)$$

\noi where $A$ and $B$ are related by (6.14).} \\

\noi {\bf Proof.} The proof follows from Proposition 5.3 supplemented 
with the reconstruction of $\varphi$ described above in this
section, except for the estimates (6.15)-(6.17) on $u$. In particular 
the estimates (6.9)-(6.12) are the estimates (5.29)-(5.32)
supplemented with (6.5), while (6.18) follows from (5.32) and (5.40).\par

We next prove the estimates (6.15)-(6.17) on $u$. From (6.13) with 
$\varphi = \phi + \psi$ and from (6.7), it follows that
$$ \parallel |J|^m ( \exp (i D_0 \phi ) u - M D W) \parallel_2 \ = \ 
\parallel \omega^m (\exp ( - i \psi ) w - W)
\parallel_2\ , \eqno(6.19)$$
$$ \parallel |J|^m (|x|/t) ( \exp (i D_0 \phi ) u - M D W) 
\parallel_2 \ = \ \parallel \omega^m |x| (\exp ( - i \psi ) w - W)
\parallel_2 \ .\eqno(6.20)$$

\noi We next estimate for $0 \leq m \leq k+1$
$$\parallel \omega^m (\exp ( - i \psi ) w - W)\parallel_2 \ \leq C 
\parallel \omega^m (\exp ( - i \psi ) -1)\parallel_{r_1} \ \parallel
w \parallel_{r_2}$$
$$+\  C \parallel \exp ( - i \psi )  - 1\parallel_{\infty} \ 
\parallel \omega^m  w \parallel_2 \ + \ C \parallel \omega^m (w -
W)\parallel_2$$
$$\leq C \exp (C \parallel \psi \parallel_{\infty}) \parallel\omega^m 
\psi\parallel_{r_1} \ \parallel w\parallel_{r_2} + C \left ( 
\parallel \psi
\parallel_{\infty}\ |w|_m + |w-W|_m \right ) \eqno(6.21)$$

\noi with $1/r_1 + 1/r_2 = 1/2$, $r_1 < \infty$, by Lemmas 3.2 and 
3.3, and similarly for $0 \leq m \leq k$

$$\parallel \omega^m |x| (\exp (-i \psi ) w - W)\parallel_2 \ \leq C 
\exp (C \parallel \psi \parallel_{\infty}) \parallel \omega^m
\psi \parallel_{r_1} \ \parallel xw \parallel_{r_2}$$
$$+ C \left ( \parallel \psi \parallel_{\infty}\ |xw|_m + |x(w-W)|_m 
\right ) \ . \eqno(6.22)$$

\noi Taking $r_1 = 6$, $r_2 = 3$ for $m = 0$ and $r_1 = 2$, $r_2 = 
\infty$ for $m \geq 1$, using the Sobolev inequality
$$\parallel \psi \parallel_{\infty} \ \leq C \left ( \parallel \sigma 
\parallel_{2}\ \parallel \nabla \sigma \parallel_{2}\right )^{1/2}
$$

\noi and using (6.9)-(6.11) yields (6.15)-(6.16). \par

The estimate (6.17) follows immediately from (6.15) and from the inequality
$$\begin{array}{ll}
\parallel f \parallel_r \ = t^{-\delta (r)} \parallel D^*M^* 
f\parallel_r \ &\leq C\ t^{-\delta (r)} \parallel \omega^{\delta (r)} 
D^*M^*f
\parallel_2\\ \\
&= C\ t^{-\delta (r)} \parallel |J(t)|^{\delta (r)} f \parallel_2 \\
\end{array}$$

\noi for $2 \leq r < \infty$ and from a similar inequality for $r = 
\infty$.\par\nobreak
\hfill $\sq$\par

\noi {\bf Remark 6.1.} The leading term in the asymptotic behaviour 
of $A$ is $t^{-1} D_0 B_a(W)$. Replacing $W$ by $w_+$ as a first
approximation, one obtains
$$A \sim  t^{-1} \ D_0 \ B_a(w_+)$$

\noi and since $B_a(w_+)$ is constant in time, that term spreads by 
dilation by $t$ and decays as $t^{-1}$ in $L^{\infty}$ norm. In the
norms considered in (6.18), that term is $O(1)$, so that the 
remainder is smaller than the leading term by $t^{-1} \ell nt$. We 
have
stated the remainder estimates in terms of $B$ rather than $A$ 
because they are simpler for $B$, since for $A$ the dilation $D_0$
induces a dependence of the time decay on the order of derivation. In 
fact (6.18) is equivalent to
$$\parallel \omega^m (A - t^{-1} D_0 B_a (W)) \parallel_2 \ \vee \ 
\parallel \omega^m \nabla x \cdot (A - t^{-1} D_0 B_a(W)) 
\parallel_2$$
$$\leq C\ t^{-m-1/2} \ \ell n \ t \eqno(6.23)$$

\noi for the relevant values of $m$, namely $1 \leq m \leq k+1$ for 
the first norm and $1 \leq m \leq k$ for the second one.\\

\noi {\large \bf Acknowledgments.} One of us (G.V.) is grateful to 
Professor D. Schiff for the hospitality extended to him at the
Laboratoire de Physique Th\'eorique in Orsay, where part of this work was done.

\end{document}